\documentclass[10pt]{article} 

\usepackage[margin=1in]{geometry}
\usepackage{fancyhdr}
\fancypagestyle{plain}{%
\fancyhf{} 
\fancyfoot[R]{-\,\fontsize{8pt}{8pt}\selectfont\thepage\,-} 

}
\RequirePackage{amsthm,amsmath,amsfonts,amssymb}
\RequirePackage[numbers]{natbib}
\RequirePackage[colorlinks,citecolor=blue,urlcolor=blue]{hyperref}
\RequirePackage{graphicx}

\theoremstyle{plain}

\newtheorem{theorem}{Theorem}[section]
\newtheorem{lemma}[theorem]{Lemma}
\newtheorem{corollary}[theorem]{Corollary}
\newtheorem{proposition}{Proposition}[section]
\theoremstyle{remark}
\newtheorem{definition}[theorem]{Definition}

\usepackage{amsmath}
\usepackage{amsfonts}
\usepackage{amssymb}

\usepackage{dsfont}

\usepackage{graphicx}
\usepackage{sidecap}
\usepackage{xspace} 
\usepackage{mathtools}
\usepackage{multirow}
\usepackage{hyperref}

\usepackage{thmtools,thm-restate}

\usepackage{letltxmacro}
\LetLtxMacro{\originaleqref}{\eqref}
\renewcommand{\eqref}{Eq.~\originaleqref}

\newcommand{\multidistro}{\text{mult}\xspace}
\newcommand{\dd}{\text d}
\newcommand{\equ}[1]{\begin{gather} #1 \end{gather}}
\newcommand{\sfrac}[2]{\mbox{$\frac{#1}{#2}$}}
\newcommand{\order}[1]{\mathcal O\left(#1\right)} 
\newcommand{\abss}[1]{\vert#1\vert} 
\newcommand{\quads}[1]{\quad #1 \quad}
\newcommand{\qand}{\quad \text{and} \quad}
\newcommand{\maxent}{\textsc{m}ax\textsc{e}nt\xspace}

\newcommand{\iproj}{$I$-projection\xspace}

\newcommand{\mle}{\textsc{mle}\xspace}

\newcommand{\idiv}{$I$-divergence\xspace}

\newcommand{\vspan}{\text{span}\xspace}

\newcommand{\prob}[1]{\mathfrak{#1}}
\newcommand{\refP}{\boldsymbol \upsilon} 

\makeatletter
\newcommand*\bigcdot{\mathpalette\bigcdot@{.7}}
\newcommand*\bigcdot@[2]{\mathbin{\vcenter{\hbox{\scalebox{#2}{$\m@th#1\bullet$}}}}}
\makeatother

\DeclarePairedDelimiterX\infdivx[2]{(}{)}{
  #1\;\delimsize\|\;#2
}
\newcommand{\infbraket}[2]{D\hspace{0.1mm}\infdivx{#2}{#1}}   
\DeclarePairedDelimiterX\braket[2]{\langle}{\rangle}{
  #1\;\vert\;#2
}
\DeclarePairedDelimiterX\braketO[3]{\langle}{\rangle}{
  #1\;\vert#2\vert\;#3
}
\newcommand{\ket}[1]{\vert#1\rangle}
\newcommand{\bra}[1]{\langle#1\vert}
\newcommand{\vev}[1]{\langle#1\rangle}

\newcommand{\infdiv}[2]{D\hspace{0.1mm}\infdivx{#1}{#2}}

\usepackage{xargs}

\newcommand{\config}{{\boldsymbol e}}
\newcommand{\statespace}{\mathcal E}
\newcommand{\statespaceR}{\statespace_{\star}}

\newcommand{\cartprod}{\boldsymbol{\mathcal D}}
\newcommand{\cartprodR}{{\boldsymbol{\mathcal D}_\star}}

\newcommandx{\equivclass}[3][1=\vecOp, 2=\prob f, 3=\statespace_0,usedefault]{[#2;{#1}]}
\newcommandx{\iprojSymbol}[1][1=\equivclass]{\mathbb P_{#1}}

\newcommand{\rank}{\text{rank}\,\xspace}
\newcommand{\Id}{\text{\small 1}\hspace{-3.5pt}\text{1}}

\newcommand{\op}[1]{\,\hat{#1}}
\newcommandx{\proj}[1][1=I]{\,\mathbb{#1}}
\newcommandx{\vecOp}[1][1=M]{\mathbf{\op#1}}

\newcommand{\totem}{\textsc{Tot}\textsc{Em}\xspace}
\newcommand{\totemplex}{\totem{plex}\xspace}

\usepackage{clipboard}
\usepackage[colorinlistoftodos,prependcaption,textsize=small]{todonotes}

\begin{document}


\thispagestyle{empty}

\begin{flushright}
\phantom{Version: \today}
\\
\end{flushright}
\vskip .2 cm
\subsection*{}
\begin{center}
{\Large {\bf Total Empiricism: Learning from Data
} }
\\[0pt]

\bigskip
\bigskip {\large
{\bf Orestis Loukas}\footnote{E-mail: orestis.loukas@uni-marburg.de},\,
{\bf Ho Ryun Chung}\footnote{E-mail: ho.chung@uni-marburg.de}\bigskip}\\[0pt]
\vspace{0.23cm}
{\it Institute for Medical Bioinformatics and Biostatistics\\
Philipps-Universität Marburg\\
Hans-Meerwein-Straße 6, 35032 Germany}

\bigskip
\end{center}

\begin{abstract}
Statistical analysis is an important tool to distinguish systematic from chance findings. Current statistical analyses rely on distributional assumptions reflecting the  structure of some underlying model, which if not met lead to problems in the analysis and interpretation of the results. Instead of trying to fix the model or ``correct'' the data, we here describe a totally empirical statistical approach that does not rely on ad hoc distributional assumptions in order to overcome many  problems in contemporary statistics. 
Starting from elementary combinatorics, we motivate  an information-guided formalism to quantify knowledge extracted from the given data.
Subsequently, we derive model-agnostic methods to identify patterns that are solely evidenced by the data based on our prior knowledge.
The data-centric character of empiricism allows for its universal applicability, particularly as sample size grows larger. 
In this comprehensive framework, we re-interpret and extend
model distributions, scores and statistical tests 
used in different schools of statistics.
\end{abstract}


\section{Introduction}

Statistical analysis helps distinguish ``true'' from chance findings. Statistics derives much of its power and usefulness from some distributional assumptions on the data. 
Model-centric learning essentially consists of verifying or slightly improving correctly suspected distributional assumptions  that accurately capture underlying mechanisms which generate data. In such convenient scenario, we are mostly aware (even if it is not readily acknowledged) of statistics that are sufficient to ideally describe the data, every deviation usually attributed to sampling noise~\cite{hullermeier2021aleatoric}. 
However, this attitude towards modeling does not exclude that our model might be plagued by both intentional and unintentional misconceptions~\cite{Lehmann1990ModelST} introducing tension with the reality of data.

As datasets grow larger and more versatile, they refuse to adhere to the distributional assumptions made by statisticians. 
The violation of distributional assumptions~\cite{10.1093/aje/kwx259}, which increasingly manifests  in any statistical analysis, points towards what is often described as model misspecification~\cite{10.3389/fevo.2019.00372}.
In fact, there exists a large body of literature~\cite{mayo2018statistical,https://doi.org/10.1111/cobi.13984} concerned about correcting the adverse effects of unmet distributional assumptions in order to rescue the statistical analysis based on them.
%
%
In this work, we take an alternative route. Instead of ``repairing'' the adverse effects of unmet distributional assumptions, we remove them altogether and focus on the observed data itself.
 Given a dataset or a series of datasets on some subjects,  the question then arises regarding ``good'' statistics that adequately summarize  the data~\cite{380f7170-a649-307c-9495-f3b3298846ff}. 
 
Abandoning distributional assumptions to answer this question paves the way for a novel approach to statistics based solely on data, which we refer to as \textsc{Tot}al \textsc{Em}piricism (\totem). Starting from the concept of the \emph{reality of data} gives rise to the notion that all attributes multivariately describing the studied subjects inevitably take a finite number of manifestations bounded by the sample size $N$. 
In particular, metric attributes would only take a finite number of manifestations, even if we believed that they descend from (in)finite interval on the real line.
Instead of thus augmenting observations with infinitely many unobserved values to approximate the distribution of an attribute by e.g.\ a normal distribution, we concentrate on the \textit{observed} manifestations.
The finite and countable number of manifestations of attributes --\,the so-called \textit{entities}\,-- in a dataset defines the resolution power of our empirical analysis. Any statement beyond entities could lead to daydreaming away from data.

Starting from the  entities always allows to compute their frequencies in a dataset, which we summarize by the empirical distribution $\prob f$.
Therefore, the subjects in the dataset form a closed system where all information about them lies within the empirical distribution assigning a weight to each entity. More generally, a linear operator $\hat X$ assigns $c$-numbers to entities. 
In this way, $\hat X$
extracts statistical information about the subjects through its expectation on the data (summarized by $\prob f$) 
which represents a collective characteristic. Depending on the application, collective characteristics can correspond to marginal distributions of subsets of attributes, moments of metric attributes or generalizations thereof. 
Statistics is generally concerned with demarcating those characteristics that are required to describe the subjects in the dataset from those collective characteristics that capture 
random noise in the dataset. Any demarcating statement inevitably refers back to our prior knowledge about the investigated subjects, which we encode by a distribution $\refP$.

To foster the aforementioned demarcation, we propose to consider probability distributions with the property that the collective characteristics extracted by some linear operators exactly match the corresponding expectations from the empirical distribution $\prob f$. Such \totem construction allows for an infinite set of probability distributions that all share with $\prob f$ the same expectations about the chosen characteristics. We refer to this set of probability distributions as the \totemplex. 
Subsequently, we select one, special probability distribution $\prob q$ from the \totem, which can be understood as the $I$-projection~\cite{Cover2006} of the reference $\refP$ onto the \totemplex.
%
Finally, we derive the probability density over distributions in the \totem, as the universal result of a large $N$ expansion around the $I$-projection. From this probability density over the \totem, we formulate the $I$-score and the $I$-test, which allow for demarcating those collective characteristics that are required to describe the subjects in a dataset from those that --\,at the current data size\,-- cannot be distinguished from sample noise.

All in all, we provide a novel approach for the statistical analysis of data that does not require distributional assumptions based on and generalizing the fruitful path laid out in~\cite{loukas2022categorical,loukas2022entropy}. 
In our data-centric treatment, we do not need to commit to any particular model or try to be near some (potentially unknown and unreachable) ground truth~\cite{kato2022view}. 
Instead of regarding continuous probability densities as fundamental, we 
derive them --\,at most\,-- as an approximation to the fundamentally discrete problem of combinatorics given $N$  data points. 
Mainly, our study is guided by empiricism~\cite{Wold1987-WOLTEA,5c1cde5d-facd-306d-9616-595d010829b9,popper2014conjectures} without excluding though, testing and  extending of established models. 
The universality of empiricism at larger sample sizes ensures that the derived methodology would self-consistently report on ``good'' observations, i.e.\ summary statistics, that are needed to enhance our prior knowledge about investigated subjects, no matter how misinformed this prior knowledge might be.

\paragraph{Outline of the paper.}
The paper is structured as follows: First, we introduce in Section~\ref{sc:TOTEM_algebra} an operator formalism to deal with probability distributions and summarize the phenomenology of data via expectations in a parameter-agnostic way. 
Subsequently, we develop in Section~\ref{sc:MIDIV} the totally empirical theory of sampling at large sample size $N$.
After variationally motivating the effective description of the large-$N$ combinatorics through a probability density over distributions that comply with phenomenology, we exemplify (Section~\ref{ssc:MIDIV}) the special role of the the $I$-projection in empirical studies. Taking advantage of this special distribution  enables to derive 
the universal $N$-leading contribution to the effective description of sampling in Section~\ref{ssc:DistroOverDistros}.
Using the introduced probability density in Section~\ref{sc:ClassificationOfKnowledge} helps classify any additional knowledge gained by the phenomenology of data from the viewpoint of prior knowledge.  
Eventually, we show in Section~\ref{sc:Examples} how to apply  the proposed methods to uncover trends in the phenomenology of familiar settings.
Appendix~\ref{app:Regularization} discusses covariant regularization in probability space and its (potential) issues.
Appendix~\ref{app:Derivations} provides lemmas and proves supporting the stated theorems. Appendix~\ref{app:ModelCentric} relates our findings to standard score and test in the model-centric literature. 

\section{The Algebra of Empiricism} 
\label{sc:TOTEM_algebra}

In a fully data-centric approach, the very concept of a model is critically  challenged, 
since observations can change not only our prior knowledge about the strength of suspected associations in the model, but the very model structure itself. 
Hence, we need a flexible formalism that would keep track of such structural variations.  
In this paragraph, we provide the algebraic characterization of a quantitative analysis that is solely based on empiricism.
%

\paragraph*{The entities}
An empirical description about some subjects 
comprises a finite set of $L$ attributes that distinctly (e.g.\ a person cannot be minor and adult at the same time) characterize these subjects. 
Each attribute $i \in \lbrace 1, \ldots, L \rbrace$ admits a finite number of countably many {categories} or levels; either unordered (e.g.\ black, white and yellow) or ordered 
(e.g.\ detector counts). In this spirit, a metric variable over a finite support like human height resolved up to centimeter scale can be thought of as an ordered categorical attribute.
The set of possible categories for an attribute constitutes its domain $\mathcal D_i$. Its specification is based on plausibility 
or prior knowledge about the investigated subjects. 

Any subject is represented by some configuration 
where each attribute assumes a definite category  from its respective domain, $e_i\in\mathcal D_i$ for all $i=1,\ldots,L$.
For ease of reference, we can think of a subject as an entity\footnote{Given a well-defined notion of a system, its entities would be called \textit{micro}states in statistical physics.} 
labeled by an $L$-tuple $\left(e_1,\ldots,e_L\right)$.
In more algebraic terms, an entity can be represented by an element
\equ{
\ket{\config} \equiv \ket{e_1,\ldots, e_L}
}
living in an abstract vector space. 
Two subjects represented by the same entity cannot be empirically distinguished. 
Hence, the maximal resolution power of empiricism after specifying attributes is dictated by the Cartesian product of their domains,
\equ{
\label{eq:CartesianProduct}
\cartprod = \mathcal D_1 \times\cdots\times \mathcal D_L = \left\lbrace \left(e_1,\ldots,e_L\right) ~\vert~ e_1\in \mathcal D_1,\ldots,e_L\in\mathcal D_L\right\rbrace~,
}
which indexes the basis of our vector space,
\equ{
\statespace = \left\{ \ket{e_1,\ldots,e_L} ~\vert~ \left(e_1,\ldots,e_L\right)\in \cartprod\right\}~.
}

This basis includes all entities 
both observed and unobserveded,  which are maximally anticipated to be distinguishable in the attribute-based approach.
In principle, an unobserved entity that belongs to a subset $\statespace_0\subset\statespace$ could be attributed either to finite-sampling effects on rare realizations (e.g.\ a desert receiving more than 25 centimeters of precipitation a year) or deterministically impossible situations (e.g.\ being male and pregnant). In statistics literature~\cite{bishop2007discrete}, former entity is referred to as a sampling zero and latter as a structural zero. 
We shall call both types \textit{null}entities. 

Let $\cartprodR$ denote manifestations of the Cartesian product \eqref{eq:CartesianProduct} that retain some stochasticity.
A sampled record 
from the database would be always a distinct entity $\ket{\config}\in\statespaceR$ where
\equ{
\statespaceR\equiv\statespace\setminus\statespace_0=\{\ket{\config}\in\statespace~\vert~\config\in\cartprodR \}
}
signifies the subspace spanned by admissible entities. 
After (arbitrarily) enumerating possible entities via ${enum}:\statespaceR\rightarrow\mathbb N$, one can intuitively identify them 
with the canonical basis in $\mathbb R^{\abss{\statespace}-\abss{\statespace_0}}$:
\equ{
\label{eq:BasisVectors:ColumnVectorRepresentation}
\ket{\config}
\overset{\cdot}{=} 
\begin{pmatrix}
0 \\[-1.ex]
\vdots\\[-0.6ex]
1      \\[-0.8ex]
\vdots\\[-0.5ex]
0
\end{pmatrix}
\quad\leftarrow 
 \text{ row } {enum}(\config)
~.
}
In that way, entities can be thought of as column vectors in $\mathbb R^{\abss{\statespace}-\abss{\statespace_0}}$.
Dual to entities we introduce row vectors $\bra{e}$ so that 
we can define a scalar product with
\equ{
\label{eq:OrthonormalBrakets}
\braket{\config'}{\config} = \prod_{i=1}^L \delta(e_i,e'_i) \quads{\text{for}} \ket{\config}\in\statespace\qand\bra{\config'}\in \statespace^*
}
where $\delta$ denotes $L$ indicator 
maps 
$\mathcal D_i\times\mathcal D_i\rightarrow\{0,1\}$ in the respective feature domains.

\paragraph*{The data}
In any setting, there exists (even if we do not have direct access to)
a data matrix $\texttt X$ 
with the records from a sample of size $N$ in the rows and $L$ features along the columns. In the picture of \eqref{eq:BasisVectors:ColumnVectorRepresentation}, the $r$-th row of $\texttt X$ corresponds to an entity from $\statespace$. 
%
The total number of records corresponding to rows of $\texttt X$ defines the only (hopefully large) scale we are going to use in our analysis, the sample size $N$. 
From the data matrix, we can further calculate over all distinct configurations in \eqref{eq:CartesianProduct}  relative frequencies 
\equ{
\label{eq:Empirical:RelativeFrequencies}
f_\config 
= \frac{1}{N} \sum_{r=1}^N \delta(e_1, \texttt X_{r,1}) \cdots \delta(e_L, \texttt X_{r,L})
\quads{\text{for}} \config\in\mathcal D_1\times\dots\times\mathcal D_L~.
}

\paragraph*{Probabilities}
We quantify our algebraic characterization of inference problems by only incorporating physically invariant, i.e.\ parameter-agnostic information.
As explained in the previous paragraph, 
an empirical analysis can at most draw conclusions on subjects from 
the abundance of entities, so that (relative) frequencies in \eqref{eq:Empirical:RelativeFrequencies} become the cornerstone of data-centric studies.
Algebraically, {our} quantitative {perception} of distinct entities  can be fully described by a diagonal operator 
\equ{
\label{eq:JointDistribution}
\prob p = \sum_{\config\in\cartprodR} 
p_\config\,\ket{\config}\bra{\config}~
\in \,\mathcal{B}(\statespaceR)
~,
}
acting on $\text{span}(\statespaceR)$, the space spanned by 
the set of relevant entities.
By construction, its eigenvalues 
\equ{
\prob p\ket{\config} = p_{\config} \ket{\config}
\quads{\text{for}} \ket{\config}\in\statespaceR
}
give the expansion coefficients $p_\config\geq0$ which 
weigh our anticipation about each admissible entity.
Such operators are sufficient to characterize our knowledge about subjects for any practical purpose.

At given sample size $N$, the smallest possible variation of a relative frequency
\equ{
\label{eq:Diophantine:CountVariation}
p_\config \quads{\rightarrow} p_\config + \delta p_\config
}
is naturally of order $1/N$ corresponding to a change of one data entry. 
Statistics as an effective description is expected to become sensible only at sufficiently large $N$.
In anticipation of the large-$N$ expansion motivated in Section~\ref{ssc:DicreteContinuous}, we can thus approximate relative frequencies with non-negative real numbers, the error being variationally suppressed by $N$.
Working at fixed sample size immediately leads to normalization condition
\equ{
\label{eq:Probabilities:Normalization}
\sum_{\config\in\cartprodR} p_\config 
\overset{!}{=} 1
~.
}
Normalization and non-negativity imply in turn, $p_\config\in[0,1]$.

Consequently, the expansion coefficients of $\prob p$ can be regarded at large $N$ as probabilities making $\prob p$ a probability operator.
In statistics, 
such map from $\statespaceR$ to  $[0,1]^{\abss{\statespace}-\abss{\statespace_0}}$ which assigns to each admissible entity $\ket{\config}$ a weight $ p_\config$
gives rise to the joint probability distribution on the simplex
\equ{
\label{eq:ProbabilitySimplex}
\mathcal P(\statespaceR) 
= \left\lbrace 
\prob p
~\big\vert~   
 p_\config \in\mathbb R^+_0 ~\text{ for }~ \config\in\cartprodR
\qand
\sum_{\config\in\cartprodR} p_\config=1
\right\rbrace
~.
}
For the simplex over $\statespace$, we simply write $\mathcal P\equiv\mathcal P(\statespace)$.

\paragraph*{The empirical distribution}
At the center of our algebra, there always exists the empirical distribution ${\prob f}$ over the full $\statespace$  with components given in \eqref{eq:Empirical:RelativeFrequencies} that can be easily extracted from $\prob f$ via 
\equ{
f_e = \braketO{\config}{\prob f}{\config} 
}
using the orthonormality of \eqref{eq:OrthonormalBrakets}.
The empirical set  of unobserved entities
\equ{
\label{eq:EmpiricalStatistics:nullentities}
\left\{ 
\ket{\config}\in\statespace~\vert~ f_\config = 0
\right\}
}
dictates  the maximal set of nullentities that is logically allowed by the provided data. 
Evidently,  nullentities from a declared $\statespace_0$ must be compatible with the empirical perception, meaning  that they form a subset of the empirical nullentities, 
\equ{
\label{eq:Empirical:nullentitiesCompatibility}
\forall\,\ket{\config}\in\statespace_0\,:
\quad
f_\config = 0
~.
}
Consequently, ${\prob f}$ can be accommodated within $\mathcal B(\statespaceR)$  
by appropriately tuning $p_\config = f_\config\geq0$.

\paragraph*{Priors and reference distribution}
In addition to the provided dataset $\texttt X$, we can express any prior knowledge on the subjects under investigation through a reference distribution $\refP\in\mathcal P(\statespace)$.
For example, ${\refP}$ could incorporate deterministic statements reflecting  structural zeros or more generally fundamental laws governing the subjects under investigation. 
Falsely declaring structural zeros 
that are later observed in the data would lead to infinities in the information-theoretic metrics. 
Hence, 
\begin{definition}
\label{def:BaseKnowledge}
    A data-compatible reference distribution ${\refP}$ fulfills 
\equ{
\nonumber
f_\config > 0 \quads{\Rightarrow}
\upsilon_\config > 0 
~,
}
i.e.\ its nullentities are a subset of the empirically observed nullentities.
\end{definition}
\noindent
Throughout the statistical analysis 
in Sections~\ref{sc:MIDIV} and~\ref{sc:ClassificationOfKnowledge}, 
${\refP}$ remains the reference state of knowledge to benchmark empirical descriptions and select the most adequate for a given problem.

\paragraph*{Minimal prior knowledge} 
A very symmetric situation arises when we are most ignorant about subjects and their relations.
Lacking any prior knowledge about the underlying mechanisms that generate data,
each attribute domain $\mathcal D_i$ enjoys a permutation symmetry (mere re-ordering) of its levels $\config_i$ as discussed e.g.\ in~\cite{keynes1921treatise}. 
At first, such symmetry might seem counter-intuitive for counts and metric variables where a natural order exists. However, such physical order is later re-instated via moment observations and must not be confused with the arbitrary order of the labels $e_i$.

Since the set of entities  inherits the aforementioned permutation symmetry, minimal knowledge about $\statespace$ is necessarily described by 
uniformly assigning probabilities to all
entities 
via
\equ{
\label{eq:UniformDistro:DEF}
{\prob u} = 
\frac{1}{\abss{\statespace}} \proj_{\statespace}
= \frac{1}{\abss{\statespace}}\sum_{\config\in\cartprod}\ket{\config}\bra{\config}
~.
}
The most symmetric probability assignment 
reflects the same degree of ignorance for all entities, as any other distribution would discriminate against some entity, hence breaking  permutation symmetry.

\paragraph*{Characterizing entities}
To quantitatively describe subjects, the associated entities must be related 
to some observable quantity that characterizes them.
At the algebraic level, 
a sensible characteristic
can be inferred 
using a suitable operator $\op{M}\in\mathcal B(\statespaceR)$ which assigns a definite value $m(\config)\in\mathbb R$ to each admissible entity
via a well-defined eigenvalue equation 
\equ{
\label{eq:Pheno:Measurement}
\forall\,\ket{\config}\in\statespaceR:\quad\op{M} \ket{\config} = m(\config)\ket{\config}
~.
}
Colloquially, we also refer to $\op{M}$ as the characteristic. 
Notice that linear action \eqref{eq:Pheno:Measurement} does not need to uniquely determine an entity.  
Most notably the resolution of identity 
\equ{
\label{eq:IdentityOperator}
\proj_{\statespaceR} = \sum_{\config\in\cartprodR}\ket{\config}\bra{\config}
}
assigns the same value to all admissible entities.
In contrast, projector 
operator
\equ{
\label{eq:ProjectorOperator}
\proj_{\config} = \ket{\config}\bra{\config}
}
fully adheres to the unique characteristic of entity $\ket{\config}$, which is unambiguously defined  by specifying the levels $e_i$ for all attributes, as already discussed. 

Most interestingly,
characteristics which non-trivially group subjects under common eigenvalues enable a quantitative  investigation of the induced collective effect. 
When a distribution $\prob p$ describes the knowledge about entities,
the statistics of $\op M$ --\,which collectively characterize subjects\,-- 
can be extracted as the $\prob p$-expectation 
\equ{
\label{eq:Expectation_DEF}
\vev{\prob p\op M} \equiv\braketO{\statespaceR}{\prob p \op M}{\statespaceR} = \sum_{\config\in\cartprodR}  p_\config\,m(\config) = \sum_{\config\in\cartprodR}\braketO{\config}{\prob p \op M}{\config}
}
where 
\equ{
\ket{\statespaceR} = \sum_{\config\in\cartprodR} \ket{\config} ~\in~ \text{span}(\statespaceR)~.
}
Whenever clear from context or if $\statespace_0=\emptyset$, we shall drop $\ket{\statespaceR}$ from expectations and use the compact notation $\vev{\prob p\op M}$.
Trivially, 
normalization condition~\eqref{eq:Probabilities:Normalization} can be written  as the $\prob p$-expectation of identity \eqref{eq:IdentityOperator}
\equ{
\vev{\prob p\proj_\statespace} = \vev{\prob p} \overset{!}{=} 1
~,
}
while the probability of entity $\ket{\config}$ may be understood 
as the $\prob p$-expectation of the corresponding  projector \eqref{eq:ProjectorOperator}
\equ{
\vev{\prob p \proj_\config} = \braketO{\config}{\prob p}{\config} = p_\config~.
}

In principle, there exist various --\,possibly unaccountably many, not necessarily distinct\,-- ways to characterize entities. This results in a multitude of statistics.   
Collective characterizations
like the average weight in a group of pupils or the prevalence of single-author papers in a field,
usually offer a concise  way to succinctly describe trends and patterns among attributes 
without committing to a certain entity.
Any quantity that is easily observed like net magnetization (against orientation of constituent spins) or easily accessed like prevalence of a disease (against patients' microdata)
could sensibly serve as a collective characterization of entities.   

Different types of characteristics can be analyzed contingent on the form of $m(\config)$ in \eqref{eq:Pheno:Measurement}. 
For instance, the eigenvalue problem 
\equ{
\label{eq:MarginalConstraintFkt}
\proj_{z_i} \ket{\config} = \delta(e_i, z_i) \ket{\config}
    \quads{,}z_i\in\mathcal D_i
}
would simply result in the marginal probability $\vev{\prob p\proj_{z_i}}$ of feature $i$ being on state $e_i = z_i$.
Also, moment measurements are easily realized, e.g. the $n$-th moment of feature $j$ 
can be extracted via 
\equ{
    \op L_{j}^n\ket{\config}
    = 
    e_j^n  \ket{\config}
    ~, 
}
invoking an operator which gets the level of desired attribute, $\op L_{j}\ket{\config} = e_j  \ket{\config}$.
Any 
combination of marginal and moment operators is evidently admissible.
%


Considering multiple characteristics $\op M_\alpha$ provides a summary of entities.
To obtain the summary statistics, i.e.\ the statistics of each
characteristic, 
we use:
\begin{definition}[constructing elements]
\label{def:Irreducible_VectorOperator}
    An irreducible constructing element (also thought as a vector operator)
    \equ{
        \vecOp \equiv \vspan \left\{
        \op M_1, 
        \ldots,
        \op M_D
        \right\}
    }
    is a collection of $D$ diagonalizable 
    operators $\op{M}_\alpha$ defined by linearly independent actions on $\text{span}(\statespaceR)$, 
    i.e.\ fulfill as an operator identity:
    \equ{
    \label{eq:OperatorIndependence}
         \sum_{\alpha=1}^Dc_\alpha\op{M}_\alpha = 0 \quads{\Rightarrow} c_\alpha = 0~\,\forall \,\alpha=1,\ldots,D
         ~.
    }
\end{definition}
\noindent
Within maximally diagonalizable operators in $\mathcal B(\statespaceR)$, 
we express 
through 
irreducible vector operator $\vecOp[N]$ satisfying 
\equ{
\label{eq:VectorOperator:OrthogonalComplement}
\vev{\op N_a \op M_\alpha} = 0 \quad\forall~\alpha=1,\ldots,D\qand a=1,\ldots \abss{\statespace}-\abss{\statespace_0}-D
}
the orthogonal complement to the space spanned by $\vecOp$.
In other words, $\vecOp[N]$ spans the kernel of $\vecOp$, $\ker\vecOp=\vspan(\vecOp[N])$.  
By definition \eqref{eq:OperatorIndependence}, the eigenvalues of the constructing element on 
admissible entities,
\equ{
\alpha=1,\ldots D \qand \ket{\config}\in\statespaceR\,:\quad\op{M}_\alpha\ket{\config} = m_\alpha(\config)\ket{\config}
~,
}
can be viewed as elements of a $D\times\left(\abss{\statespace}-\abss{\statespace_0}\right)$ matrix of {full row-rank}.
In that way, the $\prob p$-expectation of $\vecOp$ induces a  linear map 
\equ{
\label{eq:pheno:LinearMap}
\braketO{\statespaceR}{\prob p\vecOp}{\statespaceR}: ~~[0,1]^{\abss{\statespace}-\abss{\statespace_0}} \rightarrow \mathbb R^D
}
from the probabilities 
of distinct entities to 
the observable world of collective characteristics.

Of course, nothing prevents us from learning the data $\texttt X$ ``by heart'', when $\vecOp$ is equivalent to a collection of 
$\abss{\statespace}-\abss{\statespace_0}$ projector operators \eqref{eq:ProjectorOperator},
\equ{
\label{eq:LearningByHeart}
\Id 
=
\left\{\ket{\config}\bra{\config}\right\}_{\config\in\cartprodR}~,
}
represented by identity matrix 
in $\mathbb R^{\abss{\statespace}-\abss{\statespace_0}}$ (viz.\ identification of \eqref{eq:BasisVectors:ColumnVectorRepresentation}).
In that way, we recover the relative frequencies $f_\config$ from $\texttt X$.
Committing to empirical distribution ${\prob f}$, 
 trace action \eqref{eq:Expectation_DEF} can be understood as an abstraction for {measurement} (e.g.\ in the lab) or {observation} (e.g.\ in astronomy). 
 As both are formally realized by the same algebraic operation as an $\prob f$-expectation of characteristics, 
we interchangeably use the notions of measurement and observation. 
Due to \eqref{eq:Empirical:nullentitiesCompatibility}, 
the empirical statistics $\vev{\prob f\op M_\alpha}\equiv\braketO{\statespace}{\prob f\op M_\alpha}{\statespace} = \braketO{\statespaceR}{\prob f\op M_\alpha}{\statespaceR}$ of any operator $\op M_\alpha\in\mathcal B(\statespaceR)$ remain insensitive to nullentities. 

\begin{definition}[phenomenology]
\label{def:Phenomenology}
    Using irreducible 
    vector operator $\vecOp[M]$
    to infer summary statistics, 
    the set of 
    $D$ empirical expectations $\vev{\prob f\vecOp}$ that are  calculated from data
    constitutes a {phenomenology}.
\end{definition}

Defined as a set of observed values referring back to some distinct entities, phenomenology naturally induces an equivalence orbit on $\mathcal P(\statespaceR)$, 
\equ{
{\prob p} \quads{\overset{\vecOp}{\sim}} {\prob f} 
\nonumber
}
that contains all distributions with expectations
$\braketO{\statespaceR}{\prob p\vecOp}{\statespaceR}$
dictated by the empirical estimates 
$\vev{\prob f\vecOp}$.
Hence, phenomenology becomes a set of \textit{guidelines} to select adequate distributions based on the empirical summary statistics, 
%
such as the macroscopic effects of pressure and temperature hint towards the kinematic distribution of molecules in an atomic gas.
In most interesting applications, linear map \eqref{eq:pheno:LinearMap} --\,being usually rank-deficit\,-- is not one-to-one, so that 
phenomenology $\vev{\prob f\vecOp}$ cannot  fix the probabilities of all entities. 
For example, there could exist many atomic gases with different microscopic properties all exhibiting the same pressure at the same temperature.

Ultimately, we arrive at
\begin{definition}[\totem]
\label{def:TOTEM}
After declaring nullentities in $\statespace_0$ compatible with prior knowledge,
\equ{
\ket{\config}\in\statespaceR\,:\quad
\refP_\config > 0
~,
}
a \totem of the data  
is defined by 
selecting  
distributions $\prob p\in\mathcal P(\statespaceR)$ 
which support a phenomenology $\vev{\prob f\vecOp}$: 
\equ{
\label{eq:Pheno:LinearSystem}
\braketO{\statespaceR}{\prob p\vecOp}{\statespaceR}
\overset{!}{=} 
\vev{\prob f\vecOp}
~.
}
The induced equivalence class $\equivclass$ 
which comprises all distributions on $\mathcal P(\statespaceR)$ that satisfy  
the phenomenology 
of the \totem description 
forms the \totemplex. 
\end{definition}
\noindent
Two irreducible vector operators 
related by full-rank\footnote{$\mathbf T$ must be a $D\times D$  matrix  of full rank, otherwise the transformed vector operator would be reducible.} transformation matrix $\mathbf T$ via
\equ{
\label{eq:TOTEMs:Reparametrization}
\op{M}_\alpha\quads{\rightarrow}\sum_{\beta=1}^D T_{\alpha\beta}\, \op{M}_\beta \quads{\text{for}}\alpha=1,\ldots,D
}
construct the same \totemplex being for all practical purposes (\textsc{fapp}) equivalent.

Normalization condition \eqref{eq:Probabilities:Normalization} implemented by $\proj_{\statespaceR}$ 
is automatically taken  to be part of any $\vecOp$ constructing a \totem, either explicitly or implicitly e.g.\ by fixing two complementary marginals which empirically sum to one. 
In absence of observations besides the specification of $\statespace_0$, 
the \totemplex thus reduces to
the simplex \eqref{eq:ProbabilitySimplex}, $\equivclass[{\proj_{\statespaceR}}]=\mathcal P(\statespaceR)$
where joint distributions 
are still empirically indistinguishable points. 
In the opposite limit, whenever $\abss{\statespace}-\abss{\statespace_0}$ independent observations are specified so that  $\vecOp$ is \textsc{fapp} equivalent to \eqref{eq:LearningByHeart},
we have $\equivclass[\Id][\prob f] = \{\prob f\}$. 

Trivially, phenomenological constraints \eqref{eq:Pheno:LinearSystem} always admit at least one solution, the empirical distribution $\prob f$ itself.
As long as $D < \abss{\statespace}-\abss{\statespace_0}$, so that $\abss{\ker\vecOp}=\abss{\statespace}-\abss{\statespace_0}-D>0$, 
there exist infinitely many distributions on the \totemplex
generated by convexity through
\equ{
\label{eq:TOTEMplex:Convexity}
\lambda{\prob p} + \left(1-\lambda\right) {\prob p'}\in\equivclass
\quads{\text{for}}
{\prob p}, {\prob p'}\in\equivclass\subseteq\mathcal P(\statespaceR)
\qand\lambda\in[0,1]~,
}
where ${\prob p} - {\prob p'}\in\ker\vecOp$. 
In Section~\ref{ssc:DistroOverDistros}, we come back to the study of such kernel fluctuations. 

Using linear 
action \eqref{eq:Pheno:Measurement} of operators on admissible entities,
the defining operation \eqref{eq:Pheno:LinearSystem} of a \totem can be written in the spirit of \eqref{eq:pheno:LinearMap} as a linear program
\equ{
\label{eq:TOTEM:LinearProgram}
p_\config\geq 0 ~\forall\,\config\in\cartprodR
\qand
\sum_{\config\in\cartprodR} m_\alpha(\config) 
\left(p_\config - f_\config\right)
\overset{!}{=} 0
\quad\forall\, \alpha=1,\ldots D
}
in the non-negative probabilities $ p_\config$ 
with the components of coefficient matrix given by the eigenvalues 
$m_\alpha(\config)$.
In particular, 
full-rank linear transformation \eqref{eq:TOTEMs:Reparametrization} is
realized by  elementary row operations on the matrix of eigenvalues. 
Relating \textsc{fapp} equivalent constructing elements,
such transformations do not change the solution space of linear program \eqref{eq:TOTEM:LinearProgram}.
To robustly classify inequivalent \totem descriptions
the reduced row-echelon form $\mathbf R$
of the coefficient matrix in \eqref{eq:TOTEM:LinearProgram} could be invoked,
leading to generalized expectations. 
In particular, normalization on the simplex~\eqref{eq:ProbabilitySimplex} implies
\equ{
\sum_{\alpha=1}^{D} R_{\alpha,\config} = 1 ~~\forall~ \config\in\cartprodR 
~,
}
since there are by construction $D=\rank\mathbf \vecOp$ columns in $\mathbf R$ with all zeros but the leading one.

\subsection{Information-theoretic metrics}
\label{ssc:InformationTheoreticMetrics}
 
To algebraically compare different 
states of knowledge about entities 
we need a notion of distance. 
On  simplex $\mathcal P(\statespaceR)$,
an information-compatible, non-commutative form is given by the negative $\prob p$-expectation of the logarithm of probability operator $\prob q$
\equ{
\label{eq:LikelihoodForm:DEF}
\ell(\prob p; \prob q)\equiv-\vev{\prob p\log \prob q} = -\braketO{\statespaceR}{\prob p\log \prob q}{\statespaceR} = 
-\sum_{\config\in\cartprodR} p_\config\log q_\config 
~.
}
Evidently, it is linear in its first argument. In particular, 
the series definition of logarithm implies due to orthonormality \eqref{eq:OrthonormalBrakets} the expression  
\begin{align*}
\log \prob q = &\,\,
\log \left(\proj_{\statespaceR}+\sum_{\config\in\cartprodR}( q_\config-1)\ket{\config}\bra{\config}\right) 
=
\sum_{\config\in\cartprodR}(q_\config-1)\ket{\config}\bra{\config}
- \sfrac12\left(\sum_{\config\in\cartprodR}(\prob q_\config-1)\ket{\config}\bra{\config}\right)^2
+\ldots
\\
=&\,\,
\sum_{\config\in\cartprodR} \left[(q_\config-1)
- \sfrac12\left(\prob q_e-1\right)^2 + \ldots \right]\ket{\config}\bra{\config}
=
\sum_{\config\in\cartprodR} \log q_\config \ket{\config}\bra{\config}
~,
\end{align*}
where the resolution of identity \eqref{eq:IdentityOperator} has been used in rewriting the argument within the logarithm.
The induced norm squared 
\equ{
\label{eq:Norm:Entropy}
 H[\prob p] 
 \equiv 
 \ell(\prob p;\prob p)
 =
 -\sum_{\config\in\cartprodR} p_\config\log p_\config
%
}
coincides~\cite{6773024} with Shannon's {differential entropy} and is semi-positive definite. 

Notice that 
\equ{
\label{eq:likelihood}
\exp\left(-\braketO{\statespaceR}{\prob p\log \prob q}{\statespaceR}\right)
}
determines the \textit{likelihood} of ${\prob p}$ under ${\prob q}$. 
Any entity with  $ p_\config=0$ 
does not contribute to the $\log$-likelihood leaving thus the likelihood invariant. 
On the other hand, an entity with $ p_\config \neq 0$ for which $ q_\config = 0$ would result into vanishing likelihood  
as contradictory to logic. 
In particular, a compatible in the sense of~\ref{def:BaseKnowledge} reference distribution 
asserts that information-theoretic measures like \eqref{eq:LikelihoodForm:DEF} 
remain finite throughout our explorations. 
A safe strategy to avoid breakdowns of the information-theoretic description is to select a reference distribution which assigns\footnote{The opposite direction that $\refP_\config > 0$ implies $\prob f_\config > 0$ does not need  to hold; and will not for many realistic datasets.} non-zero probability to all possible entities, but (at most) structural zeros.

Based on the $\log$-likelihood, we now introduce    
a ``distance'' via 
\equ{
\label{eq:I_Divergence:DEF}
\infdiv{\prob p}{\prob q} :=  
-\ell(\prob p;\prob p) + \ell(\prob p;\prob q)
= \sum_{\config\in\cartprodR} p_\config\log\frac{ p_\config}{ q_\config}
}
to measure~\cite{amari2000methods}  the information divergence (in short \idiv) between two distributions on $\mathcal P(\statespaceR)$. Noting that $0\cdot\log0=0$, the \idiv can be lifted on the ambient simplex $\mathcal P$
to become  
the well-known~\cite{KL_divergence_originalPaper,kullback1997information} Kullback-Leibler divergence or relative entropy. 
This information-theoretic metric remains~\cite{Cover2006} non-negative 
and vanishes iff 
${\prob p} = {\prob q}$.
On the other hand, it is neither symmetric nor linear in its arguments.


\section{Knowledge given phenomenology}
\label{sc:MIDIV}

In this Section, we systematize the study of trends and patterns 
given a \totem~\ref{def:TOTEM} of the data triggered by phenomenology~\ref{def:Phenomenology}.  
%
Starting from a sensible reference distribution, we show the existence and uniqueness of a special distribution compatible with the phenomenology. 
Subsequently, we show how to estimate the probabilities of distributions on the \totemplex as a large-$N$ expansion around this special distribution. 


\subsection{From discrete to continuous}
\label{ssc:DicreteContinuous}

So far, we have not discussed the total number of entries in $\texttt X$, the sample size $N>0$.
In fact, $N$ plays a crucial role in any \totem, as it becomes the controlling parameter of the effective theory that assigns probabilities to distributions on the \totemplex. 
To appreciate the significance of this scale we emphasize the discrete nature of  \eqref{eq:TOTEM:LinearProgram}.
The actual empirical problem at hand consists of the Diophantine system of linear equations
\equ{
\label{eq:Diophantine:LinearSystem}
\sum_{\config\in\cartprodR} m_\alpha(\config) \left(N p_\config-N f_\config\right) = 0
\quads{\text{with}} 
N f_\config\,,\, N p_\config  \in \mathbb N_0^+ \quad\forall\,\config\in\cartprodR
}
to determine appropriate counts $N p_\config$ 
satisfying the phenomenology $\vev{\prob f\vecOp}$. 
For the purposes of this paragraph, we recycle operator notation $\prob p$ to represent estimates about relative frequencies. 
The solution set of \eqref{eq:Diophantine:LinearSystem} is an integral lattice polytope, i.e.\ a convex bounded polyhedron with vertices over the integer lattice.

Notice that the reference distribution itself does not need to belong to the polytope (or the \totemplex), since it would generically neither share the same nullentities with the \totem nor adhere to the defining expectations of phenomenology~\ref{def:Phenomenology}.
Nevertheless, $\refP\in\mathcal P$ as part of prior knowledge expresses all we know about the subjects when we receive the dataset. 
Consequently, counts compatible with the phenomenology 
are naturally\footnote{For finite but sufficiently large populations, it is permissible~\cite{1572261550647894016} to use instead of the hypergeometric the multinomial distribution.} sampled~\cite{hajek1960limiting,beran1983bootstrap,politis1999subsampling} from this reference point on the simplex via ($a!=\Gamma(a+1)$ denotes factorial)
\equ{
\label{eq:MutlinomialDistro}
\multidistro(N\prob p;\refP) = 
\exp \left\{ \log N! + \sum_{\config\in\cartprodR} \left[N p_\config \log\upsilon_\config - \log\left(N p_\config\right)! \right]\right\}
~.
}

\paragraph*{The large-$N$ expansion}
Appealing to large-$N$ and large-deviation techniques~\cite{dembo1993large,zinn2021quantum}, we can smoothly extend \eqref{eq:MutlinomialDistro} from relative counts to more general probabilities.
At the operational level, this means that rational relative frequencies in the Diophantine system 
are lifted over to the non-negative real numbers roughly
approximating the lattice polytope by the \totemplex
to deduce probabilities over phenomenology-compatible distributions.

At sufficiently large sample size $N$, 
the probability density to select a distribution in the differential region around $\prob p$   
becomes
\equ{
\label{eq:MutlinomialDistro:LargeN}
\multidistro\left(N\prob p;\refP\right) = 
\left(2\pi N\right)^{-\frac{\abss{\statespaceR}-1}{2}}
\exp \left\{ - N \infdiv{\prob p}{\refP}
+\frac{\abss{\statespace}-\abss{\statespace_0}}{2}\, \ell(\prob u; \prob p) 
+ \order{N^{-1}} 
\right\} 
}
using Stirling's formula. 
The large-$N$ expansion of the multinomial distribution 
thus reveals that $\refP$-based sampling  
is governed by the 
\idiv \eqref{eq:I_Divergence:DEF}
of phenomenology-compatible distribution $\prob p$ 
from reference distribution $\refP$,
\equ{
\label{eq:MutlinomialDistro:LargeN_leading}
\log\multidistro\left(N\prob p;\refP\right) \approx  - N \infdiv{\prob p}{\refP}  
~.
}

To demonstrate the logic of continuum limit, we consider the trivial setting in absence of phenomenological constraints. 
Taking for simplicity $\statespace_0=\emptyset$, 
we can replace 
discrete summation over counts $N  p_\config \in\{ 0,1,\ldots,N\}$ with integration on the simplex, 
\equ{
\label{eq:FromDiscreteToContinuous}
\prod_{\config\in\cartprod}\sum_{N p_\config} \delta\left(\sum_{\config'\in\cartprod} N p_{\config'}, N\right) 
\quads{\rightarrow} 
N^{\abss{\statespace}-1}
\int_{\mathcal P} 
~,
}
the relative mistake generically being of $\order{1/N}$ for variations \eqref{eq:Diophantine:CountVariation}.
The $\delta$ symbol on the l.h.s.\ 
stands for the Kronecker Delta. 
%
Along these lines, the volume of the 
$N$-dilation of the simplex $\mathcal P$, 
\equ{
\text{vol}\left(N \mathcal P\right) 
=
N^{\abss{\statespace}-1}\,\text{vol}\left(\mathcal P\right) 
=
\frac{N^{\abss{\statespace}-1}}{\left(\abss{\statespace}-1\right)!}
~,
}
indeed shares the same large-$N$ scaling 
\equ{
\frac{\text{vol}\left(N \mathcal P\right)}{\text{vol}\left(\mathcal P\cap\mathbb N_0^{\abss{\statespace}}\right)}
= 1 - \frac{(L-1) L}{2 N} + \order{N^{-2}}
}
with the volume~\cite{feller1968introduction} of its discrete counterpart
\equ{
\text{vol}\left(\mathcal P\cap\mathbb N_0^{\abss{\statespace}}\right) 
= 
\begin{pmatrix}
N + \abss{\statespace} - 1\\
\abss{\statespace} - 1
\end{pmatrix}
~,
}
which counts the non-negative solutions of Diophantine system \eqref{eq:Diophantine:LinearSystem} in the trivial case of $\vecOp=\{\proj_{\statespace}\}$.  
One can explicitly verify (see Appendix~\ref{app:integration_largeN}) integrating order-by-order the large-$N$ expansion of \eqref{eq:MutlinomialDistro:LargeN} the consistency of the suggested continuum approximation \eqref{eq:FromDiscreteToContinuous}.

The bounded convex polytope associated~\cite{ziegler2012lectures} to a more generic Diophantine linear system  \eqref{eq:Diophantine:LinearSystem} over $\abss{\statespaceR}$ counts
becomes complicated to deal with. Already the computation of its volume
is \textsc{np}-complete~\cite{schrijver2003combinatorial}. 
At large $N$ and in absence of many low-frequency entities nevertheless, 
we anticipate the continuum limit 
to give stable relative estimates 
about the scaling and tendencies of the discrete problem. 
In this spirit, we formally write based on \eqref{eq:MutlinomialDistro:LargeN_leading} the leading probability density 
\equ{
\label{eq:TOTEMplex:KnowledgeDistribution}
f(\prob p) = 
Z^{-1}\exp\left\{-N \infdiv{\prob p}{\refP}\right\} 
}
to select a distribution ${\prob p}$ from the \totemplex $\equivclass$  
normalized by the cumulative probability
\equ{
\label{eq:TOTEMplex:CumulativeProbability}
Z = 
\int_{\equivclass} e^{-N \infdiv{\prob p}{\refP} 
}
}
to sample any distribution on $\equivclass$. 
Notice that the probabilities of distributions on the \totemplex depend on operator $\refP\proj_{\statespaceR}$ --\,not necessarily a distribution\,-- due to $\infdiv{\prob p}{\refP\proj_{\statespaceR}} = \infdiv{\prob p}{\refP}$.

At this point, there a priori exists no clear integral prescription over generic polytopes. 
However, the $N$-scaling in the exponential of the integrand invites us to perform Laplace-type approximations of increased quality as the sample size grows. 
In Section~\ref{ssc:DistroOverDistros}, we argue how to perform the integration order-by-order in large $N$ around a stable saddle point  to estimate probability density \eqref{eq:TOTEMplex:KnowledgeDistribution}. 

\subsection{The minimum information divergence given phenomenology}
\label{ssc:MIDIV}

Clearly, leading probability density \eqref{eq:TOTEMplex:KnowledgeDistribution} of distributions on a \totemplex and 
 cumulative probability \eqref{eq:TOTEMplex:CumulativeProbability} 
are dominated with increasing $N$  by 
the minimum  \idiv of phenomenology-compatible distribution
from the reference. In any realistic problem (proof given in~\ref{app:iProjection}), 
\begin{theorem}[the \iproj]
\label{th:iProjection}
\Copy{th:iProjection}{
    Among distributions ${\prob p}\in\equivclass$ satisfying the summary statistics of  phenomenology $\vev{\prob f\vecOp}$,  
    there always exists a unique distribution signified by $\iprojSymbol[\equivclass]\refP \equiv \prob q$ 
    of minimum \idiv from reference distribution ${\refP}\in\mathcal P$:
    \equ{
        \infdiv{\prob q}{\refP}\leq \infdiv{\prob p}{\refP}
        \quad \forall~ {\prob p} \in \equivclass
    }
}
\end{theorem}
\noindent
In information geometry, the point on the simplex\footnote{$\vev{\prob p}=1$ for all  probability operators $\prob p$, thus also $\vev{\prob q}=1$.} corresponding to $\prob q\in\equivclass$ best assimilates~\cite{nielsen2018information} the properties of a projection of reference point ${\refP}$ onto the \totemplex, when the \idiv plays the role of ``distance''. %
Hence, we shall call $\iprojSymbol[\equivclass]\refP$ 
the \iproj of ${\refP}$ onto $\equivclass$. 

Let us now focus on an ``extreme'' distribution $ {\prob b}$ from the boundaries of the \totemplex $\equivclass$ where at least one $\prob b_\config$ vanishes, saturating the lower bound 
in \eqref{eq:TOTEM:LinearProgram}. 
Using any other distribution ${\prob p}\in\equivclass$  
we introduce --\,exploiting the convexity \eqref{eq:TOTEMplex:Convexity}\,--
an one-parameter family of distributions 
\equ{
\prob{l} = {\prob b} + \lambda \left({\prob p} - {\prob b}\right)
\quads{\text{for}} \lambda\in[0,1]~,
}
living on the investigated \totemplex.
Next, we look at the \idiv of $\prob{l}$ from ${\refP}$ as a convex function of $\lambda$,
whose right-derivative at the boundary equals in terms of \eqref{eq:LikelihoodForm:DEF}
\equ{
\label{eq:iProj:RightDerivative}
\lim_{\lambda\rightarrow0^+}
\frac{\dd \infdiv{\prob l}{\refP}}{\dd \lambda}
= \sum_{\statespaceR}\left( p_\config -  b_\config \right) \log\frac{ b_\config}{\upsilon_\config}
=
-\ell(\prob p;\prob b) + \ell(\prob p;\refP) - \infdiv{\prob b}{\refP}
~.
}
It is readily recognized that $-\ell(\prob p;\prob b)$ diverges as $-\infty$, if any $ b_\config$ vanishes whenever $ p_\config>0$, while the rest of the metrics remain finite (for compatible reference distribution). Hence, the rate of change of the \idiv from ${\refP}$ while leaving 
the boundary  
becomes arbitrarily large and negative. This forces us to conclude that the minimum must correspond to a distribution well in the interior of the \totemplex where all of the admissible probabilities are non-vanishing.

Having established that the \iproj cannot live on the boundaries of the \totemplex, i.e.\ its probabilities must be strictly positive, we can apply variational techniques (well established in fields such as Bayesian inference~\cite{jaakkola2000bayesian}) to find its form.
\begin{restatable}{theorem}{iProjExponentialForm}
\label{th:iProj:ExponentialForm}
The \iproj of ${\refP}$ onto $\equivclass$ defined by~\ref{def:TOTEM}, which lives in the interior of the \totemplex, can be parametrized as 
\equ{
\label{eq:MaxEnt_parametricForm}
\iprojSymbol\refP =
\exp \left( \boldsymbol \theta \cdot \vecOp \right) \refP
=
\sum_{\config\in\cartprodR}
\upsilon_\config \exp \left( \sum\limits_{\alpha=1}^{D}  \theta_\alpha\,m_\alpha(\config) \right) \ket{\config}\bra{\config}
~,
}
given the observed values in \eqref{eq:TOTEM:LinearProgram} and a vector $\boldsymbol\theta\in\mathbb R^{\abss{\statespace}-\abss{\statespace_0}}$ of appropriately fixed Lagrange multipliers.
The converse also holds: any distribution on the \totemplex of this form is the \iproj of $\refP$ onto the \totemplex.
\end{restatable}
\noindent
As before with \eqref{eq:LikelihoodForm:DEF}, we use the series definition to  resolve the exponential of an operator by
\equ{
\exp \left( \boldsymbol \theta \cdot \vecOp \right)  = \proj_{\statespaceR} + \boldsymbol \theta \cdot \vecOp + \sfrac12\left(\boldsymbol \theta \cdot \vecOp \right)^2 + \ldots~.
}
The details about the exponential form of \iproj{s}  are given in~\ref{app:ExponentialForm}.

From the variational treatment of the linear program, we verify that the \iproj does not introduce any zero probabilities, besides the declared nullentities associated to $\statespace_0$ in the \totem. 
Despite that the \iproj itself {is} unique given a phenomenology $\vev{\prob f\vecOp}$ and prior knowledge ${\refP}$,  its parameterized form \eqref{eq:MaxEnt_parametricForm} 
is not generally unique, since the parametrization in terms of Lagrange multipliers  depends on physically equivalent observations realized as linear transformations of coupled system \eqref{eq:Pheno:LinearSystem}.
Using the derived parametric form \eqref{eq:MaxEnt_parametricForm} the full version of ``Pythagorean'' relation is easily checked:
\begin{corollary}
\label{cor:Pythagoras}
For any distribution ${\prob p}\in\equivclass$ the equivalence applies:
\equ{
\infbraket{\refP}{\prob p} \geq \infbraket{\refP}{\prob q}
\quads{\Leftrightarrow}
\infbraket{\refP}{\prob p} = \infbraket{\refP}{\prob q} +  \infbraket{\prob q}{\prob p} 
~,
}
$\prob q\in\equivclass$ being the \iproj of $\refP$ onto the \totemplex.
\end{corollary}

\paragraph*{The emergence of entropy maximization}

%
Starting from minimal prior knowledge realized~\cite{keynes1921chapter} by the uniform distribution, ${\refP}={\prob u}$, 
where 
\equ{
\multidistro(N\prob p;\prob u) = \abss{\statespace}^{-N} \, W[N\prob p]
~,
}
the combinatorics of \textit{multinomial coefficient} $W[N\prob p]$
determines the probability density  of distributions on the \totemplex. $W[N\prob p]$ counts the combinatoric ways to partition $N$ identities into $\abss{\statespace}$ teams such that the counts $N \prob p_\config$ are thermodynamically reproduced.
At sufficiently large $N$, this is governed by the norm squared \eqref{eq:Norm:Entropy} coinciding with the entropy of ${\prob p}$, 
so that distributions of higher entropy, conversely closer to ${\prob u}$,
\equ{
\label{eq:Idivergence_Entropy}
H[\prob p] = - \infdiv{\prob p}{\prob u} + \log\abss{\statespace}
~,
}
contribute more to the cumulative probability \eqref{eq:TOTEMplex:CumulativeProbability} 
over the \totemplex.
In this special case that we have started from minimal prior knowledge over $\statespace$, we recover~\cite{jaynes1957information,csiszar1996maxent,jaynes2003probability} the Principle of Maximum Entropy (\maxent).

Despite that the \maxent principle has been extensively studied in data-driven literature and beyond~\cite{kapur1989maximum}, there is no 
mathematical prior in favor of the uniform distribution as a reference distribution, besides axiomatic formulation(s)~\cite{shore1980axiomatic,paris1990note,csiszar1991least}. 
In fact, any sensible ${\refP}\in\mathcal P$ that respects~\ref{def:BaseKnowledge} could mathematically serve as a valid starting point to apply the theoretical framework developed in the following.
In~\cite{loukas2023demographic} and future work, we show by applying \totem methodology to versatile problems
the utility of reference distributions well away from the \maxent prior.

\subsection{Iterative minimization of information divergence}
\label{ssc:Idiv_Minimizer}

As we have shown in the previous paragraph, the linear program \eqref{eq:Pheno:LinearSystem} supplemented by the objective of \idiv minimization (entropy maximization in the case of uniform prior knowledge) poses a well-defined optimization problem. Given phenomenological constraints induced by the $\prob f$-expectation of $\vecOp$, 
we now want to obtain the \iproj of $\refP$ on $\equivclass$
in a covariant way  that does not commit to any simplifying conditions of a particular problem.

In the context of statistical learning, the procedure of inferring from data values that parameterize a certain distribution is known~\cite{nguyen2017inverse} as the \textit{inverse} problem. In presence of (exponentially) many entities, the inverse --\,similar to the forward\,-- problem can quickly become intractable. 
In practice however, the \textsc{np}-complete bound usually lies far above a pragmatic bound, as noted in~\cite{psf2022005028}. Latter bound is posed by the size of realistic datasets (in turn bounded by finite populations of species in the universe) which severely constrains the ability to formulate statistical tests concerning the structure of more complex~\cite{PhysRevE.102.053314} relationships among subjects. 

For small kernels \eqref{eq:VectorOperator:OrthogonalComplement},  one could well parameterize $\ker\vecOp$  to analytically 
minimize the \idiv in the solution space of linear system \eqref{eq:Pheno:LinearSystem}. 
In the \maxent context, a wealth of perturbative methods~\cite{Plefka_1982,tanaka1998mean,montanari2005compute,sessak2009small} might prove useful to approximate $\iprojSymbol\prob u$. 
In special setups where (usually marginal) constraints relate attributes according to the topology of a tree, there even exists~\cite{cocco2012adaptive} a closed-form solution for the \maxent distribution. 
Beyond tree topologies iterative algorithms such as (generalized) belief propagation~\cite{yedidia2005constructing} have been used in the \maxent context for approximate inference in higher dimensional problems. 

More generally, starting from some reference distribution ${\refP}$ we 
iteratively impose constraints~\eqref{eq:TOTEM:LinearProgram} that we summarize in a vector function $\mathbf F$ with 
\equ{
\label{eq:Newton:constraintFunctions}
F_{\alpha}[\prob p] = 
\vev{(\prob p-\prob f)\op M_\alpha}
=\sum_{\config\in\cartprodR}  m_\alpha(\config) \left( p_{\config} -  f_{\config}\right)
\quads{\text{for}}\alpha=1,\ldots D
~,
}
where $ p_\config$ are parametrized according to \eqref{eq:MaxEnt_parametricForm}.
A well-known numerical method to find the root in the resulting  non-linear system of coupled equations is the Newton-Raphson algorithm, 
which has inspired many iterative methods~\cite{malouf2002comparison} for approximate inference.
%
In probability space, we can covariantly formulate  (convergence discussed in~\ref{app:Newton})
\begin{restatable}{theorem}{NewtonRaphson}
\label{th:NewtonRaphson}
The iterative algorithm to covariantly  approach the probabilities of the \iproj of ${\refP}$ onto $\equivclass$, expressed for $n\in\mathbb N$ in matrix-operator notation 
\equ{
\prob p^{(0)} = \refP \proj_{\statespaceR}
\quads{,}
\mathbf J^{(n)} = \vev{\vecOp\, \prob p^{(n)} \vecOp[M^T]}
\qand
\prob p^{(n+1)} =\exp\left\{- \vecOp[M^T] (\mathbf J^{(n)})^{-1} \mathbf F[\prob p^{(n)}]
\right\} \prob p^{(n)}
\nonumber
}
or component-wise as
\begin{align}
\label{eq:NewtonStep}
 p^{(0)}_\config = &\,\,\upsilon_\config
&\quads{\text{for}}& \config\in\cartprodR
\\[1ex]
J_{\alpha\beta} =&\,\, 
\vev{\op M_\alpha\, \prob p^{(n)} \op M_\beta}
=
\sum_{\config\in\cartprodR} m_{\alpha}(\config)\, p_\config^{(n)} \,m_{\beta}(\config)
&\quads{\text{for}}&\alpha,\beta=1,\ldots,D
\nonumber
\\[0.8ex]
 p^{(n+1)}_\config =&\,\, p^{(n)}_\config \exp\left\{ - \sum_{\alpha,\beta=1}^{D} m_\alpha(\config) \left(\mathbf J^{-1} \right)_{\alpha\beta} 
F_\beta[\prob p^{(n)}]
\right\}
&\quads{\text{for}}& \config\in\cartprodR
\nonumber
\end{align}
is well-behaved, whenever $\iprojSymbol\refP$ is reasonbly close to  reference distribution ${\refP}$.
\end{restatable}
\noindent
The corresponding update rule to numerically determine the Lagrange multipliers in \eqref{eq:MaxEnt_parametricForm} reads 
\equ{
\label{eq:Newton:LagrangeMultipliers}
\theta^{(n+1)}_\alpha = \theta^{(n)}_\alpha - \sum_{\beta=1}^{D}  \left(\mathbf J^{-1}\right)_{\alpha\beta} F_\beta[\prob p^{(n)}] 
~.
}
Concerning the value $F^{(n)}_\alpha \equiv F_\alpha[\prob p^{(n)}]$ of the vector function \eqref{eq:Newton:constraintFunctions} whose unique root we are after,  we recognize after $n$ iterations by the scalar product
\equ{
\sum_{\alpha=1}^D
F_\alpha^{(n)}\left(\theta^{(n+1)}_\alpha - \theta^{(n)}_\alpha\right)
=
- \sum_{\alpha,\beta=1}^D F_\alpha^{(n)} \left(\mathbf J^{-1}\right)_{\alpha\beta} F_\beta^{(n)} < 0 \quad\forall \,n\quads{\text{s.t.}} F^{(n)}_\alpha\neq0
\nonumber
}
using \eqref{eq:Newton:LagrangeMultipliers} and positive-definiteness of $\mathbf J^{-1}$
that the procedure \eqref{eq:NewtonStep} 
moves along a descent direction until the desired convergence.

The covariance of algorithm~\ref{th:NewtonRaphson} manifests due to the affine invariance~\cite{boyd_vandenberghe_2004} of Newton-Raphson method. 
A physically equivalent re-parametrization of measuring operations performed via $\mathbf T$ in  \eqref{eq:TOTEMs:Reparametrization} relates two \textsc{fapp}-equivalent vector operators constructing the same \totemplex, hence leading to the same \iproj of $\refP$.
%
Accordingly, constraining function \eqref{eq:Newton:constraintFunctions} covariantly transforms under $\mathbf T$ as a vector,
\equ{
F_\alpha 
\quads{\rightarrow} \sum_{\beta=1}^D T_{\alpha\beta} F_{\beta} 
~,
}
while the Jacobian  $\mathbf J$ as a matrix,
\equ{
\mathbf J_{\alpha\beta} \quads{\rightarrow} \sum_{\alpha'\beta'=1}^D T_{\alpha,\alpha'}\, J_{\alpha'\beta'}\,T_{\beta,\beta'}
~.
}
In total, we readily verify that the update rule of \eqref{eq:NewtonStep} remains invariant under re-parametrization of measuring operations \eqref{eq:TOTEMs:Reparametrization} exemplifying the merits of working in probability space.
Not only this property programmatically streamlines the generation of model distributions, but helps setup approximate algorithms~\cite{cocco2012adaptive,psf2022005028} --\,independently of the current framework\,-- for inference in higher-dimensional problems where the \textsc{np}-completeness manifests.

\paragraph*{Chain of \iproj{s} in nested descriptions}
Since the Newton-Raphson estimate of~\ref{th:NewtonRaphson} remains agnostic to the concrete parametrization of Theorem~\ref{th:iProj:ExponentialForm}, convergence of the iterative procedure solely\footnote{At the technical level, the inversion of the Jacobian at each Newton-Raphson step might be poorly conditioned for some parametrizations of constraints eventually de-railing the iterative estimate.} depends on the physical setup. 
If the reference distribution  ${\refP}$ poorly captures phenomenology $\vev{\prob f\vecOp}$  resulting into higher \idiv  $\infdiv{\iprojSymbol\refP}{\refP}$, 
the algorithm might fail to converge to the \iproj, as generally expected by Newton-based methods~\cite{boyd_vandenberghe_2004,nesterov2003introductory}.

\begin{definition}[nested descriptions]
    \label{def:NestedTOTEMs}
    Consider two irreducible constructing elements $\vecOp$ and $\vecOp[M']$  over $\statespaceR$ related  for $D<D'$  through
    \equ{
        \label{eq:LRT:NestedTOTEMs}
        \op M_{\alpha} =  \sum_{\alpha'=1}^{D'} S_{\alpha\alpha'} \op M^\prime_{\alpha'} 
        \quads{\text{for}} \alpha=1,\ldots,D 
    }
    where $\rank\mathbf S=D$ (otherwise $\vecOp$ would be row-rank deficit).
    The \totem{s} defined by phenomenologies $\vev{\prob f\vecOp}$ and $\vev{\prob f\vecOp[M']}$ 
    and the resulting \totem{lices}    $\equivclass[{\vecOp[M']}]\subseteq\equivclass[\vecOp]$ are called {nested}.
\end{definition}
\noindent
The fact that the induced \totem{plices} are subsets can be readily seen from linearity of expectations: 
\begin{align*}
{\prob p}\in\equivclass[{\vecOp[M']}]:\quad
\vev{\left({\prob p}-{\prob f}\right)\op M_\alpha}
= \sum_{\alpha'=1}^{D'} S_{\alpha\alpha'} \underbrace{\vev{\left({\prob p}-{\prob f}\right)\op M'_{\alpha'}}}_{=0} 
=0
\quad\Rightarrow\quad
{\prob p}\in\equivclass[\vecOp]
~.
\end{align*}
The more knowledgeable we become about the provided data via additional measurements performed as linearly independent $\prob f$-expectations,
the further we generically move away from ${\refP}$, since 
\equ{
\label{eq:Newton:MoreKnowledgeable}
\infbraket{\refP}{\iprojSymbol[{[\prob f;\vecOp']}]\refP} \geq \infbraket{\refP}{\iprojSymbol[\equivclass]\refP} 
}
by construction~\ref{th:iProjection} of the \iproj of ${\refP}$ onto ambient \totemplex $[\prob f;\vecOp]$ where $\iprojSymbol[{[\prob f;\vecOp']}]\refP$ also belongs to.

%
Exploiting this fact, we might still try to numerically determine the \iproj of ${\refP}$ onto $\equivclass$ constructed by some vector operator $\vecOp$,
for which the iterative algorithm of~\ref{th:NewtonRaphson} failed to converge, by progressively projecting onto
nested \totem{plices}.
Specifically,  
\begin{restatable}{theorem}{nestedNewtons}
\label{th:nested_Newtons}
Partition $\vecOp$ into $d\leq D$ constructing elements $\vecOp[M^{(\nu)}]$
inducing nested \totem{s} in the sense of~\ref{def:NestedTOTEMs},
\equ{
\equivclass[\vecOp] \subseteq \ldots \equivclass[{\vecOp[{M^{(\nu)}}]}] \subseteq \ldots 
\equivclass[{\mathbf{\hat{M}}^{(1)}}] \subseteq \mathcal P(\statespaceR) 
~,
}
where eventually $\mathbf{\hat{M}}^{(d)}$ is \textsc{fapp} equivalent to the given $\vecOp$.
Starting from reference $\refP^{(0)} = \refP\proj_{\statespaceR}\in\mathcal P(\statespaceR)$, 
we successively use iterative algorithm~\ref{th:NewtonRaphson} 
to determine the \iproj of 
\equ{
\label{eq:Newton:SubRoutines}
\refP^{(\nu-1)} 
=
\iprojSymbol[{[\prob f;\vecOp[{M^{(\nu-1)}}]]}]\cdots\iprojSymbol[{[\prob f;\mathbf{\hat M}^{(1)}]}]\refP
\,\in \equivclass[{\vecOp[{M^{(\nu-1)}}]}]
}
onto  inner \totemplex $\equivclass[{\vecOp[M^{(\nu)}]}]$ for $\nu=2,\ldots,d$. Eventually, the last estimate
coincides with the desired \iproj, $
\refP^{(d)} = \iprojSymbol \refP
$.
\end{restatable}
\noindent
After $\nu$ steps, trying to numerically find the distribution on $\equivclass[{\mathbf{\hat M}^{(\nu+1)}}]$ closest to updated reference distribution ${\refP^{(\nu)}}$ (itself living on the ambient $\equivclass[{\vecOp[M^{(\nu)}]}]$) is expected to exhibit better convergence properties than  directly projecting from ${\refP^{(\nu-1)}}$ onto inner \totemplex $\equivclass[{\mathbf{\hat M}^{(\nu+1)}}]$. 
By merit of
Corollary~\ref{cor:Pythagoras} applied within 
$\equivclass[{\vecOp[M^{(\nu)}]}]$
where $\refP^{(\nu)} = \iprojSymbol[{\equivclass[\vecOp^{(\nu)}]}]\refP^{(\nu-1)}$ is the $I$-projection of $\refP^{(\nu-1)}$ and  
$\refP^{(\nu+1)}=\iprojSymbol[{\equivclass[{\vecOp[M^{(\nu+1)}]}]}]\refP^{(\nu)}\in\equivclass[{\vecOp[M^{(\nu+1)}]}] \subset\equivclass[{\vecOp[M^{(\nu)}]}]$ we recognize that 
\equ{
\infdiv{\refP^{(\nu+1)}}{\refP^{(\nu-1)}} 
=
\infdiv{\refP^{(\nu)}}{\refP^{(\nu-1)}} + \infdiv{\refP^{(\nu+1)}}{\refP^{(\nu)}}
\geq \infdiv{\refP^{(\nu+1)}}{\refP^{(\nu)}}
\nonumber
~.
}

Pipeline \eqref{eq:Newton:SubRoutines} breaks the original problem of constrained minimization of the \idiv into root-finding 
sub-routines whose target lies closer to 
initial estimate.
This only makes sense 
iff the consecutive determination of \iproj{s} onto nested \totem{plices} is equivalent to directly taking the \iproj of outermost reference distribution onto  the inner-most \totemplex. 
Indeed, this is asserted in~\ref{app:nestedTOTEMs}. 

\paragraph*{Iterative Proportional Fitting}
There exist many variations of the main iterative scheme.  
Suppose that 
$M$ marginal, not necessarily linearly independent operators $\proj_\alpha$ with eigenvalues like in \eqref{eq:MarginalConstraintFkt}  are provided. Their combined action on entities can always be encoded by 
a binary matrix $\mathbf C$, such that 
\equ{
\vev{\prob p\proj_\alpha} = \sum_{\config\in\cartprodR} C_{\alpha,\config}\,  p_\config
~.
}
Cycling through each row depicting a marginal constraint for $\alpha=1,\ldots M$ to iteratively update probabilities via one-dimensional Newton method (as the multidimensional Newton method would be plagued by the singularity of redundancies) gives after $n$ cycles the  rule
\equ{
\label{eq:IPF}
 p^{(nM+\alpha)}_\config =  p^{(nM+\alpha-1)}_\config  \exp\left\lbrace C_{\alpha, \config} \left(\frac{\vev{\prob f\proj_\alpha}}{\vev{\prob p^{(nM+\alpha-1)}\proj_\alpha}} -1\right)\right\rbrace 
\quad\forall\,\config\in\cartprodR
~.
}

Using $C_{\alpha,\config}^2=C_{\alpha,\config}$
in the series expansion of the exponential up to first order\footnote{Equivalently, we could set $\vev{\prob p^{(nM+\alpha-1)}\proj_\alpha}=\vev{\prob f\proj_\alpha}/x$ to Taylor-expand around $x=1$ up to first order.} gives
\equ{
 p^{(nM+\alpha)}_\config \approx  p^{(nM+\alpha-1)}_\config 
\left(\frac{\vev{\prob f\proj_\alpha}}{\vev{\prob p^{(nM+\alpha-1)}\proj_\alpha}}\right)^{C_{\alpha, \config}}
~.
}
One can show~\cite{csiszar1975divergence,ireland1968contingency} that the Taylor-expanded version
always converges to the desired \iproj starting from a sensible $\prob p^{(0)}=\refP$. This approach leads to the algorithm of iterative proportional fitting~\cite{kruithof1937telefoonverkeersrekening} with its own long history~\cite{darroch1972generalized,haberman1974log,10.2307/2242759,bishop2007discrete,phdthesis_IPF,loukas2022categorical}.

\subsection{The distribution of phenomenology-compatible probability fluctuations}
\label{ssc:DistroOverDistros}

After establishing existence and uniqueness of the \iproj $\prob q\equiv\iprojSymbol\refP$ given prior knowledge and a \totem of the data,
we estimate the probability to sample  distributions on the \totemplex by performing a large-$N$ expansion in the probability fluctuations around the \iproj. 

First, let us examine --\,again in the discrete problem\,-- 
the relative abundance of a count vector predominantly sampled from $\refP$
over another  count vector in the Diophantine class. Following the logic of Section~\ref{ssc:DicreteContinuous},  such ratio can be well represented by 
\equ{
\label{eq:Motivation:ExponentialDominanance}
\max_{\prob p'\in\equivclass}\log \frac{{\multidistro}(N\prob p';\refP)}{{\multidistro}(N\prob p;\refP)}
\overset{N\gg1}{\approx} N\left[\infdiv{\prob p}{\refP} - \infdiv{\prob q}{\refP}\right] \geq 0
~.
}
This clearly 
demonstrates that counts close to  
$N\prob q$ 
should be overwhelmingly anticipated given the empirical constraints from phenomenology~\ref{def:Phenomenology} and compatible  prior knowledge~\ref{def:BaseKnowledge}, both at finite\footnote{The finite-$N$ combinatorics of the Diophantine problem and the associated notion of typical count vector(s) are explored in~\cite{loukas2022categorical}.} and large $N$ and the latest asymptotically as $N\rightarrow\infty$. 
The probability of any other member vector $N\prob p$ 
whose \idiv from $\refP$ differs by more than $1/N$ from $\infdiv{\prob q}{\refP}$ would be exponentially suppressed. 

Intuitively, the dominating role of $\prob q$ 
in the equivalence class explains why the \iproj is the most typical of the \totem based on our prior knowledge.
Intuitively, this is perhaps best understood starting from minimal prior  knowledge, when the minimization of \idiv reduces to the \maxent. In that case, the ratio of multinomial coefficients for some count vector in the Diophantine class,
\equ{
\log\frac{W[N\iprojSymbol\prob u)]}{W[N\prob p]} \overset{N\gg1}{\approx} N\left(H[\iprojSymbol\prob u] - H[\prob p]\right)
~,
\nonumber
}
clearly shows that 
any count vector closest to $\iprojSymbol\prob u$, which is the least biased given the constraints~\cite{jaynes1968prior}, can be realized in the most combinatorial ways, 
as the entropy difference will be by definition  positive.

\newcommand{\dimchi}{k} 
\newcommand{\probflac}{\boldsymbol{\pi}}

Regarding the probability density at some $\prob p\in\equivclass$ in the continuum,
the combinatorics of Section~\ref{ssc:DicreteContinuous} have (approximately) motivated 
\eqref{eq:TOTEMplex:KnowledgeDistribution},
 whose $N$-leading behaviour can  be rewritten  by Corollary~\ref{cor:Pythagoras} 
as 
\equ{
\label{eq:TOTEMplex:KnowledgeDistribution2}
f(\prob p) = Z^{-1}\,\exp\left\{ - N \infbraket{\prob q}{\prob p} \right\} + \ldots
}
with redefined partition function
\equ{
\label{eq:TOTEMplex:PartitionFunction}
Z = \int_{{\prob p}\in\equivclass} \exp\left\{ - N \infbraket{\prob q}{\prob p} \right\}
~.
}
To systematically investigate this probability density function $(\textsc{pdf})$ over member distributions on the \totemplex constructed by $\vecOp$ on data described by $\prob f$, we exploit the moral of the discrete observations above.

\paragraph*{Probability fluctuations}
 Controlled deviations around the special point on the \totemplex represented by the \iproj 
would prove particularly useful in estimating the aforementioned \textsc{pdf}.
Such probability fluctuations $\probflac$
around the \iproj of ${\refP}$ onto the \totemplex  span  
the appropriate polyhedral cone such that
\equ{
 \label{eq:Fluctuations:PolyhedralCone}
\prob p = \iprojSymbol\refP + {\boldsymbol\pi}
}
describes a valid distribution from $\mathcal P(\statespaceR)$ with non-negative probabilities $p_\config\geq0$. 
Since by definition $\vev{\probflac\vecOp}=0$, fluctuation operator $\probflac\in\ker\vecOp$ can be accommodated along the ``normal'' basis of \eqref{eq:VectorOperator:OrthogonalComplement},
\equ{
\label{eq:KernelFluctuations:Expansion}
\probflac = \sum_{a=1}^{\abss{\ker\vecOp}}\probflac_a \op N_a
~.
}
As we are going to assert below, a large-$N$ expansion in the neighbourhood of the \iproj enables us to extract the universal behavior of metrics in Section~\ref{ssc:InformationTheoreticMetrics} without worrying about the conic boundaries.
Since the \iproj lives well in the interior of the \totemplex (cf.\ Theorem~\ref{th:iProj:ExponentialForm}), we anticipate to always find some parametrization that does not degenerate, thus ensuring that probability fluctuations $\probflac_a=\order{1}$.

%


Due to convexity, any distribution on the \totemplex can be reached  by fluctuating around some base point. 
%
%
Consequently, the study of probability density  \eqref{eq:TOTEMplex:KnowledgeDistribution2} becomes equivalent to investigating how fluctuations distribute,
$f(\prob p) \equiv f(\probflac)$.
Studying fluctuations around some arbitrary base distribution on the \totemplex however, would not help much in computing the \textsc{pdf}. 
On the other hand, fluctuating around the \iproj 
of our reference distribution onto the \totemplex 
enables to formulate
\begin{theorem}[Gaussian approximation]
\label{th:Universality:GaussianApproximation}
Consider a \totem constructed by $\vecOp$ over $\statespaceR$ 
and the \iproj  of reference distribution $\refP$ onto the \totemplex. 
At large $N$, 
 probability  fluctuations in 
$\ker\vecOp$ 
around the \iproj of $\refP$ are controlled by Gaussian approximation in the sense of 
predominantly following
a
$\abss{\ker\vecOp}$-dimensional 
normal distribution of zero mean 
and $\order{1/N}$ variance.
\end{theorem}
\paragraph*{Proof} 
After deriving 
the Gaussian kernel that drives probability fluctuations around the \iproj, an integral representation for the partition function \eqref{eq:TOTEMplex:PartitionFunction} of the $N$-leading \textsc{pdf} would be naturally motivated. Subsequently, the universal character of the expansion would be self-consistently  guaranteed by showing 
 suppression at large $N$ of any higher-order corrections 
 to the \textsc{pdf} over the \totemplex.
For that, we need (proof in~\ref{app:UniversalKernel})
\begin{restatable}{lemma}{UniversalKernel}
\label{lm:UniversalKernel}
    The \idiv of any distribution ${\prob p}\in\equivclass$ from the \iproj of ${\refP}$ onto $\equivclass$ is universally (i.e.\ independently of the phenomenological problem) determined in the vicinity of the \iproj by a quadratic form $\mathcal K$ over the kernel of $\vecOp$:
    \equ{
    \label{eq:KnowledgeDistro:UniversalKernel}
        \infbraket{\iprojSymbol\refP}{\prob p} =
        \sfrac{1}{2} \mathcal K(\probflac) 
        + \order{\probflac^3}
        \quads{\text{with}} \mathcal K(\probflac) = \vev{\probflac \prob q^{-1} \probflac} 
    }
\end{restatable}

The resulting precision matrix $\mathcal K(\probflac)$ 
is positive definite
for the reasons outlined in Lemma~\ref{lm:Jacobian} given irreducible vector operator $\vecOp[N]$ spanning the kernel of $\vecOp$. Hence, a Gaussian approximation to $f(\prob p)\equiv f(\probflac)$ in kernel space around the \iproj is always well-defined. 
In turn, this motivates the evaluation of the partition function $Z$ over the \totemplex via a multidimensional Gaussian integral 
whose limits  can be extended  to infinity, the error being exponentially suppressed in $N$.

Furthermore, Laplace approximation in \eqref{eq:TOTEMplex:KnowledgeDistribution2} ensures that higher-order fluctuations are polynomially suppressed in $N$, turning the expansion in fluctuations \eqref{eq:KnowledgeDistro:UniversalKernel} around the \iproj into a large-$N$ expansion.
The same applies for higher probability-dependent terms in the large-$N$ expansion \eqref{eq:MutlinomialDistro:LargeN} of multinomial distribution. Generically, these are expected to appear at $\order{N^{-1/2}}$. 
In total, we can self-consistently verify order-by-order in large-$N$ the universal form 
\equ{
\label{eq:KnowledgeFluctuations:UniversalGaussian}
f(\probflac) = Z^{-1} \exp\left\{-N\mathcal K(\probflac)/2\right\} + \order{N^{-1/2}}
}
of $N$-leading \textsc{pdf} of probability fluctuations around the \iproj on the \totemplex. 
$\qed$\\

As perhaps anticipated by the law of large numbers, probability fluctuations in a \totem become at sufficiently large sample size almost normally distributed, once expanded around the \iproj of our reference distribution onto the \totemplex. The dominating role of the \iproj discussed in \eqref{eq:Motivation:ExponentialDominanance} is quantitatively reflected on the variance of the multivariate normal distribution which quickly penalizes distributions whose \idiv from $\iprojSymbol\refP$ is larger than  $\order{1/N}$. In other words using Lemma~\ref{lm:UniversalKernel} together with \eqref{eq:KnowledgeFluctuations:UniversalGaussian} and \eqref{eq:appGaussian:QuadraticExpectation}, the expected \idiv of distributions on the \totemplex from the \iproj scales as
\equ{
\label{eq:TOTEMplex:expected_iDivergence}
\int_{\prob p\in\equivclass}f(\prob p) \infdiv{\prob p}{\iprojSymbol\refP}  \approx 
\prod_{a=1}^{\abss{\ker\vecOp}} \int \dd\probflac_a f(\probflac)\, \sfrac12\mathcal K(\probflac) 
=
\order{1/N}
~.
}
In the context of \maxent principle, this observation translated by Corollary~\ref{cor:Pythagoras} and \eqref{eq:Idivergence_Entropy} into an entropy difference,
\equ{
\infbraket{\iprojSymbol\prob u}{\prob p} = \infbraket{\prob u}{\prob p} - \infbraket{\prob u}{\iprojSymbol\prob u} =  H[\iprojSymbol\prob u] - H[\prob p]~,
}
led authors talk~\cite{Rosenkrantz1989} about the \textit{concentration} of entropies with increasing $N$ around the \maxent distribution.

According to Lemma~\ref{lm:UniversalKernel}, the quadratic form which controls probability fluctuations around the \iproj only depends on the \idiv of \totem distributions from the \iproj. 
This universal fact always becomes relevant as the sample size grows. 
In the reminder of this paragraph, we exploit the geometric symmetry of the Gaussian approximation to translate Theorem~\ref{th:Universality:GaussianApproximation} into operationally useful prescriptions to eventually select appropriate operators in the following section. 

\begin{restatable}{corollary}{SphericalSymmetry} 
\label{cor:ShpericalSymmetry}
On the \totemplex, 
consider
the probability of a distribution ${\prob p}\in\equivclass$ with definite \idiv 
from the \iproj $\prob q\equiv\iprojSymbol\refP$ given prior  knowledge ${\refP}$.  
  To leading order in the large-$N$ expansion, the probability density enjoys 
spherical symmetry 
in $k\equiv\abss{\ker\vecOp}$ dimensions with radial density  
\equ{
\label{eq:ToTEMplex:spherical_ProbabilityDensity}
\rho(\prob p) \equiv \rho(\infbraket{\prob q}{\prob p}) = \frac12\left(\frac{N}{\pi}\right)^{k/2} \infbraket{\prob q}{\prob p}^{k/2-1} \exp\left\{-N \infbraket{\prob q}{\prob p}\right\}
~,
}
dictated by the specified \idiv $\infbraket{\prob q}{\prob p}$ 
which plays the role of radius squared in $\ker\vecOp$.
\end{restatable}
\noindent
The way to read off the radial density from \eqref{eq:TOTEMplex:KnowledgeDistribution2} is outlined in~\ref{app:shpericalSymmetry}.
According to probability density $\rho$  any distribution on the \totemplex living on the $(k-1)$-sphere of squared radius $\infbraket{\prob q}{\prob p}$ will be assigned the same weight from the perspective of reference distribution.

To connect Corollary~\ref{cor:ShpericalSymmetry} with more established conventions in the literature, 
we 
trivially execute the angular integration in \eqref{eq:ToTEMplex:spherical_ProbabilityDensity} (cf.\ \eqref{eq:KnowledgeFluctuations:RadialProbability_K})
to obtain a cumulative radial probability
\equ{
\label{eq:TOTEMplex:RadialProbability}
\dd\infbraket{\prob q}{\prob p}\, \frac{N^{k/2}}{\Gamma(k/2)}\infbraket{\prob q}{\prob p} ^{k/2-1}\, \exp\left\{-N \infbraket{\prob q}{\prob p}\right\}
}
 on the $(k-1)$-sphere.
Subsequently, we set 
\equ{
\label{eq:ChiSqaureStatistics}
Q = 2N\infbraket{\prob q}{\prob p}
}
to immediately 
arrive 
at the chi-squared probability in $k$ dimensions
\equ{
\label{eq:TOTEMplex:ChiSquaredDistro}
\dd Q \,
\frac{1}{2^{k/2}\Gamma(k/2)}Q^{k/2-1}\, e^{-Q/2}
\quads{\text{for}} Q\geq0
~,
}
which selects any distribution ${\prob p}\in\equivclass$ whose \idiv from the \iproj lies on a spherical shell of radius $\frac{Q}{2N}$ and thickness $\dd Q$. 
This describes a chi-squared distribution with $k\equiv\abss{\ker\vecOp}$ degrees of freedom and statistics given by \eqref{eq:ChiSqaureStatistics}.

%
Due to the non-linear optimization problem leading to Theorem~\ref{th:iProj:ExponentialForm}, it becomes evident that the \iproj, which crucially depends on reference distribution ${\refP}$, does not generically coincide with the notion of a mean  distribution over the \totemplex. 
Already with $\abss{\statespace}=3$ entities, it is easy to find counter-examples both in the discrete and continuous case. Consequently, the developed large-$N$ expansion cannot be understood as mere application of central limit theorem~\cite{van2000asymptotic} on  probabilities --\,treated themselves as $\abss{\statespace}$ random variables. 
By the same token, the variance of probability distributions implied by \eqref{eq:TOTEMplex:expected_iDivergence} does not coincide with the ``naive'' variance from a mean probability on the \totemplex.

\section{The classification of phenomenologies}
\label{sc:ClassificationOfKnowledge}

Starting from prior knowledge about subjects, we have so far motivated  
--\,given some phenomenology~\ref{def:Phenomenology} of the data\,-- 
the distribution of compatible 
probability fluctuations around the \iproj of reference distribution onto the \totemplex.
Performing however different measurements on the same dataset, i.e.\ taking $\prob f$-expectations of different vector operators $\vecOp$
would
yield a different phenomenology, generically constructing a distinct \totem. 
In this Section, we describe a way to select appropriate operators inducing optimally descriptive phenomenologies given prior knowledge and the data.
Ultimately, we relate to familiar notions of statistical testing 
in our unified framework of large-$N$ description.

\paragraph*{Quotient space}
An often encountered task in data analysis is to assess the amount of knowledge gained by performing additional measurements through $\prob f$-expectations.
In other words, 
we wish to decide whether a more detailed phenomenology $\vev{\prob f\vecOp[M']}$ extracts from the data additional valuable knowledge 
on the investigated subjects  compared to  more plain  phenomenology $\vev{\prob f\vecOp}$.
From the perspective of our prior  knowledge, we would have to compare the ambient \totemplex constructed by vector operator $\vecOp$ 
to the nested \totemplex constructed by $\vecOp[M']$ which \textsc{fapp} implies the former in the sense of definition~\ref{def:NestedTOTEMs}.

Under the action of $\vecOp[M']$ probability operators in the ambient \totemplex $\equivclass$  get identified according to 
\equ{
\label{eq:QuotientSpace:EquivalenceRelation}
{\prob p}, {\prob p'} \in \equivclass~:\quad
{\prob p'} \sim {\prob p} \quads{\Longleftrightarrow}  
{\prob p'}\in\equivclass[{\vecOp[M']}][\prob p]\subset\equivclass.
}
This means that we should not distinguish between two distributions 
that both satisfy the phenomenology $\vev{\prob f\vecOp}$, if they also agree on the higher statistics captured by the phenomenology $\vev{\prob f\vecOp[M']}$. 
Conversely, two distributions that differ by some expectations in $\vev{\prob f\vecOp[M']}$ are called  $\vecOp[M']$-distinct.
The set of $\vecOp[M']$-induced equivalence classes defines the quotient space $\equivclass\slash\vecOp[M']$ of the ambient \totemplex.
Theorem~\ref{th:iProj:ExponentialForm} motivates a parametrization of the set of $\vecOp[M']$-distinct operators
where Lagrange multipliers continuously vary over representatives in the equivalence classes of the quotient space.


Having access to the full data $\texttt X$ we can construct orbits of nested \totem{s} over $\statespaceR$, 
\equ{
\label{eq:chain_nestedTOTEMs}
\{\prob f\}\equiv \equivclass[\Id] \subseteq\ldots \subseteq\equivclass[{\vecOp[M']}][]\subseteq\equivclass \subseteq \ldots \subseteq\equivclass[\hat I]\equiv\mathcal P(\statespaceR)
~,
}
which interpolate between learning empirical frequencies by heart via the identity operator \eqref{eq:IdentityOperator} and not learning anything from the data except normalization condition~\eqref{eq:Probabilities:Normalization}. 
Traversing the chain of nested \totem{plices} from right to left, the amount of valuable knowledge that each phenomenology~\ref{def:Phenomenology} at fixed sample size $N$ additionally conveys depends on the underlying (usually unknown) laws that govern the investigated subjects as well as on our prior knowledge about them.


\paragraph*{The information score}
Given the data, we first motivate a score 
to systematically benchmark  from the perspective of prior knowledge distinct \totem descriptions within an orbit like \eqref{eq:chain_nestedTOTEMs}, but also among different not necessarily nested \totem{s}.
For that, we proceed with 
\begin{definition}
\label{def:iScore}
    To each \totem we associate at sufficiently large sample size $N$ the information score 
    \equ{
    \label{eq:MLE:InformationScore}
        \textsc{s}core_I = -N \infbraket{\iprojSymbol\refP}{\prob f} + \frac{\abss{\ker\vecOp}}{2}\log N  
        + \order{1} 
    ~,
    }
    in terms of the \iproj $\iprojSymbol\refP$ onto the \totemplex $\equivclass$. 
\end{definition}
Clearly, this information score corresponds to the logarithm of the weight assigned to empirical distribution $\prob f$
by the  radial density~\ref{cor:ShpericalSymmetry} over $\ker\vecOp$ that governs probability fluctuations around the \iproj
%
onto the \totemplex. 
Concretely, the two leading terms determine how radial probability \eqref{eq:TOTEMplex:RadialProbability} would change by an infinitesimal variation in the vicinity of $\infbraket{\iprojSymbol\refP}{\prob f}$.

In our data-centric approach, 
the information criterion implied by~\ref{def:iScore} crucially depends on the \iproj as well as the degrees of freedom on the \totemplex to covariantly 
assign weights to \totem{s}, i.e.\ to operator collections given $\prob f$ and $\refP$.   
Pragmatically, this becomes useful by 
\begin{proposition}[$I$-score]
\label{prop:Score}
    Over a set of \textsc{fapp} inequivalent vector operators~\ref{def:Irreducible_VectorOperator}, 
    the principle of Maximum Likelihood dictates to select the vector operator with highest information score,
    \equ{
        \vecOp_\textsc{mle} = \max_{\vecOp} \textsc{s}core_I
        \nonumber
    }
    which maximizes the likelihood of the data in the respective \totem. 
\end{proposition}

\paragraph*{The distribution over the quotient space}
Comparing two nested (in the sense of Definition~\ref{def:NestedTOTEMs}) \totem{s}  constructed by distinct phenomenologies $\vev{\prob f\vecOp}$ and $\vev{\prob f\vecOp[M']}$ can be additionally performed via hypothesis testing.
For that, we need to compute the probability of a \totemplex, i.e.\ the cumulative probability to sample any distribution on the \totemplex. In particular, we consider inner \totemplex $\equivclass[{\vecOp[M']}][\prob p]\subset\equivclass$ represented by $\vecOp[M']$-distinct distribution $\prob p$.
Emphasizing the importance of the \iproj to leading order in $N$, 
the probability to sample $\equivclass[{\vecOp[M']}][\prob p]$  in the quotient space 
$\equivclass\slash\vecOp[M']$ is dictated by  
\equ{
\label{eq:Quotient:KnowledgeDistribution}
Z^{-1} 
\exp\left\{-N \infdiv{\iprojSymbol[{\equivclass[{\vecOp[M']}][\prob p]}]\refP}{\iprojSymbol\refP}\right\}
~,
}
as derived in Lemma~\ref{lm:Quotient:KnowledgeDistribution}. The exponential factor can also be understood as a ratio of $N$-leading probability mass \eqref{eq:TOTEMplex:ProbabilityMass} of inner \totemplex to the probability mass of ambient \totemplex upon using Corollary~\ref{cor:Pythagoras}.

In accordance with the quotient action, 
probability fluctuations in \eqref{eq:KernelFluctuations:Expansion} along $\vecOp[N]$
(spanning the ambient kernel of $\vecOp$) can be always decomposed into $\abss{\ker\vecOp[M']}$ directions $\op N_{a'}$ along inner $\ker\vecOp[M']$ and remaining normal directions $\op N_a$ with $\vev{\op N_a\op N_{a'}}=0$. 
Investigating the distribution of probability fluctuations along latter directions it is straight-forward to verify \eqref{eq:Quotient:KnowledgeDistribution} and formulate (proven in~\ref{app:LRT})
\begin{restatable}{theorem}{LRT}
\label{th:LRT}
    In the quotient space induced by $\vecOp[M']$ over $\equivclass$ consider
    inner \totemplex $\equivclass[{\vecOp[M']}][\prob p]$ 
    with \iproj $\iprojSymbol[{\equivclass[{\vecOp[M']}][\prob p]}]\refP$. 
    We examine the cumulative probability to sample from $\refP$ any 
    inner \totemplex  whose 
    \iproj has a fixed \idiv $\infdiv{\iprojSymbol[{\equivclass[{\vecOp[M']}][\prob p]}]\refP}{\iprojSymbol\refP}$ from the ambient \iproj.
    At large $N$, 
    this probability is dictated by the chi-squared distribution of $2N\infdiv{\iprojSymbol[{\equivclass[{\vecOp[M']}][\prob p]}]\refP}{\iprojSymbol\refP}$ with $D'-D$ degrees of freedom.
\end{restatable}
\noindent
Whenever $\vecOp[M']=\Id$, $\equivclass[{\Id}][{\prob p}]=\{\prob p\}$ so that probability density \eqref{eq:Quotient:KnowledgeDistribution} reduces back to \eqref{eq:TOTEMplex:KnowledgeDistribution2} and Theorem~\ref{th:LRT} the chi-squared distribution \eqref{eq:TOTEMplex:ChiSquaredDistro} below Corollary~\ref{cor:ShpericalSymmetry}.
Again, we are interested in the weight assigned within $\equivclass$ to the empirical equivalence class $\equivclass[{\vecOp[M']}]$, so that 
\begin{proposition}[$I$-test]
\label{prop:Test}
    At significance level $\alpha$, we declare nested $\vecOp[M']$ to extract significantly new information from data 
    over constructing element $\vecOp$, when 
    \equ{
    1 - \textsc{cdf}_{\chi^2_{D'-D}}(2N\infdiv{\iprojSymbol[{\equivclass[{\vecOp[M']}][\prob f]}]\refP}{\iprojSymbol\refP}) < \alpha
    ~.
    }
\end{proposition}

\subsection{The (self-)consistency of large-$N$ expansion}

The developed formalism  needs neither to assume that a reference distribution is close to the true 
distribution underlying the generation of data,  nor relies on appropriate measurements to unveil the truth. 
In fact, it does not need to adhere to the very notion of some ground truth. 
In that sense, it is a purely effective description driven by data after deciding about  attributes and their domains. 
Propositions~\ref{prop:Score} and~\ref{prop:Test} simply indicate those measurements that inform us about systematic trends --\,beyond our prior knowledge\,-- as evidenced by the data at a given sample size. 

To further support this model- and truth-agnostic formulation  of operator selection at  large-$N$, we now adopt the more familiar viewpoint that  there exists a well-defined system --\,governed by some stochastic laws\,-- which generates data. 
In the following, we 
denote by ${\prob t}$ the target distribution that underlies data generation.
Furthermore, we signify by $\mathcal T_0$ the set of non-realizable  entities so that $\mathcal T_{\star} =\statespace\setminus\mathcal T_0$. 
Given a compatible ($\refP_\config=0$ means $\ket{\config}\in\mathcal T_0$) reference distribution ${\refP}\in\mathcal P$, we can always move from $\refP$ to $\prob t$ by applying the translation operator (c.f.\ \eqref{th:iProj:ExponentialForm}), 
\equ{
    \label{def:vTruth}
    \prob t = \exp\left\{\boldsymbol\tau \cdot \vecOp[T]\right\} \refP =
    \left(\proj_{\mathcal T_{\star}} + \boldsymbol\tau \cdot \vecOp[T] + \sfrac12 (\boldsymbol\tau \cdot \vecOp[T])^2 + \ldots\right)\refP
    ~.
}
Irreducible vector operator $\vecOp[T]\in\mathcal B(\mathcal T_{\star})$ describes the number of measurements $\Delta=\rank\vecOp[T]$ that are required to learn the target  on top of reference distribution $\refP\in\mathcal P$. 
$\boldsymbol\tau\in\mathbb R^{\Delta}$ encode displacement parameters fixed by the details of the underlying system. 

The rank $\Delta$ of the vector operator 
demonstrates the significance of reference distribution in inferring laws from data.
If we know beforehand or correctly guess the law, then $\Delta=1$ corresponding  to normalization on $\mathcal P(\mathcal T_{\star})$. 
A \totem over properly declared entities with $\statespace_0=\mathcal T_0$ (so that $\statespaceR=\mathcal T_{\star}$) is already sufficient after projecting prior knowledge onto $\equivclass[{\proj_{\statespaceR}}]$, $\iprojSymbol[{\equivclass[{\proj_{\statespaceR}}]}]\refP = e^{\tau} \proj_{\mathcal T_{\star}}\refP=\prob t$, i.e.\ the data provides no additional valuable information. 
On the other hand, if our prior knowledge fails to capture aspects of the underlying law governing subjects, $1<\Delta\leq\abss{\statespace}$.
In the poorest state of prior knowledge when $\Delta=\abss{\statespace}$, a \totem would need to infer the probability of every entity from measurements on the data (cf.\ expectations $\braketO{\config}{\prob f}{\config}$) whipping out all probabilities assumed by $\refP$. 

Ultimately, the efficiency of \mle estimate is assured by (proof in~\ref{app:MLE})

\begin{restatable}{theorem}{MLE}
\label{th:InformationCriterion}
    Among various vector operators $\vecOp\in\mathcal B(\statespaceR)$ including one that is \textsc{fapp}-equivalent to $\vecOp[T]$ from \eqref{def:vTruth}, we benchmark the constructed \totem{s} using information score~\ref{def:iScore}. 
    With increasing sample size, $\vecOp_\textsc{mle}$ is expected  to be \textsc{fapp}-equivalent to $\vecOp[T]$.
\end{restatable}

\section{Illustrative examples}
\label{sc:Examples}

\subsection{Logistic regression}

As a straight-forward application, consider a regression setting with $m$ predictors $x_i\in\mathcal D_{\texttt{X}_i}$ (summarized by vector $\mathbf x$) and one response variable $y\in\mathcal D_\texttt{Y}$. 
A \totem constructed by 
\equ{
\label{eq:LogisticRegression:ConstructingElement}
\vecOp[R] = \left\{\proj_\statespace, \proj_y, \proj_{\mathbf x}, \proj_{y,\mathbf x}
\right\}
}
learns the prevalence of response, all $m$ pairwise marginals $\vev{\prob f \proj_{y, x_i}}$ alongside the joint predictor distribution $\vev{\prob f \proj_{\mathbf x}}$ through linear system \eqref{eq:pheno:LinearMap}.
Assuming (to avoid regularizing as in~\ref{app:Regularization}) that none of the aforementioned marginals vanishes and declaring $\statespace_0=\emptyset$ unambiguously  leads to \maxent distribution
 according to the exponential form in Theorem~\ref{th:iProj:ExponentialForm}:
\equ{
\label{eq:LogisticRegression_MaxEnt}
\iprojSymbol[{\vecOp[R]}] \prob u = 
\prod_{i=1}^m \sum_{x_i\in\mathcal D_{\texttt X_i}} 
Z^{-1}
\sum_{y\in\mathcal D_{\texttt Y}}
\vev{\prob f\proj_{\mathbf x}}
\exp\left\{\theta_{\texttt Y}(y) + \sum_{i=1}^m\theta^i_{\texttt Y,\texttt X}(y,x_i)\right\}
\ket{y,\mathbf x}\bra{y,\mathbf x}
}
with logistic partition function
\equ{
Z \equiv Z(\mathbf x)
=
\sum_{y\in\mathcal D_\texttt{Y}}\exp\left\{\theta_{\texttt Y}(y) + \theta_i(y,x_i)\right\}
~.
}

As explicitly shown in~\ref{app:LogisticRegression} for the case of binary variables, when Lagrange multipliers are promoted to parameters, the \maxent distribution associated to constructing element \eqref{eq:LogisticRegression:ConstructingElement} is structurally identical to the model distribution of logistic regression --\,up to reparametrization symmetry.
Exclusively considering marginal constraints, 
we do not need to worry about such redundancies (which get enhanced for categorical variables with multiple levels) when running \textsc{ipf} algorithm \eqref{eq:IPF}. 
In fact, we are free to specify the (finite) domain of our variables.
Most  importantly, $D=\rank\vecOp[R]$ robustly informs us about the degrees of freedom in classifying regression architectures via Propositions~\ref{prop:Score} and~\ref{prop:Test}.

\subsection{Binomial distribution - revisited}

To illustrate the insights unveiled by \totem we consider the {coin experiment} under the logic of~\cite{jaynes1968prior}, where 
$N$ agents perform $L$ Bernoulli experiments and we obtain (portions of) data to decide about the coin(s) used.
The outcome of a Bernoulli experiment is unambiguously  described by a binary variable $s_i\in\mathcal D_i = \{\texttt{head}, \texttt{tail}\}$, so that 
 entity $\ket{s_1,\ldots,s_L}$ results from fully specifying the sequence $\boldsymbol s\equiv(s_1,\ldots,s_L)$ of outcomes from the Cartesian product $\mathcal S =\{\texttt{head}, \texttt{tail}\}^L$. 
%

Let us declare \texttt{head} as ``success'' and define\footnote{As usual, $\delta$ denotes the indicator map in the $i$-th domain $\mathcal D_i$.} a success operator 
\equ{
\label{eq:binomial:success_operator}
\op H \ket{\boldsymbol s} = \frac{l(\boldsymbol s)}{L} \ket{\boldsymbol s} \quads{\text{with}} l(\boldsymbol s) =  \sum_{i=1}^L \delta(s_i, \texttt{head})
~,
}
whose eigenvalues are individually given by the success rates of each Bernoulli realization.
Based on success operator \eqref{eq:binomial:success_operator}, average moments $\vev{\prob f \op H^n}$ for $n\in\mathbb N$ can be initiated.  
Compatible with this information is the marginal distribution over the number of successes $k$. 
Recording in the experiment only the frequency of number of \texttt{head} can be realized by a $k$-marginal operator with eigenvalue problem 
\equ{
\label{eq:binomial:k_marginalOperator}
\proj_k \ket{\boldsymbol s} = \delta_{l(\boldsymbol s), k} \ket{\boldsymbol s}
~.
}
Here, $\delta$ stands  for the Kronecker delta. 
$\proj_k $ is an indicator operator that projects onto the subspace of Bernoulli realizations with fixed number of \texttt{head}, i.e.\ success number $k$.
Specifying phenomenology $\phi_k\equiv\vev{\prob f \proj_k}$ for $k=0, \ldots, L$
partially breaks
permutation symmetry over the $2^L$ Bernoulli entities into degenerate subspaces with fixed $k$
whose dimension 
is given by 
binomial coefficient 
\equ{
\label{eq:BinomialCoefficient}
b_{L,k} \equiv
\begin{pmatrix}
L\\
k
\end{pmatrix}
~.
}

\paragraph*{The perception of one coin}
In the case that only success rate $\eta\equiv\vev{\prob f \hat H}$ is specified from data, there is a particularly simple form (derived in~\ref{app:Coins}) for the \maxent distribution,
\begin{align}
\label{eq:binomial:MaxEnt_globalMean}
\iprojSymbol[{\equivclass[{\vecOp[H]}]}] \prob u  
= 
\sum_{\boldsymbol s\in\mathcal S} \eta^{l(\boldsymbol s)} \left(1-\eta\right)^{L-l(\boldsymbol s)} \ket{\boldsymbol s}\bra{\boldsymbol s}
~,
\end{align}
associated to constructing element $\vecOp[H] = \{ \proj_\statespace, \op H\}$.
Eventually, it is straight-forward to take the expectation w.r.t.\ \maxent distribution 
of marginal operator \eqref{eq:binomial:k_marginalOperator}, 
\equ{
\vev{\iprojSymbol[{\equivclass[{\vecOp[H]}]}]\prob u\, \proj_k} = 
b_{L,k}
\,\eta^{k} \left(1-\eta\right)^{L-k}~.
}
This coincides with the probability estimate for event $k$ from the familiar binomial distribution $B(L,\eta)$.

Following the same steps, we can determine the \iproj of $\prob u$ onto the $k$-marginal \totemplex constructed by 
$\vecOp[K]=\{\proj_0, \ldots, \proj_L\}$, 
\equ{
\label{eq:binomial:k_marginals}
\iprojSymbol[{\equivclass[{\vecOp[K]}]}] \prob u 
=
\prod_{k=0}^L \exp\left\{ \log \frac{\phi_k}{b_{L,k}}\, \proj_k \right\} \proj_{\mathcal S}
=
\sum_{\boldsymbol s\in\mathcal S} \frac{\phi_{l(\boldsymbol s)}}{b_{L,l(\boldsymbol s)}} \ket{\boldsymbol s}\bra{\boldsymbol s}
~.
}
that reproduces the empirical relative frequencies $\phi_k\equiv\vev{\prob f\proj_k}$ of the possible success outcomes $k=0,\ldots,L$.
Normalization is implied by $\sum_{k=0}^L \phi_k=1$.
Incidentally, success operator \eqref{eq:binomial:success_operator} can be written as the superposition of $k$-marginal operators \eqref{eq:binomial:k_marginalOperator}:
\equ{
\label{eq:binomial:Nestedness}
\op H = L^{-1}\sum_{k=0}^L k\,\proj_k ~.
}
Hence, \totemplex $\equivclass[{\vecOp[K]}]\subseteq\equivclass[{\vecOp[H]}]$ is nested in the sense of Definition~\ref{def:NestedTOTEMs}.

Using the distribution over the quotient space $\equivclass[{\vecOp[H]}]\slash\vecOp[K]$ in Theorem~\ref{th:LRT}, we can compare how representative is the inner  \totemplex of $k$-marginals within the ambient  \totemplex mean success rate.
At large $N$, the $I$-test results into  a chi-squared distribution with $\rank\vecOp[K]-\rank\vecOp[H]=L-1$ degrees of freedom and  statistics dictated by $2N$ times 
\begin{align}
\label{eq:binomial:TestStatistics}
\infdiv{\iprojSymbol[{\equivclass[{\vecOp[K]}]}] \prob u}{\iprojSymbol[{\equivclass[{\vecOp[H]}]}] \prob u } 
= 
\sum_{k=0}^L 
\phi_k
\log
\frac{
\phi_k
}
{
b_{L,k}
\eta^k \left(1-\eta\right)^{L-k}
}
~.
\end{align}
As intuitively anticipated, the $I$-test statistics boils down to  the \textsc{kl} distance of the binomial distribution with empirical mean success rate $\eta=\vev{\prob f \op H}$ from the 
empirical marginal distribution $\phi_k=\vev{\prob f\proj_k}$ of success number $k$.

To perform the $I$-test,  we do not need to assume anything more concrete about the experimental setup or the empirical distribution, beyond specifying the entity space of $L$ Bernoulli trials. 
In particular, we have never assumed independence of those Bernoulli trials. 
Therefore, it is straight-forward to apply our \totem formalism to explore the significance of higher phenomenologies $\vev{\prob f\op H^n}$ for $n>1$ against the $k$-marginal \totem.

\paragraph*{Updating prior knowledge}

After deciding in favor of a \totem on \texttt{day1} constructed say by $\op H$, we can use the \iproj of $\prob u$ onto $\equivclass[{\op H}][{\prob f_\texttt{day1}}]$ as our updated reference distribution: 
\equ{
\refP_1 = \iprojSymbol[{\equivclass[{\vecOp[H]}][{\prob f_\texttt{day1}}]}] \prob u
~.
}
On the following day, the \iproj of $\refP_1$ onto the \totemplex constructed by $\op H$ on the new data boils down to simply updating average success rate $\eta_1=\vev{\prob f_\texttt{day1}\op H}$ in \eqref{eq:binomial:MaxEnt_globalMean} with the new estimate $\eta_2=\vev{\prob f_\texttt{day2}\op H}$. This is easily recognized by the form of \iproj{s} in Theorem~\ref{th:iProj:ExponentialForm},
\equ{
\iprojSymbol[{\equivclass[{\vecOp[H]}][{\prob f_\texttt{day2}}]}] \refP_1 = 
    \iprojSymbol[{\equivclass[{\vecOp[H]}][{\prob f_\texttt{day2}}]}] \iprojSymbol[{\equivclass[{\vecOp[H]}][{\prob f_\texttt{day1}}]}]  \prob u = 
\exp\{(\theta_1 + \theta_2)\op H - \log(Z_1Z_2)\proj_\statespace\} \prob u
~.
}
which leads to the same  structure \eqref{eq:Binomial:ExponentialForm_MaxEnt} by redefining Lagrange multipliers $\theta=\theta_1+\theta_2$ and $Z=Z_1Z_2$.  
Similarly, the \iproj of $\refP_1$ onto $\equivclass[{\vecOp[K]}][{\prob f_\texttt{day2}}]$ now learns the new $k$-marginals by heart,
since the form
\equ{
\iprojSymbol[{\equivclass[{\vecOp[K]}][{\prob f_\texttt{day2}}]}]\refP_1 
= \iprojSymbol[{\equivclass[{\vecOp[K]}][{\prob f_\texttt{day2}}]}] \iprojSymbol[{\equivclass[{\vecOp[H]}][{\prob f_\texttt{day1}}]}]\prob u
=
\exp\left\{ \sum_k \vartheta_k \proj_k + \theta_1 \op H - \log Z \right\} \prob u
~,
}
upon using \eqref{eq:binomial:Nestedness} and redefining $\tilde\vartheta_k=\vartheta_k + \theta_1 k / L - \log Z$, remains structurally the same leading to \eqref{eq:binomial:k_marginals} with $\phi_k = \vev{\prob f_\texttt{day2}\proj_{k}}$.
The test statistics of 
$
\infdiv{\iprojSymbol[{\equivclass[{\vecOp[K]}][{\prob f_\texttt{day2}}]}]\refP_1}{\iprojSymbol[{\equivclass[{\vecOp[H]}][{\prob f_\texttt{day2}}]}]\refP_1}
$
inform us whereas we need to update the mean success rate or not given the new $k$-marginals.

Furthermore, the covariance of Algorithm~\ref{th:NewtonRaphson} enables to numerically deduce the \iproj associated to any higher-degree moment $\op H^n$ 
beyond the mean. 
Consequently, we can decide about updating not only the value of the mean, but also the very structure of our prior knowledge from \texttt{day1}, augmenting it with the value of higher moments estimated on \texttt{day2}.

\paragraph*{The perception of two coins}
As a little twist, let us assume that there is a confounding attribute \texttt{Group} $g\in \mathcal D_\texttt{G} = \{\texttt{A,B}\}$ 
to which we suspect that slightly different coins are assigned.
Accordingly, the entity space is augmented by this binary variable, $(g;\boldsymbol s)\equiv(g, s_1,\ldots, s_L)\in \mathcal D_\texttt{G}\times \mathcal S$. 
In particular, success operator \eqref{eq:binomial:success_operator} is trivially lifted onto the augmented entity space, 
\equ{
\label{eq:TwoCoins:SuccessOperator}
\op H \ket{\config} = L^{-1} l(\boldsymbol s) \ket{\config} 
\quads{\text{where}}
\ket{\config} \equiv \ket{g;\boldsymbol s}
~.
}
To measure \texttt{Group} prevalence of team \texttt{A} within the entirety of $N$ agents, we use indicator operator
\equ{
\proj_\texttt{A} \ket{\config} = \delta(A, g)\ket{\config}
~,
}
and similarly for \texttt{Group} prevalence of team \texttt{B}.
Analogously to \eqref{eq:binomial:MaxEnt_globalMean}, we find (details in~\ref{app:Coins}) the \iproj of $\prob u$ onto a \totemplex constructed by\footnote{Since $\proj_\texttt{A}+\proj_\texttt{B}=\proj_\statespace$ (reflecting that $\Phi_\texttt{A}+\Phi_\texttt{B}=1$), $\{\proj_\statespace, \proj_\texttt{A}, \op H\}$ and $\{\proj_\texttt{A}, \proj_\texttt{B}, \op H\}$ are \textsc{fapp} equivalent.} $\vecOp[H]=\{\proj_\statespace, \proj_\texttt{A}, \op H\}$, 
\equ{
\label{eq:binomial:two_coins:oneMean_MaxEnt}
\iprojSymbol[{\equivclass[{\vecOp[H]}]}]\prob u
=
\sum_{g\in\{\texttt{A,B}\}}\Phi_\texttt{g} \sum_{\boldsymbol s\in\mathcal S} \eta^{l(\boldsymbol s)}(1-\eta)^{L-l(\boldsymbol s)} \ket{g;\boldsymbol s}\bra{g;\boldsymbol s}
}
to be a straight-forward generalization of preceding discussion on the global mean success rate additionally taking \texttt{Group} prevalence $\Phi_g\equiv\vev{\prob f\proj_g}$ for $g=\texttt{A,B}$ into account. 

The empirical success rate in \texttt{Group} $g$ can be estimated by the conditional mean
\equ{
\eta_{g}\equiv
\frac{\vev{\prob f \op H\proj_g}}{\vev{\prob f \proj_g}}
\quads{\text{for}} g=\texttt{A,B}
~,
}
using success operator \eqref{eq:TwoCoins:SuccessOperator}.
In addition to the phenomenology of global mean \eqref{eq:binomial:two_coins:oneMean_MaxEnt}, the \maxent distribution that incorporates knowledge from \texttt{Group} prevalences and \texttt{Group} mean success rates via linear system \eqref{eq:Pheno:LinearSystem},
\equ{
\vev{\prob p \proj_g} \overset{!}{=} \vev{\prob f \proj_g} = \Phi_g 
\qand
\vev{\prob p \op H \proj_g} \overset{!}{=} \vev{\prob f \op H \proj_g} = \Phi_g \eta_g 
\quads{\text{for}} g=\texttt{A,B}
}
is given by the \iproj 
\equ{
\label{eq:binomial:two_coins:twoMeans_MaxEnt}
\iprojSymbol[{\equivclass[{\vecOp[H']}]}] \prob u
=
\sum_{g\in\{\texttt{A,B}\}} \Phi_g \sum_{\boldsymbol s\in\mathcal S} \eta_g^{l(\boldsymbol s)} \left(1-\eta_g\right)^{L-l(\boldsymbol s)} \ket{g;\boldsymbol s}\bra{g;\boldsymbol s}
~.
}
By merit of 
$
\proj_\texttt{A} + \proj_\texttt{B} = \proj_\statespace
$
and 
$
\proj_\texttt{A}\op H + \proj_\texttt{B}\op H = \op H
$,
vector operator $\vecOp[H']=\{\proj_\texttt{A}, \proj_\texttt{B}, \op H\proj_\texttt{A}, \op H\proj_\texttt{B}\}$ (which treats the two teams $\texttt{A,B}$ symmetrically) is \textsc{fapp} equivalent to 
$\{\proj_\statespace, \proj_\texttt{A}, \op H, \op H\proj_\texttt{A}\}$.
Since the latter form readily implies $\vecOp[H]$,
we could immediately appeal to the $I$-test~\ref{prop:Test}.
At large $N$, we use chi-squared distribution with one degree of freedom and test statistics given by $2N \infdiv{\iprojSymbol[{\equivclass[{\vecOp[H']}]}] \prob u}{\iprojSymbol[{\equivclass[{\vecOp[H]}]}] \prob u}$  to compare whether $\eta_\texttt{A}$ is significantly different from $\eta_\texttt{B}$ or if both are well described by global $\eta$.

Again, our model-agnostic formalism with minimal assumptions allows the straight-forward comparison of higher phenomenologies $\vev{\prob f \op H^n \proj_g}$ constructing nested \totem{plices}. 
Whenever $n>1$ or multiple categorical attributes act as confounders quickly leading to analytically intractable expressions for the \iproj, the iterative minimization of $I$-divergence (entropy maximization in the case of $\refP=\prob u$) presented in~\ref{ssc:Idiv_Minimizer} turns out to be practically beneficial. 

\begin{appendix}

\section{The curse of undersampling}
\protect\label{app:Regularization}

Many realistic datasets suffer from sampling zeros which reduce the statistical power of any inferential method, usually whenever $N \leq \abss{\statespace}$. 
At the level of contingency tables and histograms (used to deduce empirical distribution \eqref{eq:Empirical:RelativeFrequencies}), sampling zeros are manifested by a wealth of vanishing entries\footnote{In particular, $ f_\config=0$ for any structural zero 
otherwise the provided dataset would have violated known deterministic relationships.}.
To counteract ``annoying'' undersampling effects, regularization has become a common practice in the literature~\cite{friedman2001elements}. 

Working in probability space, 
any sensible form of regularization boils down to the method~\cite{MorcosE1293} of \textit{pseudocounts}. 
The core idea relies on augmenting the data matrix 
$\texttt X$
by adding ``unobserved'' rows 
\equ{
\config\in\statespace\,:\quad
\texttt X_{\config,i} = e_i \quad\text{for }~ i=1,\ldots L
 ~.
}
This would amount to modifying e.g.\ in a uniform manner the empirical distribution $\prob f\rightarrow\prob f_\lambda$  such that $\vev{\prob f_\lambda}=1$ according to 
\equ{
\label{eq:RegularizedFrequencies}
{ f_\config} \quads{\rightarrow} \frac{{ f_\config} + \lambda/N}{1+\abss{\statespaceR}\cdot\lambda/N}
\quads{\text{for}} \config \in \statespaceR
~.
}
The introduced hyper-parameter $\lambda>0$, often dumped regularization strength 
helps interpolate from the original towards regularized data.

Despite that regularization methods lift sampling zeros and related singularity problems for parametrizations 
of certain phenomenologies, 
the --\,even uniform\,-- modification of input data $\texttt X$ in critical regions where $N$ and $\abss{\statespace}$ are of comparable order usually introduces artifacts 
that undermine the statistical methods presented in Section~\ref{sc:ClassificationOfKnowledge}.
Besides the extreme scenario of ``strong'' regularization 
which erases all data information in favor of full agnosticism, $\prob f_\lambda \overset{\lambda\rightarrow\infty}{\rightarrow}\prob u$, there exist less evident, serious situations.
Starting for concreteness from uniform prior, the test statistics of Corollary~\ref{cor:ShpericalSymmetry} are determined by the deviation of empirical entropy from the maximum entropy (entropy of the \iproj of $\prob u$). 
Whenever regularization \eqref{eq:RegularizedFrequencies} is applied,  
the critical entropy difference 
\equ{
\label{eq:regularized_entropy_difference}
\Delta H_c(\lambda) \equiv  H[\iprojSymbol[{\equivclass[][{\prob f_\lambda}]}]\prob u] -  H[\prob f_\lambda]
}
is distorted, especially if $\order{\abss{\statespaceR}} = \order{N}$. Unless $\lambda=0$, the distortion in that regime cannot be generically remedied 
even for ``mild'' regularization schemes where
\equ{
H[\prob f_\lambda] 
=H[\prob f] + \lambda \frac{\abss{\statespaceR}}{N} \left(1 - H[\prob f] \right)  + \order{\lambda^2}
~.
}
Most alarmingly, even if both entropies in \eqref{eq:regularized_entropy_difference} smoothly go to their non-regularized counterparts as $\lambda\rightarrow0$, the critical statistics governing \eqref{eq:ChiSqaureStatistics}, 
\equ{
2N  \Delta H_c(\lambda)  = 2N \Delta H_c(0) + \order{N \lambda} 
}
imply that any finite $\lambda$ could control (at fixed threshold $\alpha$ in Proposition~\ref{prop:Test}) alongside the physical scale $N$ the significance outcome in experimentally decisive settings where $\Delta H_c(0)$ is below $\order{1/N}$. 

To avoid such artifacts and 
stay data-centric retaining full control over the information flow in our framework, we do \textit{not} modify the data at any stage of the statistical analysis. Therefore, we must learn to tolerate (nearly) non-realizable entities. 

\section{Detailed derivations}
\label{app:Derivations}

\subsection{The information projection}
\label{app:iProjection}
\begin{theorem}
\Paste{th:iProjection}
\end{theorem}
\paragraph*{Proof}
In our data-centric approach, we recycle familiar~\cite{amari2000methods} notions of information geometry.
The uniqueness of a minimum of \idiv immediately follows in probability space by the convexity of the feasible region $\equivclass$ of the phenomenological problem \eqref{eq:Pheno:LinearSystem} at hand, combined with strict 
convexity of the \idiv~\cite{cover2012elements, joram_soch_2022_5820411} thought as a function $[0,1]^{\abss{\statespace}-\abss{\statespace_0}}\rightarrow\mathbb R_0^+$.
Hence,  the \idiv \eqref{eq:I_Divergence:DEF} possesses one global minimum in $\equivclass$, at most.

At the same time, a \totemplex in the identification of \eqref{eq:BasisVectors:ColumnVectorRepresentation} is represented  by a non-empty, convex, bounded and closed --\,hence compact\,-- subset of $[0,1]^{\abss{\statespace}-\abss{\statespace_0}}$, viz.\ \eqref{eq:TOTEM:LinearProgram}, over which any continuous function necessarily attains a minimum by the extreme value theorem. In total, we conclude that the \idiv must attain its global minimum in the \totemplex.

\subsection{Form of the information projection}
\label{app:ExponentialForm}

\iProjExponentialForm*

\paragraph*{Proof}
Following an established procedure in optimization theory,
we incorporate phenomenological
constraints on 
expectations 
\eqref{eq:Pheno:LinearSystem}
via Lagrange multipliers $\boldsymbol\theta$ in the \textit{Lagrangian}   
\equ{
\label{eq:iProj:Lagrangian}
\mathcal L 
= \infdiv{\prob p}{\refP}
- \sum_{\alpha=1}^{D} \theta_\alpha  \sum_{\config\in\cartprodR} m_\alpha(\config) \left(  p_\config -  f_\config \right) 
~.
}
Subsequently, we extremize the Lagrangian
in the probabilities $\prob p_\config>0$ and the parameters $\boldsymbol \theta\in\mathbb R^{D}$. 
%
Using exactly $D=\rank\vecOp$ Lagrange multipliers --\,in accordance with Definition~\ref{def:Irreducible_VectorOperator}\,-- ensures the independence of the incorporated constraints, which becomes essential for the convergence of many iterative algorithms and important for their efficiency.  

In the basis of admissible entities \eqref{eq:BasisVectors:ColumnVectorRepresentation}, we directly arrive at the exponential  expression
\equ{
\braketO{\config}{\iprojSymbol}{\config} = \refP_\config \exp \left\{\sum_{\alpha=1}^D \theta_\alpha m_\alpha(\config) - 1 \right\} 
}
for some vector of Lagrange multipliers such that \eqref{eq:Pheno:LinearSystem} is fulfilled.
Noting that our constructing elements~\ref{def:Irreducible_VectorOperator} always imply normalization \eqref{eq:Probabilities:Normalization} means that we can always transform via  \eqref{eq:TOTEMs:Reparametrization} performing a ``rotation'' in the space of expectations, so that  
\equ{
\exists\, T_{0\alpha} ~: \quad \sum_{\alpha=1}^D T_{0\alpha} m_\alpha(\config) = 1~.
}
Shifting accordingly
\equ{
\theta_\alpha \quads{\rightarrow} \theta_\alpha + T_{0\alpha}
}
identically results into the r.h.s.\ of
\eqref{eq:MaxEnt_parametricForm}. 

To show the converse of the Theorem we need a 
\begin{lemma}
\label{Lemma:Pythagoras}
    If  the relation  
    \equ{
    \label{eq:PythagorasFromLemma}
        \infdiv{\prob p}{\refP} = \infdiv{\prob q}{\refP} +  \infdiv{\prob p}{\prob q} 
        \quads{\text{for all}} {\prob p}\in\equivclass 
    }
    for a distribution ${\prob q}\in\equivclass$ holds,
    then this probability operator 
    is the \iproj of ${\refP}$ onto the \totemplex $\equivclass$. 
\end{lemma}
In fact, this is an analogue of 
Pythagorean theorem for $I$-divergences, see e.g.~\cite{csiszar2004information,nielsen2018information}. 
\paragraph*{Proof.}
Due to the non-negativity of the \idiv \eqref{eq:I_Divergence:DEF}, the forward direction of the Lemma follows from \eqref{eq:PythagorasFromLemma} rewritten as  
\equ{
\infdiv{\prob q}{\refP} = \infdiv{\prob p}{\refP} - \infdiv{\prob p}{\prob q}
\leq \infdiv{\prob p}{\refP},
}
for all distributions ${\prob p}\in\equivclass$,
which reproduces the definition of the \iproj in the \totemplex.
$\qed$\\

Back to the proof of the Theorem,
substituting the assumed exponential form into the r.h.s.\ of \eqref{eq:PythagorasFromLemma} we see that 
\equ{
\infbraket{\refP}{\prob q} + \infbraket{\prob q}{\prob p} =
\infbraket{\refP}{\prob p} + \sum_{\alpha=1}^D\theta_\alpha\underbrace{\sum_{\config\in\statespaceR} \left( q_\config -  p_\config \right) m_\alpha(\config)}_{=0\text{ in }[\prob f;\vecOp]}
}
on any probability operator ${\prob p}\in\equivclass$.
Hence, the variationally optimal solution \eqref{eq:MaxEnt_parametricForm} fulfills the relation of Lemma~\ref{Lemma:Pythagoras} and is the unique \iproj of ${\refP}$ onto $\equivclass$.

\subsection{Constrained minimization of information divergence}
\label{app:Newton}

\NewtonRaphson*

\paragraph*{Proof}
At the $n$-th step,  the $D\times D$ Jacobian $\mathbf J$ is invertible --\,which is a necessary condition for convergence\,-- due to  following
\begin{lemma}
\label{lm:Jacobian}
    Let $\mathbf P$ be a real symmetric positive definite matrix (signified by $\mathbf P\succ0$) 
    and $\mathbf C$ a real matrix of full-row rank with a column-structure compatible with $\mathbf P$. 
    Then, the adjoint map from the column to the row space of $\mathbf C$ is also given by 
    a symmetric positive definite and thus invertible matrix
    \equ{
        \mathbf C \mathbf P \mathbf C^T
        °.
    }
\end{lemma}
\paragraph*{Proof}
For linearly independent rows in $\mathbf C$, it is easy to show that the adjoint map is positive definite (and thus invertible):
\begin{align}
\forall\,\mathbf u \in \mathbb R^{\rank\mathbf C}\,:
0 \overset{!}{=} &\,\,
\mathbf u^T \left(\mathbf C \mathbf P \mathbf C^T\right)\mathbf u
=
\left(\mathbf C^T\mathbf u\right)^T \mathbf P \left(\mathbf C^T\mathbf u\right)
\\[1ex]
&\overset{\mathbf P\succ0}{\Longrightarrow} \quad
~
\mathbf u^T\mathbf C = \mathbf 0_{\rank\mathbf P}
\quads{\overset{\mathbf C\succ0}{\Longrightarrow}} 
\mathbf u = \mathbf 0_{\rank\mathbf C}
~.\phantom{...}\qed
\nonumber
\end{align}

Since an \textit{irreducible} vector operator~\ref{def:Irreducible_VectorOperator} that measures collective characteristics  constructs a \totem, the eigenvalues $m_\alpha(\config)$ on admissible  microstates seen as a linear map \eqref{eq:pheno:LinearMap} constitute a matrix of full row-rank equaling $D$.
A diagonal operator with elements $\prob p^{(n)}_e>0$ can be thought of in the representation of \eqref{eq:BasisVectors:ColumnVectorRepresentation} as a diagonal matrix 
in $\mathbb R^{\abss{\statespace}-\abss{\statespace_0}}$ that remains positive definite throughout the optimization. 
Consequently, Lemma~\ref{lm:Jacobian} applies and the Jacobian is indeed at each Newton-Raphson step  positive-definite and thus invertible.

\subsection{The information divergence in nested descriptions}
\label{app:nestedTOTEMs}

\nestedNewtons*
\noindent
To assert that successive determination of \iproj{s} on nested \totem{plices} is equivalent to the \iproj of referrence distribution to the innermost, target \totemplex, we need 
\begin{lemma}
For nested \totem{lices} $\equivclass[{\vecOp[B]}]\subseteq\equivclass[{\vecOp[A]}]$ spanned by irreducible vector operators  ${\vecOp[A]}$ and ${\vecOp[B]}$  
the \iproj of ${\refP}$ onto the inner \totemplex $\equivclass[{\vecOp[B]}]$ satisfies chain rule
\equ{
\label{eq:iProj:ChainRule}
\iprojSymbol[{\equivclass[{\vecOp[B]}]}]\refP
=
 \iprojSymbol[{\equivclass[{\vecOp[B]}]}] \iprojSymbol[{\equivclass[{\vecOp[A]}]}] \refP
}
coinciding with the \iproj onto the inner \totemplex of reference distribution $\iprojSymbol[{\equivclass[{\vecOp[A]}]}]\refP$ which is in turn
the \iproj   of ${\refP}$ onto the ambient \totemplex $\equivclass[{\vecOp[A]}]$.
\end{lemma}
\noindent
Reading \eqref{eq:iProj:ChainRule} from right to left reveals that the theorem can help eliminate the intermediate \iproj of ${\refP}$ onto the ambient \totemplex.
Repeated application of this notion starting at the terminating step $\nu=d$ of \eqref{eq:Newton:SubRoutines} to eliminate intermediate \iproj{s} 
from inner towards outer \totem{plices},
\begin{align}
\refP^{(d)} =&\,\,
\iprojSymbol[{\equivclass[{\vecOp[M]^{(d)}}]}]  \iprojSymbol[{\equivclass[{\vecOp[M]^{(d-1)}}]}] \iprojSymbol[{\equivclass[{\vecOp[M]^{(d-2)}}]}]  \cdots \iprojSymbol[{\equivclass[{\vecOp[M]^{(1)}}]}] \refP
\\[1ex]
=&\,\,
\iprojSymbol[{\equivclass[{\vecOp[M]^{(d)}}]}]  \iprojSymbol[{\equivclass[{\vecOp[M]^{(d-2)}}]}] \cdots \iprojSymbol[{\equivclass[{\vecOp[M]^{(1)}}]}] \refP
= \ldots = 
\iprojSymbol[{\equivclass[{\vecOp[M]^{(d)}}]}]  \refP
~,
\nonumber
\end{align}
we 
verify that the suggested pipeline is equivalent to directly determining the \iproj of ${\refP}$ onto the \totemplex defined by $\vecOp^{(d)}\equiv\vecOp$.

%

\paragraph*{Proof}
We start from the definition of the \iproj of $\iprojSymbol[{\equivclass[{\vecOp[A]}]}]\refP$ onto the inner \totemplex $[\prob f; \vecOp[B]]$:
\equ{
\label{eq:Newton:nestedTOTEMs:Proof}
\forall\,\,{\prob p}\in [\prob f; \vecOp[B]]\, : \quad 
\infdiv{\prob p}{\iprojSymbol[{\equivclass[{\vecOp[A]}]}]\refP}
\geq 
\infdiv{\iprojSymbol[{\equivclass[{\vecOp[B]}]}] \iprojSymbol[{\equivclass[{\vecOp[A]}]}]\refP }{ \iprojSymbol[{\equivclass[{\vecOp[A]}]}]\refP }
~.
}
Applying Corollary~\ref{cor:Pythagoras} in the ambient \totemplex, onto which $\iprojSymbol[{\equivclass[{\vecOp[A]}]}]\refP$ is the \iproj of ${\refP}$, for distributions ${\prob p}~,~\iprojSymbol[{\equivclass[{\vecOp[B]}]}]\iprojSymbol[{\equivclass[{\vecOp[A]}]}]\refP \in[\prob f]_{\hat{\mathbf A}}$ 
\begin{align}
\infdiv{\prob p}{ \iprojSymbol[{\equivclass[{\vecOp[A]}]}]\refP } =&\,\, 
\infdiv{\prob p}{\refP} - 
\infdiv{ \iprojSymbol[{\equivclass[{\vecOp[A]}]}]\refP }{\refP}
\quad \text{and}
\nonumber
\\[1ex]
\infdiv{\iprojSymbol[{\equivclass[{\vecOp[B]}]}]\iprojSymbol[{\equivclass[{\vecOp[A]}]}]\refP}{\iprojSymbol[\refP][{[\prob f]_{\hat{\mathbf A}}}]}
=&\,\,
\infdiv{\iprojSymbol[{\equivclass[{\vecOp[B]}]}]\iprojSymbol[{\equivclass[{\vecOp[A]}]}]\refP}{\refP} - \infdiv{\iprojSymbol[{\equivclass[{\vecOp[A]}]}]\refP}{\refP}
\nonumber
~,
\end{align}
cancels the common term $\infdiv{\iprojSymbol[{\equivclass[{\vecOp[A]}]}]\refP}{\refP}$ on both sides of defining inequality \eqref{eq:Newton:nestedTOTEMs:Proof}
leading to 
\equ{
\infdiv{\prob p}{\refP} 
\geq 
\infdiv{\iprojSymbol[{\equivclass[{\vecOp[B]}]}]\iprojSymbol[{\equivclass[{\vecOp[A]}]}]\refP}{\refP} 
~.
}
Since this applies for all ${\prob p}$ in the inner \totemplex $\equivclass[{\vecOp[B]}]$, we conclude by uniqueness (Theorem~\ref{th:iProjection}) that the distribution  $\iprojSymbol[{\equivclass[{\vecOp[B]}]}]\iprojSymbol[{\equivclass[{\vecOp[A]}]}]\refP$ defined as the \iproj of $\iprojSymbol[{\equivclass[{\vecOp[A]}]}]\refP$ onto $\equivclass[{\vecOp[B]}]$ is also the \iproj of ${\refP}$ onto $\equivclass[{\vecOp[B]}]$.
\phantom{...}$\qed$

\subsection{The Gaussian kernel}\protect\label{app:UniversalKernel}

\UniversalKernel*
\paragraph*{Proof}
First, we pass from perception of entities in the \totemplex to fluctuations of perception
which appropriately live in the kernel of $\vecOp$, as described above \eqref{eq:Fluctuations:PolyhedralCone}. 
Next, we perform 
a Taylor expansion of the \idiv of generic ${\prob p}\in\equivclass$ from the \iproj in the \totemplex around base distribution ${\prob p}=\iprojSymbol\refP\equiv \prob q$
in ``small'' fluctuations, 
\begin{align}
\label{eq:Idiv:TaylorExpansion1}
\infdiv{\prob p}{\prob q} = &\,\,
\vev{\left(\prob q + \probflac\right)\left(\log\left(\prob q + \probflac\right) - \log\prob q\right)}
\\[1ex]
=&\,\,
\vev{\left(\prob q + \probflac\right)\log\left(1+\prob q^{-1}\probflac\right)}
=
\vev{\probflac} + \sfrac12 \vev{\probflac^2q^{-1}} +  \order{\probflac^3}~.
\nonumber
\end{align}
The linear term vanishes, since we immediately recognize from \eqref{eq:Fluctuations:PolyhedralCone} that $\vev{\prob p}=1$ implies $\vev{\probflac}=0$.
This is expected, as we expand around the minimum of the \idiv in the convex space of the \totemplex.
As indicated in \eqref{eq:KnowledgeDistro:UniversalKernel}, there thus remains to leading order in the fluctuations 
a quadratic form 
\equ{
\label{eq:GaussianKernel}
\mathcal K(\probflac) = 
\vev{\probflac \prob q^{-1} \probflac} 
=
\sum_{a,b=1}^{\abss{\ker\vecOp}}\pi_a \vev{\op N_a \prob q^{-1} \op N_b} \pi_b
~.
}

\subsection{Spherical symmetry in probability fluctuations}
\label{app:shpericalSymmetry}

\SphericalSymmetry*
\paragraph*{Proof}
Any distribution on the \totemplex is sampled from ${\refP}$ as a fluctuation away from its \iproj 
by the  probability 
\equ{
\label{eq:KnowledgeDistro:ProbabilityToSample}
f(\probflac) \prod_{a=1}^{k}\dd \pi_a~,
}
in terms of an $N$-leading Gaussian viz.\ \eqref{eq:KnowledgeFluctuations:UniversalGaussian}.
Let $\lambda_a>0$ denote the $a$-th eigenvalue of the quadratic form in \eqref{eq:GaussianKernel} 
and group its orthonormal eigenvectors into a $k\times k$ orthogonal matrix $\mathbf U$.
Setting 
\equ{
\pi_a =  \sum_{b=1}^{k}U_{ab}  x_b 
\nonumber
}
we diagonalize the Gaussian in \eqref{eq:GaussianKernel} by the canonical fluctuations $\boldsymbol x\in\mathbb R^{k}$, 
\equ{
\mathcal K(\boldsymbol x) = \sum_{a=1}^k \frac{ x_a^2}{\lambda_a}~,
\nonumber
}
with unit-$\det$ Jacobian. 
Subsequently we eliminate the positive eigenvalues by setting
\equ{
 x_a = \sqrt{\lambda_a}\, y_a 
\quads{\text{so that}}
\mathcal K(\boldsymbol y) = \sum_{a=1}^k  y_a^2~,
\nonumber
}
where the Jacobian $\sqrt{\prod_a \lambda_a}=\order{1}$ of the total differential in \eqref{eq:KnowledgeDistro:ProbabilityToSample} gets canceled  by the partition function $Z$. 
Eventually, switching to spherical coordinates with unit-$\det$ Jacobian via
\equ{
 r^2 = \sum_{a=1}^k  y_a^2
\nonumber
}
such that the probability 
becomes with Gaussian function $f( r)\equiv f(\probflac)$
\equ{
\label{eq:KnowledgeFluctuations:ShpericalMeasure}
\dd\Omega_k\, \dd r\,  r^{k-1} f( r)~,
}
the spherical symmetry of the $N$-leading Gaussian term readily manifests,
\equ{
\mathcal K( r) =  r^2
~.
\nonumber
}
It is trivial to execute the angular integration, where the solid angle in $k$ dimensions equals 
\equ{
\Omega_k = \frac{2\pi^{k/2}}{\Gamma(k/2)}
~.
}

Finally, 
we express \eqref{eq:KnowledgeFluctuations:ShpericalMeasure} in terms of the quadratic form,
\equ{
\label{eq:KnowledgeFluctuations:RadialProbability_K}
\Omega_k\, \dd r \,  r^{k-1} \, f( r) = \Omega_k\, \dd \mathcal K
\underbrace{\mathcal K^{k/2-1}\, f(\mathcal K)}_{\rho(\mathcal K)}
}
using the squared radius in kernel space as integration variable, $\dd( r^2/2)=\dd\mathcal K$.  
From \eqref{eq:KnowledgeDistro:UniversalKernel}, we know that
\equ{
\mathcal K(\boldsymbol r) = \infbraket{\prob q}{\prob p} + \order{N^{-1/2}}
}
so that the r.h.s.\ of probability \eqref{eq:KnowledgeFluctuations:RadialProbability_K} can be re-expressed in terms of $\infbraket{\prob q}{\prob p}$, with
the density $\rho$ taking the desired form of \eqref{eq:ToTEMplex:spherical_ProbabilityDensity} up to $N$-suppressed corrections. 

\subsection{The probability of equivalence classes}
\label{app:LRT}

\LRT*

\paragraph*{Proof}
As outlined in \eqref{eq:KernelFluctuations:Expansion}, 
a distribution 
on inner \totemplex  
$\equivclass[{\vecOp[M']}][\prob p]\subset\equivclass$
can be parameterized through probability fluctuations 
around the \iproj  
of reference  distribution onto the ambient \totemplex $\equivclass$.  
According to Theorem~\ref{th:Universality:GaussianApproximation},  a Gaussian density always approximates  at sufficiently large-$N$ the  
distribution of such probability fluctuations around $\iprojSymbol\refP$  on the ambient \totemplex.

Next, we formulate 
\begin{lemma}
\label{lm:Quotient:KnowledgeDistribution}
    Given the ambient \totemplex $\equivclass$, the probability at large $N$  to sample from reference $\refP$ any distribution on the inner \totemplex $\equivclass[{\vecOp[M']}][\prob p]$ 
    is determined to $N$-leading order by $ \infdiv{\iprojSymbol[{\equivclass[{\vecOp[M']}][\prob p]}]\refP}{\iprojSymbol\refP}$. 
    Parametrizing this \idiv by probability fluctuations in the normal space complementing $\ker\vecOp[M']$ to $\ker\vecOp$, the associated probability density      
    follows
    a multivariate Gaussian distribution
    with 
    $D'-D$ degrees of freedom. 
\end{lemma}
\paragraph*{Proof}
Using the definition of quotient space, 
the quadratic form of Lemma~\ref{lm:UniversalKernel} that governs the distribution of distributions \eqref{eq:TOTEMplex:KnowledgeDistribution2} on the ambient \totemplex decomposes with $\prob q\equiv \iprojSymbol\refP$ according to
\begin{align}
&\,\,2\infbraket{\prob q}{\prob p} = 2\vev{\probflac\prob q^{-1}\probflac} 
\\
=&\,\,
\sum_{a',b'=1}^{\abss{\ker\vecOp'}}\pi_{a'}\pi_{b'} \vev{\op N_{a'}\prob q^{-1}\op N_{b'}} + \sum_{a,b=1}^{D'-D}\pi_{a}\pi_{b} \vev{\op N_{a}\prob q^{-1}\op N_{b}} 
+ 
2\sum_{a,a'}\pi_{a'}\vev{\op N_{a'}\prob q^{-1}\op N_{a}} \pi_{a}
\nonumber
\end{align}
for any $\prob p\in\equivclass$.
Integrating \eqref{eq:TOTEMplex:KnowledgeDistribution2} over all fluctuations  $\pi_{a'}$ with positive-definite precision matrix $\mathcal K'_{a'b'}= \vev{\op N_{a'}\prob q^{-1}\op N_{b'}}$ 
we can calculate 
the probability density to sample from reference ${\refP}$ any distribution ${\prob p'}$ on the inner \totemplex. 
In the resulting marginal, 
there still remains a Gaussian in $\abss{\ker\vecOp}- \abss{\ker\vecOp'}=D'-D$ dimensions defined by quadratic form
\begin{align}
&2\infdiv{\iprojSymbol[{\equivclass[{\vecOp[M']}][\prob p]}]\refP}{\prob q} 
= 
\min\limits_{\prob p'\in\equivclass[{\vecOp[M']}][\prob p]} 2\infdiv{\prob p'}{\prob q} 
\\
& =
\sum_{a,b=1}^{D'-D} \pi_{a} \left(\vev{\op N_a\prob q^{-1}\op N_b} - \sum_{a',b'=1}^{\abss{\ker\vecOp'}}\vev{\op N_{a'}\prob q^{-1}\op N_{a}}(\mathcal K')^{-1}_{a'b'} \vev{\op N_{b'}\prob q^{-1}\op N_{b}}\right) \pi_{b} 
~,
\nonumber
\end{align}
see Appendix~\ref{apps:MarginalGaussian} for technical details.

Descending from a positive-definite quadratic form $\mathcal K$ on the ambient \totemplex (due to Lemma~\ref{lm:UniversalKernel}),  the resulting quadratic form 
is also positive-definite. 
In total, we find that in the quotient space $\equivclass\slash\vecOp[M']$ 
the conjectured density \eqref{eq:Quotient:KnowledgeDistribution}
that governs the probability of inner \totemplex $\equivclass[\vecOp'][\prob p]$ is self-consistently approximated by a Gaussian distribution in the normal fluctuations.
\phantom{...}$\qed$

Along the lines of Corollary~\ref{cor:ShpericalSymmetry}, we can transform to spherical coordinates in order to integrate-out the angular dependence. After  rescaling the radius squared to achieve $Q=2N \infdiv{\iprojSymbol[{\equivclass[{\vecOp[M']}][\prob p]}]\refP}{\prob q} $, we can eventually express the radial probability to select from reference ${\refP}$ any \totemplex from the quotient space with \idiv in the vicinity of $Q/2N$ by a chi-squared distribution in $D'-D$ dimensions.   

\subsection{The efficiency of MLE}\protect\label{app:MLE}

\MLE*

\paragraph*{Proof}
In the convenient scenario that we know the target distribution, we exemplify the efficiency of 
operator selection. 
Starting from some reference distribution  ${\refP}$ 
we concretely verify that 
information score \eqref{eq:MLE:InformationScore} is expected  to be maximized for \totem{s} constructed by vector operators that are \textsc{fapp}-equivalent to $\vecOp[T]$. 
To do so we need to take expectations 
over (infinitely) many datasets.

In a hypothetical setting, 
datasets $\mathbf X$ of fixed size $N$ are generated from target distribution ${\prob t}$. Each dataset is described by empirical distributions ${\prob f}\in\mathcal P$, as outlined in \eqref{eq:Empirical:RelativeFrequencies}. 
Over sampled phenomenologies $\braketO{\statespaceR}{\prob f \vecOp}{\statespaceR}$, we take the expectation of information score \eqref{eq:MLE:InformationScore} assigned to \totem{s} 
constructed by $\vecOp$, %
\equ{
\label{eq:Expected_iScore}
\sum_{N\prob f} \multidistro\left(N\prob f;\prob t\right) \text{score}_I 
=
-N \sum_{N\prob f} \multidistro\left(N\prob f;\prob t\right) \infbraket{\iprojSymbol\refP}{\prob f} + \frac{\abss{\ker\vecOp}}{2} \log N + \order{1}
~,
}
where the sum runs over all count vectors $N\prob f\in\mathbb N_0^{\abss{\statespace}}$ that are partitions of $N$.
At this level, recall that incompatible 
prior knowledge which falsely declares nullentities from $\statespace\setminus\mathcal T_0$ would eventually lead \eqref{eq:MLE:InformationScore} to diverge to $-\infty$,
see also discussion below Definition~\ref{def:BaseKnowledge}.
Consequently, we only need to consider \totem{s} with nullentities from some $\statespace_0\subseteq\mathcal T_0$.

Henceforth, we focus on the $N$-leading contribution to the $I$-score.
Recycling the argumentation of Section~\ref{ssc:DicreteContinuous}, which becomes meaningful at larger $N$, we can write an integral representation for the expectation:
\equ{
\label{eq:iScore:Expectation_IntegralRepresentation}
Z^{-1}\int_{\prob f\in\mathcal P} \exp\left\{- N \infbraket{\prob t}{\prob f}\right\}  \infbraket{\iprojSymbol\refP}{\prob f} 
~.
}
%
%
In each \totemplex $\equivclass[][\prob q]$ from the quotient space $\mathcal P\slash\vecOp$ labeled by $\vecOp$-distinct knowledge  $\prob q$, we next apply Lemma~\ref{Lemma:Pythagoras} from the true perspective,
\equ{
\infbraket{\prob t}{\prob f} =  
\infbraket{\iprojSymbol[{\equivclass[][\prob q]}]\prob t}{\prob f} + \infbraket{\prob t}{ \iprojSymbol[{\equivclass[][\prob q]}]\prob t }
}
for sampled empirical $\prob f\in\equivclass[][\prob q]$. 
This observation covariantly breaks the integral \eqref{eq:iScore:Expectation_IntegralRepresentation} over the simplex into 
\equ{
\label{eq:integrationOverDatasets}
\int_{\prob q^*} \exp\left\{-N\infbraket{\prob t}{\iprojSymbol[{\equivclass[][\prob q]}]\prob t}\right\} 
\int_{\prob f\in\equivclass[][\prob q]} \exp\left\{-N\infbraket{\iprojSymbol[{\equivclass[][\prob q]}]\prob t}{\prob f}\right\}
 \infbraket{\iprojSymbol[{\equivclass[][\prob q]}]\refP}{\prob f}
~.
}
The first integral of the density in quotient space runs over all $\vecOp$-distinct distributions $\prob q$
in the sense of equivalence relation \eqref{eq:QuotientSpace:EquivalenceRelation}, while 
the second integral of \totemplex density runs within $\ker\vecOp$ over distributions ${\prob f}\in\equivclass[\vecOp][\prob q][\mathcal S_0]$, as in Section~\ref{ssc:DistroOverDistros}. Appropriate normalization factors are understood for both integrals.

Applying localization techniques on the inner integral, the saddle point of the exponential function is readily given by the \iproj $\iprojSymbol[{\equivclass[][\prob q]}]\prob t$ of target distribution onto the \totemplex. Considering probability fluctuations around this distribution
gives a Gaussian integral in $\ker\vecOp$  (cf.\ \eqref{eq:GaussianKernel}), which is well-defined at large $N$.
On the other hand, the $I$-score is calculated using prior knowledge $\refP$ with an a priori different \iproj on the same \totemplex.
Generically, we thus anticipate the integral representation \eqref{eq:iScore:Expectation_IntegralRepresentation} to scale as $\order{1}$, so that the  expected score \eqref{eq:Expected_iScore} scales linearly with negative sample size. 

However, if $\vecOp$ acting over $\statespace\setminus\mathcal T_0$ captures all measurements of $\vecOp[T]$, then $ \infbraket{\iprojSymbol[{\equivclass[][\prob q]}]\refP}{\prob f}= \infbraket{\iprojSymbol[{\equivclass[][\prob q]}]\prob t}{\prob f}$ 
on account of 
\begin{lemma}
\label{lm:NestedTruth}
    Consider a \totem over $\statespace\setminus\mathcal T_0$ constructed by vector operator 
    $\vecOp$ describing $D$ independent measurements
    which imply
    in the sense of Definition~\ref{def:NestedTOTEMs} the $\refP$-efficient vector operator $\vecOp[T]$ that parameterizes the target distribution,
    \equ{
    \label{eq:nestedTruth}
        \hat T_\alpha = \sum_{\alpha'=1}^D S_{\alpha\alpha'} \hat M_{\alpha'}\quads{\text{for}} \alpha=1,\ldots,D_{min}~.
    }
    The \iproj of reference distribution $\refP$ onto the \totemplex $\equivclass[][\prob q]$ constructed by $\vecOp$ coincides with the \iproj of true knowledge $\prob t$ onto $\equivclass[][\prob q]$. 
\end{lemma}
\paragraph*{Proof}
For the rest of this proof, we silently imply the set $\mathcal T_0$ of truly non-realizable entities.
We start from the definition  of the \iproj of ${\prob t}$ onto $\equivclass[][\prob q]$:
\begin{align}
\forall~\prob p\in \equivclass[][\prob q] \,: \quad
\infbraket{\prob t}{\prob p}
\geq 
\infbraket{\prob t}{\iprojSymbol[{\equivclass[][\prob q]}]\prob t}
\quad &\Leftrightarrow
\\
\infbraket{\refP}{\prob p}
\geq 
\infbraket{\refP}{\iprojSymbol[{\equivclass[][\prob q]}]\prob t}
+ \sum_{\alpha=1}^{D_{min}} \tau_{\alpha} \sum_{\alpha'=1}^D S_{\alpha\alpha'}  
\underbrace{
\vev{\left(\prob p - \iprojSymbol[{\equivclass[][\prob q]}]\refP\right)\op M_{a'}}
}_{
=0 
}
\nonumber
~.
\end{align}
The expectation difference of all $\op M_{\alpha'}$ vanishes, since both $\prob p$ and $\iprojSymbol[{\equivclass[][\prob q]}]$ belong to $\equivclass[][\prob q]$.
Hence,  
$\iprojSymbol[{\equivclass[][\prob q]}]\prob t$
is also the \iproj of $\refP$ onto $\equivclass[][\prob q]$ and since the \iproj of a reference distribution onto the \totemplex is by Theorem~\ref{th:iProjection} unique, we conclude that
$\iprojSymbol[{\equivclass[][\prob q]}]\prob t=\iprojSymbol[{\equivclass[][\prob q]}]\refP$. 
\phantom{...}$\qed$

This Lemma helps evaluate the second integral in \eqref{eq:integrationOverDatasets} within $\ker\vecOp$,
\equ{
\label{eq:Expected_iScore:Inner_KernelIntegral}
\int_{\prob f\in\equivclass[][\prob q]} \exp\left\{-N\infbraket{\iprojSymbol[{\equivclass[][\prob q]}]\prob t}{\prob f}\right\}
 N\infbraket{\iprojSymbol[{\equivclass[][\prob q]}]\refP}{\prob f} = \frac{\abss{\ker\vecOp}}{2} = \order{1}
}
using \eqref{eq:appGaussian:QuadraticExpectation}, whenever $\vecOp$ satisfies \eqref{eq:nestedTruth}.

In total, we conclude for \eqref{eq:Expected_iScore} that 
\equ{
\sum_{N\prob f} \multidistro\left(N\prob f;\prob t\right) \text{score}_I = 
\begin{cases}
   \frac{\abss{\ker\vecOp}}{2} \log N + \order{1} & \vecOp \text{ implies } \vecOp[T] \text{ over } \statespace\setminus\mathcal T_0\\[1ex]
   -\order{N} & \text{otherwise}
\end{cases}
~.
}
Among all the $\vecOp$ that 
imply vector operator 
$\vecOp[T]$ parametrizing the target knowledge, clearly the operator that is \textsc{fapp} equivalent to $\vecOp[T]$ exhibits the largest kernel and hence maximizes our information score according to \mle principle.

\subsection{The expected divergence from target}\protect\label{app:ExpectedTruth}
To provide further evidence in favor of the minimization of \idiv, 
we next evaluate the expected \idiv from target distribution ${\prob t}$ of the \iproj of reference ${\refP}$ onto the \totemplex $\equivclass[][\prob f]$ constructed on empirical $\prob f$ which was sampled from $\prob t$. 
Similar to \eqref{eq:integrationOverDatasets}, the integral representation over datasets breaks under the quotient action into integration over equivalence classes 
and integration within \totem{plices}.
Trivially executing  the inner integral within $\ker\vecOp$ like in \eqref{eq:Expected_iScore:Inner_KernelIntegral}, 
we are left this time with an integral over $\vecOp$-distinct distributions,
\equ{
\label{eq:ExpectedDivergenceTruth:QuotientSpace_Integral}
N \sum_{N\prob f}\infbraket{\prob t}{\iprojSymbol[{\equivclass[][\prob f]}]\refP}
=
Z^{-1}\int_{\prob q^*} \exp\left\{-N\infbraket{\prob t}{\iprojSymbol[{\equivclass[][\prob q]}]\prob t}\right\} N \infbraket{\prob t}{\iprojSymbol[{\equivclass[][\prob q]}]\refP} 
+ \ldots
~.
}

Saddle-point techniques indicate that the integral localizes around $\iprojSymbol[{\equivclass[][\prob q]}]\prob t$.
As noted  in the proof of Theorem~\ref{th:InformationCriterion}, whenever $\mathcal S_0\subset\mathcal T_0$ or if $\vecOp$ is not sufficient to capture all the measurements of $\vecOp[T]$, the expectation of the \idiv based on ${\refP}$ remains non-vanishing at large $N$.
Hence, we only need to focus on vector operators $\vecOp$ acting over $\statespace\setminus\mathcal T_0$ that satisfy \eqref{eq:nestedTruth}.
On such \totem{plices}, the \iproj{s} of reference and target distribution coincide due to Lemma~\ref{lm:NestedTruth}.

Introducing a displaced  target operator 
\equ{
    \hat q = 
    \exp\left\{\sum_{\alpha=1}^{D-1}\theta_\alpha\op M_\alpha\right\}\prob t
}
via linearly independent operators $\op M_\alpha$ that do not imply the identity $\proj_\statespace$, we can parameterise 
the \iproj of ${\prob t}$ onto \totem{plices} from the quotient space $\mathcal P(\statespace\setminus\mathcal T_0)\slash\vecOp$ according to 
\equ{
\iprojSymbol[{\equivclass[][\prob q]}]\refP =
\iprojSymbol[{\equivclass[][\prob q]}] \prob t = 
\vev{\hat q}^{-1} \hat q
~.
}
Recall that expectations are taken w.r.t.\ $\ket{\statespaceR}=\sum_{\config\in\cartprodR}\ket{\config}$. Expanding in ``small'' fluctuations, 
\equ{
\label{eq:QuotientSpace:Parametrization}
\iprojSymbol[{\equivclass[][\prob q]}] \prob t =
\left[1 + \sum_{\alpha=1}^{D-1}\left(\op M_\alpha - \vev{ \prob t \op M_\alpha} \right)\theta_\alpha + \order{\boldsymbol\theta^2}  \right] \prob t
~,
}
gives  rise to 
\equ{
\label{eq:QuotientSpace:iDivergence}
\infbraket{\iprojSymbol[{\equivclass[][\prob q]}] \prob t}{\prob t} = 
\sfrac12 \vartheta_\alpha 
\vev{(\op M_\alpha - \vev{ \prob t \op M_\alpha} )\prob t(\op M_\beta - \vev{ \prob t \op M_\beta})}
\vartheta_\beta
+ \order{\boldsymbol\vartheta^2} 
= \infbraket{\prob t}{\iprojSymbol[{\equivclass[][\prob q]}] \prob t}
~.
}
Due to vanishing linear term, the Taylor expansion up to quadratic order symmetrizes  the \idiv.

The measuring operators $\op M_\alpha$ are per assumption linearly independent, while $c\proj_\statespace = \sum_{a=1}^{D-1}c_\alpha\op M_\alpha$ implies $c=c_\alpha=0$. 
Hence, the shifted operators $\op M_\alpha - \vev{ \prob t \op M_\alpha}= \op M_\alpha - \vev{ \prob t \op M_\alpha}\proj_\statespace$ remain linearly independent and can be summarized by an irreducible vector operator in the sense of~\ref{def:Irreducible_VectorOperator}.  
Applying Lemma~\ref{lm:Jacobian}, we conclude that the precision matrix in \eqref{eq:QuotientSpace:iDivergence} is positive definite, 
following the argumentation in the proof of Theorem~\ref{th:LRT}.
Along similar lines, the Gaussian approximation in the small-fluctuations expansion \eqref{eq:QuotientSpace:Parametrization} becomes a large-$N$ expansion around the saddle point of \eqref{eq:ExpectedDivergenceTruth:QuotientSpace_Integral}, so that we can evaluate the expectation by \eqref{eq:appGaussian:QuadraticExpectation}. In total, we have 
\begin{align}
&\sum_{N\prob f} N\infbraket{\iprojSymbol\refP}{\prob t}  
\\
\approx& \sum_{N\prob f} N \infbraket{\prob t}{\iprojSymbol\refP}  
=\begin{cases}
   \frac{D-1}{2}+ \order{N^{-1/2}} & \vecOp \text{ implies } \vecOp[T] \text{ over } \statespace\setminus\mathcal T_0\\[1ex]
   \order{N} & \text{otherwise}
\end{cases}
~.
\nonumber
\end{align}

For all operators that imply $\vecOp[T]$ the expectation depends on the number $D$ of measurements that construct the associated phenomenology. For any other operator, the expectation diverges with $N$. Therefore, the \iproj of ${\refP}$ onto the \totemplex constructed by $\vecOp[T]$ over $\statespace\setminus\mathcal T_0$ exhibits the lowest expected \idiv from the target distribution. 
In particular, whenever ${\refP}={\prob t}$ the expected \idiv vanishes asymptotically. 
This extends the results of~\cite{loukas2022entropy} beyond uniform prior and categorical attributes.

\subsection{Logistic regression}
\protect\label{app:LogisticRegression}

To relate to the formal definition of logistic regression we assume for simplicity that both response and predictors are binary, i.e.\ $\statespace=\{\texttt{TRUE}, \texttt{FALSE}\}^{m+1}$. 
Since the constructing element \eqref{eq:LogisticRegression:ConstructingElement} is reducible, 
the \maxent distribution \eqref{eq:LogisticRegression_MaxEnt}
\equ{
\iprojSymbol[{\equivclass[{\vecOp[R]}]}]\prob u 
= \sum_{y,\mathbf x} \exp\{H(y,\mathbf x)\}
\ket{y,\mathbf x}\bra{y, \mathbf x}
}
with energy function
\equ{
\label{eq:LogisticRegression:Energy}
H(y, \mathbf x) = \theta_0 + \theta_\texttt{Y}(y) + \sum_{i=1}^m \theta^i_{\texttt Y,\texttt X}(y,x_i) + \Theta(\mathbf x) 
}
enjoys gauge symmetry that leaves $H$ invariant in Lagrange multipliers:
\begin{align}
\label{eq:LogisticRegression:GaugeTrafo}
\theta_0~{\rightarrow}~& \tilde\theta_0 = \theta_0 + \gamma_0
\quads{,}
\theta_\texttt{Y}(y)~{\rightarrow}~\tilde\theta_\texttt{Y}(y)=\theta_\texttt{Y}(y) - \gamma_0 + \sum_{i=1}^m\gamma^i_{\texttt Y}(y)~,
\nonumber
\\
\theta^i_{\texttt Y,\texttt X}(y,x_i) ~{\rightarrow}~& \tilde\theta^i_{\texttt Y,\texttt X}(y,x_i)=\theta^i_{\texttt Y,\texttt X}(y,x_i) - \gamma^i_{\texttt Y}(y) - \gamma^i_{\texttt X}(x_i)
\\
\qand
\Theta(\mathbf x)~{\rightarrow}~& \tilde\Theta(\mathbf x)=\Theta(\mathbf x)  + \sum_{i=1}^m \gamma^i_{\texttt X}(x_i)
~.
\nonumber
\end{align}
In the gauge where we take 
\begin{align*}
\gamma_0 = \theta_\texttt{Y}(0) + \sum_{i=1}^m \left( \theta^i_{\texttt Y,\texttt X}(0,0) -\gamma^i_\texttt{X}(0)  \right) 
\quads{,}&
\gamma^i_\texttt{X}(1) = \gamma^i_\texttt{X}(0) + \theta^i_{\texttt Y,\texttt X}(0,1) - \theta^i_{\texttt Y,\texttt X}(0,0)
\\
\gamma^i_\texttt{Y}(0) = -\gamma^i_\texttt{X}(0) + \theta^i_{\texttt Y,\texttt X}(0,0) 
\quads{,}&
\gamma^i_\texttt{Y}(1) = -\gamma^i_\texttt{X}(0) + \theta^i_{\texttt Y,\texttt X}(1,0) 
\\
\text{with}\quad&
\sum_{i=1}^m\gamma^i_\texttt{X}(0) = -\Theta(\mathbf 0) \qand \gamma^i_\texttt{X}(0) \in \mathbb R
~,
\end{align*}
most transformed Lagrange multipliers on the r.h.s.\ of \eqref{eq:LogisticRegression:GaugeTrafo} vanish identically:
\equ{
\tilde\theta_\texttt{Y}(0) = 0 \quads{,} \tilde\theta^i_{\texttt Y,\texttt X}(Y,0) = 0 \quads{,}  \tilde\theta^i_{\texttt Y,\texttt X}(0,X_i) = 0
\qand
\tilde\Theta(\mathbf 0)=0 ~.
}
Subsequently, we count 
\equ{
\rank\vecOp[R] = 1 + 1 + m + 2^m-1 = 1+ m + 2^m
}
remaining Lagrange multipliers encoding independent marginal constraints.

In such gauge, we can equivalently write Lagrange multipliers in terms of regression coefficients $\beta_0\in\mathbb R$ and $\boldsymbol\beta\in\mathbb R^m$ as
\equ{
   \tilde\theta_0 = -\log Z \quads{,}
   \tilde\theta_{\texttt Y}(y) = \beta_0 s_0
    \qand
    \tilde\theta^i_{\texttt Y, \texttt X}(y, x_i) =  \beta_i s_0 s_i 
}
via encoding
\equ{
\label{eq:LogisticRegression:OneHot_Encoding}
s_0(y=\text{False}) = 0~,~ s_0(y=\text{True}) = 1
\qand
s_i(x_i=\text{False}) = 0 ~,~ s_i(x_i=\text{True}) = 1
~.
\nonumber
}
This parametrizes the entity space by a binary vector $\ket{s_0;\mathbf s}$ where $s_0\in\{0,1\}$ and $\mathbf s\in \{0,1\}^m$.
Our energy function \eqref{eq:LogisticRegression:Energy} becomes
\equ{
H(s_0,\mathbf s) = \left(\beta_0 + \boldsymbol\beta^T \mathbf s\right)s_0 + \Theta(\mathbf s) -\log Z~.
}
Using the phenomenological constraint on the joint  predictor probabilities 
\equ{
\vev{(\prob p - \prob f)\proj_\mathbf{s}} = 0
}
we can eventually bring the \maxent distribution \eqref{eq:LogisticRegression_MaxEnt} to the more familiar form 
\equ{
\iprojSymbol[{\equivclass[{\vecOp[R]}]}]\prob u =
\sum_{s_0,\mathbf s} \vev{\prob f\proj_\mathbf{s}} \frac{\exp\left\{\left(\beta_0 + \boldsymbol\beta^T\mathbf s\right)s_0\right\}}{1+\exp\left\{\beta_0 + \boldsymbol\beta^T\mathbf s\right\}} \ket{s_0,\mathbf s}\bra{s_0,\mathbf s}
~.
}
where multiple summation over the encoded \eqref{eq:LogisticRegression:OneHot_Encoding} entity space is abbreviated by 
\equ{
\sum_{s_0,\mathbf s} \equiv \sum_{s_0\in\{0,1\}} \prod_{i=1}^m \sum_{s_i\in\{0,1\}}
}

In particular, the conditional probability invoked in the definition of the odds for response $s_0=1$ given predictor profile $\mathbf s$ can be directly read off:
\equ{
   \frac{\braketO{y,\mathbf s}{\iprojSymbol[{\equivclass[{\vecOp[R]}]}]\prob u}{y, \mathbf s}}{\vev{\prob f\proj_\mathbf{s}}} = %
   \frac{\exp\left\{\left(\beta_0 + \boldsymbol\beta^T\mathbf s\right)s_0\right\}}{1+\exp\left\{\beta_0 + \boldsymbol\beta^T\mathbf s\right\}}
   ~.
}
The remaining non-vanishing Lagrange multipliers $\beta_0$ and $\boldsymbol\beta$ are fixed by the empirical estimates of one- and two-site marginals 
\equ{
\vev{(\iprojSymbol[{\equivclass[{\vecOp[R]}]}]\prob u - \prob f)\proj_{s_0}} = 0 \qand \vev{(\iprojSymbol[{\equivclass[{\vecOp[R]}]}]\prob u - \prob f)\proj_{s_0,\mathbf{s}}} = 0~,
}
respectively. Due to redundancies in the marginal description, it suffices to consider e.g.\ $s_0=1$ and $s_i=1$, so that  (using the eigenvalues of projector operators $\delta(s_0,1)=s_0$ and $\delta(s_i,1)=s_1$)
\begin{align}
\sum_{\mathbf s} \frac{\vev{\prob f \proj_{\mathbf s}}}{1+ \exp\{-\beta_0-\boldsymbol\beta^T\mathbf s\}} \overset{!}{=}&\, \vev{\prob f\proj_{s_0=1}}\quad \text{and}
\\
\sum_{\mathbf s} \frac{\vev{\prob f \proj_{\mathbf s}} s_i}{1+ \exp\{-\beta_0-\boldsymbol\beta^T\mathbf s\}} \overset{!}{=}&\, \vev{\prob f\proj_{s_0=1,s_i=1}}
\nonumber
~.
\end{align}
In the  model-centric picture, the same equations arise by maximizing the likelihood 
\equ{
\ell(\prob f; \iprojSymbol[{\equivclass[{\vecOp[R]}]}]\prob u) = -H[\prob f] - \infdiv{\prob f}{\iprojSymbol[{\equivclass[{\vecOp[R]}]}]\prob u}
}
(equivalently minimizing the \idiv of empirical from model distribution) as a function of the regression coefficients, see also discussion in~\ref{app:ModelCentric}.

\subsection{Coin experiments}
\protect\label{app:Coins}

\paragraph*{Perception of one coin}
According to Theorem~\ref{th:iProj:ExponentialForm}, we anticipate 
\equ{
\label{eq:Binomial:ExponentialForm_MaxEnt}
\iprojSymbol[{\equivclass[{\vecOp[H]}]}]\prob u = Z^{-1}\exp\{ \theta \op H \} \prob u ~,
}
for the \iproj of $\prob u$ onto the \totemplex $\equivclass[{\vecOp[H]}]$ constructed by $\vecOp[H] = \{ \proj_\statespace, \op H\}$.
The partition function implementing normalization (from $\proj_\statespace$) evaluates to a closed form, 
\begin{align}
Z = \vev{\exp \{ \theta \op H \} \prob u} 
= 2^{-L}\prod_{i=1}^L\sum_{s_i\in\mathcal D_i} \exp\{ \theta l(\boldsymbol s)/L\}  
=&\,\,
2^{-L}\sum_{l=0}^L b_{L,l}
\exp\{ \theta l/L\}
\nonumber
\\
=&\,\,
2^{-L}\left(1 + \exp\{\theta/L\}\right)^L
~,
\end{align}
Substituting into the phenomenological constraint of mean success rate $\eta$,  we can solve for the Lagrange multiplier,
\equ{
\label{eq:binomial:global_mean:Derivation}
\vev{\iprojSymbol[{\equivclass[{\vecOp[H]}]}]\prob u \op H} = \frac{1}{1 + e^{-\theta/L}} \overset{!}{=}  \vev{\prob f \op H} \equiv \eta
\quads{\Rightarrow}
\theta = L\log\frac{\eta}{1-\eta}~,
}
so that in total we arrive at (recall that the resolution of identity is given by \eqref{eq:IdentityOperator})
\begin{align}
\label{eq:binomial:MaxEnt_globalMean:Deriviation}
\iprojSymbol[{\equivclass[{\vecOp[H]}]}] \prob u = &\,\,
\exp\left\{ L \log\eta \op H + L \log(1-\eta) (\hat I - \hat H)\right\} \proj_{\mathcal S}
\\[1ex]
= &\,\,
\sum_{\boldsymbol s\in\mathcal S} \eta^{l(\boldsymbol s)} \left(1-\eta\right)^{L-l(\boldsymbol s)} \ket{\boldsymbol s}\bra{\boldsymbol s}
~.
\nonumber
\end{align}
This is the least biased distribution given knowledge of entity space and the empirical success rate $\eta=\vev{\prob f \hat H}$.

Whenever the frequencies of success numbers are provided, we analogously find  the \maxent distribution 
\equ{
\label{eq:binomial:k_marginals:Derivation}
\iprojSymbol[{\equivclass[{\vecOp[K]}]}] \prob u 
=
\exp\left\{\sum_{k=0}^L\vartheta_k \proj_k \right\}\prob u 
=
\prod_{k=0}^L \exp\left\{ \log \frac{\phi_k}{b_{L,k}}\, \proj_k \right\} \proj_{\mathcal S}
=
\sum_{\boldsymbol s\in\mathcal S} \frac{\phi_{l(\boldsymbol s)}}{b_{L,l(\boldsymbol s)}} \ket{\boldsymbol s}\bra{\boldsymbol s}
}
where Lagrange multipliers are fixed according to
\equ{
\vartheta_k = \log\frac{2^L\phi_k}{b_{L,k}}
~.
}
Ultimately, we can test the ``binomial hypothesis'', i.e.\ whereas the mean success rate is sufficient to describe the distribution of success rates, using statistics \eqref{eq:binomial:TestStatistics}
or in operator language: 
\begin{align}
\label{eq:binomial:Idiv_derivation}
\infdiv{\iprojSymbol[{\equivclass[{\vecOp[K]}]}] \prob u}{\iprojSymbol[{\equivclass[{\vecOp[H]}]}] \prob u } =
\\
\left\langle 
  \exp\left\{\sum_{k=0}^L\log\frac{\phi_k}{b_{L,k}}\proj_k\right\}
 \left[\sum_{k'=0}^L\log \frac{\phi_{k'}}{b_{L,k'}} \proj_{k'} - L \log\eta \op H - L \log(1-\eta) (\proj_\statespace - \op H)  
\right]
\right\rangle
\nonumber
\\
=
\sum_{k=0}^L 
\phi_k
\log
\frac{
\phi_k
}
{
b_{L,k}
\eta^k \left(1-\eta\right)^{L-k}
}
~.
\nonumber
\end{align} 
To arrive at \eqref{eq:binomial:Idiv_derivation}, we do not need to assume independence of Bernoulli experiments. 
$k$-marginal distribution $\vev{\prob f \proj_k}$ can have any form, even evidencing  extreme violation of independence between Bernoulli experiments, which would be quickly detected by the $I$-test. 

To analytically investigate the sensitivity of \totem description when we know the target $\prob t$ beforehand, 
assume that there exists a mild violation of the independence between $i_0$-th and $j_0$-th Bernoulli experiments quantified via the two-point correlator
\equ{
\kappa\equiv \vev{\prob t\proj_{s_{i_0}=\texttt{head}}\proj_{s_{j_0}=\texttt{head}}} - \vev{\prob t\proj_{s_{i_0}=\texttt{head}}} \vev{\prob t\proj_{s_{j_0}=\texttt{head}}}
~.
}
By encoding 
\equ{
\sigma_i(s_i=\texttt{head})=1 \qand \sigma_i(s_i=\texttt{tail})=0 \quads{\text{for}} i=1,\ldots L
}
 Bernoulli configurations can be for example generated through an Ising-like system~\cite{ising1925beitrag} described by target distribution
\equ{
\prob t = Z^{-1}\sum_{\boldsymbol \sigma\in\{0,1\}^L}\exp\left\{J\sigma_{i_0} \sigma_{j_0} + h\sum_{i=1}^L\sigma_{i}\right\} \ket{\boldsymbol \sigma}\bra{\boldsymbol \sigma}
}
with partition function
\equ{
Z = \left(1+e^{h}\right)^{L-2}\left(1+2e^{h}+e^{J+2h}\right)~.
}
This implies ``free'' success rates
\equ{
\vev{\prob t \proj_{s_i=\texttt{head}}} = \frac{1}{1+e^{-h}} 
\quads{\text{for}} i\neq i_0,j_0~,
}
as before ($J=0\Leftrightarrow\kappa=0$ recovers the ``free'' theory of independent Bernoulli trials) and in addition ``interacting'' success rates
\equ{
\vev{\prob t \proj_{s_{i_0}=\texttt{head}}} = \vev{\prob t \proj_{s_{j_0}=\texttt{head}}} = \frac{e^{h}\left(1+e^{h+J}\right)}{1+2e^{h}+e^{J+2h}}
~.
}
In particular, the $k$-marginals estimate to
\equ{
\vev{\prob t\proj_k} = Z^{-1} \left[ b_{L-2,k-2} e^{h k+J} + b_{L-2,k}e^{hk} + 2 b_{L-2,k-1}e^{hk} \right]~.
}
The Ising parameters (coupling $J$ and external field $h$) are given by the inverse problem, perturbatively:
\equ{
h = \log\frac{\eta}{1-\eta} - \frac{2\kappa}{L\eta(1-\eta)^2} + \order{\kappa^2}
\qand
J = \frac{\kappa}{\eta^2(1-\eta)^2} + \order{\kappa^2}~.
}

In such scenario, the expected \idiv \eqref{eq:binomial:Idiv_derivation}, which at large $N$ is dominated by the saddle point (cf.
  \eqref{eq:iScore:Expectation_IntegralRepresentation}) around the true distribution, 
would scale 
\begin{align}
\int_{\prob f} \exp\{-N\infdiv{\prob f}{\prob t}\} \infdiv{\iprojSymbol[{\equivclass[{\vecOp[K]}]}] \prob u}{\iprojSymbol[{\equivclass[{\vecOp[H]}]}] \prob u } =&\,\,
\infdiv{\iprojSymbol[{\equivclass[{\vecOp[K]}][{\prob t}]}] \prob u}{\iprojSymbol[{\equivclass[{\vecOp[H]}][{\prob t}]}] \prob u }
+ \ldots
\nonumber
\\
=&\,\,
\frac{1}{L(L-1)}\frac{\kappa^2}{\eta^2(1-\eta)^2} + \ldots
~.
\end{align}
Consequently, the test statistics in~\ref{prop:Test} is expected to lead to significant conclusions when 
\equ{
N \gg L(L-1) \frac{\eta^2(1-\eta)^2}{2\kappa^2}~.
}

\paragraph*{Perception of two coins}

The \iproj \eqref{eq:binomial:two_coins:oneMean_MaxEnt} of $\prob u$ onto a \totemplex constructed by $\vecOp[H]=\{\proj_\statespace, \proj_\texttt{A}, \op H\}$ can be analytically expressed as
\begin{align*}
\iprojSymbol[{\equivclass[{\vecOp[H]}]}]\prob u = &\,\,
\exp\left\{ \log\Phi_\texttt{A}\,{\proj_\texttt{A}} + \log\Phi_\texttt{B}\, {\proj_\texttt{B}}
+
L \log \eta \, \op H + L \log (1-\eta) (\proj_\statespace-\op H)\right\}\Id_\statespace
\\[1ex]
 =&\,\,
\sum_{\config\in\cartprod} \Phi_\texttt{A}^{\delta(g, \texttt{A})} \Phi_\texttt{B}^{\delta(g, \texttt{B})} \eta^{l(\boldsymbol s)}(1-\eta)^{L-l(\boldsymbol s)} \ket{\config}\bra{\config}
\\[1ex]
 =&\,\,
\sum_{g\in\{\texttt{A,B}\}}\Phi_\texttt{g} \sum_{\boldsymbol s\in\mathcal S} \eta^{l(s)}(1-\eta)^{L-l(\boldsymbol s)} \ket{g;\boldsymbol s}\bra{g;\boldsymbol s}
 \nonumber
\end{align*}
in terms of the global mean success rate $\eta$ and \texttt{Group} prevalences $\Phi_g\equiv\vev{\prob f\proj_g}$ for $g=\texttt{A,B}$.

The \iproj \eqref{eq:binomial:two_coins:twoMeans_MaxEnt} of $\prob u$ onto the \totemplex $\equivclass[{\vecOp[H']}]$ constructed by  $\vecOp[H']=\{\proj_\texttt{A}, \proj_\texttt{B}, \op H\proj_\texttt{A}, \op H\proj_\texttt{B}\}$ (\textsc{fapp} equivalently by $\{\proj_\statespace, \proj_\texttt{A}, \op H, \op H\proj_\texttt{A}\}$)
can also be analytically written as 
\begin{align}
\iprojSymbol[{\equivclass[{\vecOp[H']}]}] \prob u
=&\,\,
\prod_{g\in\{\texttt{A,B}}\exp \left\{ \log\Phi_g \proj_g  + L \log\eta_g \op H\proj_g  + L\log(1-\eta_g) (\proj_\statespace-\op H) \proj_g \right\} \proj_{\statespace}
\nonumber
\\[1ex]
=&\,\,
\sum_{g\in\{\texttt{A,B}\}} \Phi_g \sum_{\boldsymbol s\in\mathcal S} \eta_g^{l(\boldsymbol s)} \left(1-\eta_g\right)^{L-l(\boldsymbol s)} \ket{g;\boldsymbol s}\bra{g;\boldsymbol s}
~,
\end{align}
in terms of \texttt{Group} prevalences and \texttt{Group} mean success rates.

\section{Relations to model-centric literature}
\label{app:ModelCentric}

\paragraph*{The BIC score}
The $I$-score of Definition~\ref{def:iScore} can be obtained in a model-centric approach, as well, which is in a sense perpendicular to the data-centric derivation previously discussed.

Conventionally, \eqref{eq:LikelihoodForm:DEF} 
is assumed 
to be a function of some parameters reflecting structural assumptions --\,called a \textit{model} $\mathcal M$\,-- about the underlying problem. 
Notice that 
the likelihood \eqref{eq:likelihood} can be thought of as the leading model-dependent term governing the large-$N$ expansion \eqref{eq:MutlinomialDistro:LargeN} of
\equ{
\label{eq:BIC:loglikelihood}
\log\multidistro\left(N\prob f;\prob p\right) = N \left[\ell(\prob f;\prob f) - \ell(\prob f;\prob p)\right] + \order{\log N}
}
A priori, lacking information about microscopic interactions, there is no appealing reason to commit to any structural assumption. 
In any case, fully working with some parameterized structure is only sensible in absence of any nullentities, 
otherwise singularities would plague the parametrization eventually leading to its break-down. 
Hence, we can only --\,without regularizing\,-- investigate  problems with  $\statespace_0=\emptyset$ in the following of this paragraph.

According to the insight previously gained, 
the model distribution could be expressed via a parameterized\footnote{In a slight abuse of notation, we let $\boldsymbol\theta$ denote both Lagrange multipliers in Section~\ref{ssc:MIDIV} that are fixed by the given data as well as generic parameters in the present context that can vary independently of the data.} distribution
\equ{
\prob p(\boldsymbol\theta) 
=
\exp\left\{\sum_{a=1}^{D}\theta_\alpha \op M_\alpha\right\}\refP
=
\sum_{\config\in\cartprod}v_\config \exp\left\{\sum_{a=1}^{D}\theta_\alpha m_\alpha(\config)\right\} \ket{\config}\bra{\config}
~,
} 
in terms of the $D\times\abss{\statespace}$ \textit{architecture} matrix $m$ (of full row-rank) whose components would correspond to the eigenvalues in \eqref{eq:Pheno:Measurement} of our measuring operators. 
In other words, model knowledge parametrizes $\vecOp$-disinct distributions over the quotient space $\mathcal P\slash\vecOp$, 
cf.\ discussion below \eqref{eq:QuotientSpace:EquivalenceRelation}.
In the model-centric approach, $D-1$ is often called the model dimension.

In a Bayesian setting under uniform prior (leading-order derivation should not depend on any sensible prior~\cite{gelman2013bayesian}), we represent the \textit{model evidence} through a truncated expansion of \eqref{eq:BIC:loglikelihood} 
\equ{
P(\texttt X\vert \mathcal M) = Z^{-1}\int \dd^{D} \boldsymbol \theta\, \delta\left(N\vev{\prob p(\boldsymbol\theta)} - N\right) \exp\left\{-N\ell(\prob f;\prob p(\boldsymbol\theta))\right\} 
~.
}
This can be seen as a large-$N$ approximation (viz.\ \eqref{eq:FromDiscreteToContinuous})  of the cumulative probability to sample the observed counts $N f_\config$ starting from all possible parametric realizations $\boldsymbol\theta\in\mathbb R^{D}$ of the considered model structure. 
Again, $Z$ is the partition function which self-consistently normalizes $P$ in the large-$N$ expansion.
Using the integral representation of Dirac-delta function, we can rewrite  the model evidence as
\equ{
\label{eq:BIC:ModelEvidence}
P(\texttt X\vert \mathcal M) = Z^{-1} \int\frac{\dd \theta_0}{2\pi}\int \dd^{D} \boldsymbol \theta\, e^{N\mathcal L}
}
in terms of Lagrangian (constants in $\boldsymbol\theta$ can be absorbed by $Z$)
\equ{
\mathcal L = \sum_{\alpha=1}^{D} \theta_\alpha \sum_{\config\in\statespace} m_\alpha(\config) f_\config 
+ i\theta_0\left( \sum_{\config\in\statespace} p_\config(\boldsymbol\theta) - 1\right)
~.
}

Due to the $\order{N}$ scaling of the exponent in the integrand, Laplace approximation applies, as in the proof of Theorem~\ref{th:Universality:GaussianApproximation}. The saddle point equations are
\equ{
\label{eq:BIC:SaddlePoint}
\sum_{\config\in\cartprod}m_\alpha(\config)\left[ f_\config - i\theta_0 \, p_\config(\boldsymbol\theta)\right] \overset{!}{=} 0
\quads{\text{for}} \alpha=1,\ldots,D
\qand
\sum_{\config\in\cartprod} p_\config(\boldsymbol\theta) - 1 \overset{!}{=} 0 
~.
}
Since normalization on $\mathcal P$ is always implied by the eigenvalues $m_\config(\config)$, it follows from the first set of equations given normalization of empirical ${\prob f}\in\mathcal P$ that 
\equ{
i\theta_0 \vev{\prob p(\boldsymbol\theta)} 
\overset{!}{=} 
\vev{\prob f} = 1~,
}
which combined with the latter equation implies $\overline\theta_0=-i$. 
Hence, \eqref{eq:BIC:SaddlePoint} reduces to 
\equ{
\sum_{\config\in\statespace}m_\alpha(\config)\left[  f_\config -  p_\config(\overline{\boldsymbol\theta})\right] \overset{!}{=} 0
}
which reproduces the phenomenological system defining a \totem in~\ref{def:TOTEM}. Thus, the Maximum likelihood Estimator (\mle) at $\overline{\boldsymbol\theta}$ coincides with the \iproj of ${\refP}$ onto the corresponding \totemplex, $\prob p(\overline{\boldsymbol\theta}) = \iprojSymbol\refP$. 

At the same time, the Hessian of the Lagrangian computes to 
\equ{
\alpha,\beta=1,\ldots,D~:\quad
\frac{\partial^2\mathcal L}{\partial\theta_\alpha\partial\theta_\beta} = 
\vev{\op M_\alpha \prob p(\boldsymbol\theta) \op M_\beta}
\quads{,}
\frac{\partial^2\mathcal L}{\partial \theta_0\partial\theta_\alpha} = i \vev{\prob p(\boldsymbol\theta) \op M_\alpha}
\qand
\frac{\partial^2\mathcal L}{\partial \theta_0^2} = 0
~.
}
At the saddle point, 
\equ{
\overline{\mathcal L} = -N \ell(\prob f;\prob q) \quads{\text{with}} \prob q\equiv \iprojSymbol\refP 
~,
}
the determinant of the Hessian breaks into block determinants 
\equ{
\det
\left( \vev{\vecOp \prob q \mathbf{\hat M}^T } \right)
\cdot 
\det\left[
\left(i\vev{\prob q\mathbf M}\right)^\dagger
\left(\vev{\vecOp \prob q \mathbf{\hat M}^T }\right)^{-1}
\left(i\vev{\prob q\mathbf M}\right)
\right]
~.
}
Transposition and hermitian conjugation act on the indices of vector operator~\ref{def:Irreducible_VectorOperator}. 
Recall from \eqref{eq:iProj:RightDerivative} that the \iproj does not introduce any zero probabilities, if not dictated by the phenomenology. 
Due to Lemma~\ref{lm:Jacobian}, $\vev{\vecOp \prob q \mathbf{\hat M}^T }$ is positive definite. Furthermore, the second determinant simply reduces to the quadratic form of the inverse of $\vev{\vecOp \prob q \mathbf{\hat M}^T }$ (recall that the inverse of a positive definite matrix is also positive definite).
The determinants of the blocks thus remain positive 
asserting the consistency of the \mle at the saddle point \eqref{eq:BIC:SaddlePoint}. 

Executing eventually the integration in \eqref{eq:BIC:ModelEvidence} we arrive at 
\equ{
\log P\left(\texttt X;\mathcal M\right) = 
-\sfrac12\textsc{bic}  + \order{1}
}
in terms of the Bayesian information criterion~\cite{10.1214/aos/1176344136}
\equ{
\label{eq:BIC}
\textsc{bic} \equiv \left(D-1\right) \log N - 2 \overline{\mathcal L}  
~.
}
In absence of sampling zeros when parametric models are expected to behave smoothly,
$\statespaceR = \statespace$ at least for all considered constructing elements $\vecOp$. 
In such setting, maximizing score \eqref{eq:MLE:InformationScore} given some data $\texttt X$ becomes equivalent to minimizing the \textsc{bic} score \eqref{eq:BIC}.
Whether a universally accepted way exists to compute $\order{1}$ corrections to the \textsc{bic} 
and whether the equivalence of the model- and data-centric pictures continuous remains open.

\paragraph*{The likelihood ratio}
The  \idiv governing the distribution over nested \totem{plices} in the quotient space $\equivclass\slash\vecOp[M']$ can be viewed as determining a likelihood ratio of simpler $\mathcal M$ to more complex $\mathcal M'$ model, $D<D'$.
For that, we consider 
some inner \totemplex labeled by $\prob p\in\equivclass$
and a test distribution $\prob q\in\equivclass[{\vecOp[M']}][\prob p]\subset\equivclass$.
Starting from \eqref{eq:BIC:loglikelihood} evaluated at the saddle point \eqref{eq:BIC:SaddlePoint} where the \mle over parameter space fixes the model distributions 
to the corresponding \iproj of $\refP$ in ambient and inner  \totemplex, $\iprojSymbol\refP$ and $\iprojSymbol[{\equivclass[{\vecOp[M']}][\prob p]}]\refP$ respectively, the likelihood ratio ($LR$) could be defined via
\equ{
\log LR \equiv \log \frac{\multidistro(N\prob q;\iprojSymbol\refP)}{\multidistro(N\prob q; \iprojSymbol[{\equivclass[{\vecOp[M']}][\prob p]}]\refP)} = 
- N\left[\infbraket{\iprojSymbol\refP}{\prob q}-\infbraket{\iprojSymbol[{\equivclass[{\vecOp[M']}][\prob p]}]\refP}{\prob q} \right]
+ \order{1}
~.
}
In the model-centric dual picture, 
this $N$-leading estimate gives the ratio of the likelihood to sample e.g.\ the data $\texttt X$ summarised by $\prob q=\prob f$ under $\mathcal M$ over the likelihood to sample ${\prob f}$ under $\mathcal M'$.

Applying Lemma~\ref{Lemma:Pythagoras} twice for  ${\prob q}\in\equivclass[{\vecOp[M']}][\prob p]$ and ${\prob q}\in\equivclass$, we see that
\equ{
\log LR = 
- N\left[\infbraket{\refP}{\iprojSymbol[{\equivclass[{\vecOp[M']}][\prob p]}]\refP} -\infbraket{\refP}{\iprojSymbol\refP} \right]
+ \order{1}
~.
}
In turn, applying again Lemma~\ref{Lemma:Pythagoras} for $\iprojSymbol[{\equivclass[{\vecOp[M']}][\prob p]}]\refP\in\equivclass$ we arrive at 
\equ{
\log LR = 
- N\infdiv{\iprojSymbol[{\equivclass[{\vecOp[M']}][\prob p]}]\refP}{\iprojSymbol\refP}
+ \order{1}
~.
}
This expression for the $LR$ confirms the form of the density appearing in 
\eqref{eq:Quotient:KnowledgeDistribution}, when interpreted as a likelihood ratio~\cite{10.1214/aoms/1177732360} among the nested \totem
constructed by $\vecOp'$ that implies the ambient \totem constructed by $\vecOp$. 
To rigorously derive the distribution of the $LR$ --\,analogously to Theorem~\ref{th:LRT}\,-- 
in the model-centric formulation requires an extended discussion, see e.g.~\cite{doi:10.1080/00031305.1982.10482817,vaart_1998}.

\section{Gaussian integrals}
\label{app:Gaussian}

For any real symmetric positive-definite matrix $\mathbf A$ the associated quadratic form 
\equ{
\label{eq:app:QuadraticForm}
Q_{\mathbf A}(\mathbf x) = \mathbf x^T \mathbf A \mathbf x = \sum_{i,j=1}^n x_iA_{i,j}x_j
}
remains strictly  positive for any $\mathbf x\in\mathbb R^{n}\setminus\lbrace0\rbrace$.
The core formula we shall need is the multivariate Gaussian integral with source term $\mathbf j\in\mathbb R^{n}$
\equ{
\label{eq:GaussianIntegrals:Sourced}
\int \dd^n \mathbf x \exp \left\{ -\sfrac12Q_{\mathbf A}(\mathbf x) - \mathbf x^T \mathbf j  \right\} 
=
\sqrt{\frac{(2\pi)^n}{\det \mathbf A}} \exp\left\{ \sfrac12 \mathbf j^T \mathbf A^{-1} \mathbf j \right\} ~.
}
In the context of multivariate statistics, $\mathbf A$ is called the \textit{precision} matrix and its inverse the covariance of the associated normal distribution
\equ{
\text{norm}\left(\mathbf x ; \mathbf A^{-1}\mathbf j , \mathbf A^{-1}\right) = \sqrt{\frac{1}{\det\left(2\pi\mathbf A^{-1}\right)}}
\exp\left\{ - \sfrac12 \left(\mathbf x - \mathbf A^{-1}\mathbf j  \right)^T  \mathbf A \left(\mathbf x - \mathbf A^{-1}\mathbf j  \right) \right\}
}
with mean $\langle\mathbf x\rangle=\mathbf A^{-1}\mathbf j$.

An often encountered expectation over a Gaussian function is of the quadratic form
\equ{
\label{eq:appGaussian:QuadraticExpectation}
\sqrt{\frac{1}{\det\left(2\pi\mathbf A^{-1}\right)}}\int \dd^n\mathbf x\, e^{-\sfrac12 Q_{\mathbf A}(\mathbf x)} \, \sfrac12 Q_{\mathbf A}(\mathbf x) = \frac{n}{2}
~.
}

\subsection{Marginals of Gaussian distribution}
\label{apps:MarginalGaussian}

To marginalize over some of the directions in $\mathbf x$ only,  
we need to partition the ambient space into subspaces
\equ{
\mathbf x = \begin{pmatrix}
\mathbf x_1\\
\mathbf x_2
\end{pmatrix}
~,
}
where $\mathbf x_a\in\mathbb R^{d_a}$  for $a=1,2$. 
This breaks the precision matrix into  so-called conformable partitions 
\equ{
\mathbf A =
\begin{pmatrix}
\mathbf A_1 &  \mathbf B \\
\mathbf B & \mathbf A_2
\end{pmatrix}
~,
}
where $\mathbf B$ is a $d_1\times d_2$ real matrix. $\mathbf A_a$ are taken to be symmetric positive-definite matrices and hence invertible. 
Accordingly, the quadratic form \eqref{eq:app:QuadraticForm} gives 
\begin{align}
\label{eq:app:QuadraticForm_Decomposition}
Q_{\mathbf A}(\mathbf x) = \mathbf x^T_1 \mathbf A_1 \mathbf x_1 + \mathbf x^T_2 \mathbf A_2 \mathbf x_2 + 2 \mathbf x^T_1 \mathbf B \mathbf x_2
\end{align}
%
For clarity, we set the source term  to zero, $\mathbf j=\mathbf 0$, which can be achieved by a shift of the coordinate system. 

Marginalizing over the Gaussian integral in say $\mathbf x_1$ gives 
\equ{
\label{eq:app:Marginal_over_x1}
\int \dd^{d_1}\mathbf x_1 \text{norm}\left(\mathbf x ; \mathbf 0 , \mathbf A^{-1}\right) =
\sqrt{\frac{\det\left(2\pi\mathbf A_1^{-1}\right)}{\det\left(2\pi\mathbf A^{-1}\right)}}
\exp\left\{ - \sfrac12 \mathbf x_2^T\left(\mathbf A_2 - \mathbf B^T\mathbf A^{-1}_1 \mathbf B \right)\mathbf x_2  \right\} 
}
which is again a $d_2$-dimensional Gaussian integral in $\mathbf x_2$.
That the resulting precision matrix is positive definite can be easily seen by completing the square for $\mathbf x_1$ in \eqref{eq:app:QuadraticForm_Decomposition}:
\equ{
Q_{\mathbf A}(\mathbf x) = \left(\mathbf x_1 + \mathbf A_1^{-1}\mathbf B\mathbf x_2\right)^T\mathbf A_1 \left(\mathbf x_1 + \mathbf A_1^{-1}\mathbf B\mathbf x_2\right)
+ 
\mathbf x_2^T \left(\mathbf A_2 - \mathbf B^T\mathbf A_1^{-1}\mathbf B \right)\mathbf x_2
~.
}
Per assumption, the l.h.s.\ always stays  positive for any $\mathbf x\in\mathbb R^{n}\setminus\lbrace0\rbrace$, thus for any vectors $\mathbf x_a\in\mathbb R^{d_a}\setminus\lbrace0\rbrace$.
Choose in particular $\mathbf x_1 = -\mathbf A_1^{-1}\mathbf B\mathbf x_2$, then 
\equ{
\forall~ \mathbf x_2 \in \mathbb R^{d_2} \setminus\lbrace0\rbrace\,:\quad
0 < Q_{\mathbf A} = \mathbf x_2^T \left(\mathbf A_2 - \mathbf B^T\mathbf A_1^{-1}\mathbf B \right)\mathbf x_2
}
Hence, $\mathbf A/\mathbf A_1 \equiv \mathbf A_2 - \mathbf B^T\mathbf A_1^{-1}\mathbf B$ has to be positive definite. 
Incidentally, $\mathbf A/\mathbf A_1$ is called the \textit{Schur} complement of block $\mathbf A_1$ for which 
%
the relation  
\equ{
\det\mathbf A = \det\mathbf A_1 \cdot \det \mathbf A/\mathbf A_1 
}
holds by the rules of determinants.
This means that the normalization of the marginal integral in \eqref{eq:app:Marginal_over_x1} reduces to 
\equ{
\label{eq:Schur:det_relation}
\frac{\det\left(2\pi\mathbf A_1^{-1}\right)}{\det\left(2\pi\mathbf A^{-1}\right)} = \left[{\det\left(2\pi \left(\mathbf A/\mathbf A_1\right)^{-1}  \right)}\right]^{-1}
}
using that $n=d_1+d_2$ for complementary spaces in $\mathbf R^n$.
Consequently, 
\equ{
\int \dd^{d_1}\mathbf x_1 \text{norm}\left(\mathbf x ; \mathbf 0 , \mathbf A^{-1}\right) =
\text{norm}\left(\mathbf x_2; \mathbf 0, \left(\mathbf A/\mathbf A_1\right)^{-1}  \right) 
~,
}
asserting the closure of normal distribution under marginalization.

\section{Constrained integrals over multinomial expansions}
\protect\label{app:integration_largeN}

In support of continuum limit \eqref{eq:FromDiscreteToContinuous}, we verify the consistency of large-$N$ expansion in a \totemplex $\equivclass[{\op N}]=\mathcal P$. 
For simplicity, we take $\statespace_0=\emptyset$.

Following the strategy of Section~\ref{ssc:DistroOverDistros},
we first express probabilities on the simplex via fluctuations  around  the \iproj which in absence of phenomenological constraints coincides with $\refP\in\mathcal P$. 
Starting from expansion \eqref{eq:MutlinomialDistro:LargeN}, we assert that   
\small{
\begin{eqnarray*}
\sum_{N\prob p\in\mathcal P} \multidistro(N\prob p;\refP) =
\\
(2\pi N)^{-\frac{\abss{\statespace}-1}{2}} \int\prod_{\config\in\statespace}\dd(N\prob p_\config) \exp \left\{ - N \sum_{\config\in\statespace} \prob p_\config\log\frac{\prob p_\config}{\refP_\config} - \sfrac12 \sum_{\config\in\statespace}\log\prob p_\config + \order{N^{-1}} \right\} \delta(N\vev{\prob p} - N)
=
\\[1ex]
(2\pi N)^{-\frac{\abss{\statespace}-1}{2}} \int\prod_{\config\in\statespace}\dd(N\probflac_\config) \exp \left\{ - N \sum_{\config\in\statespace}\frac{\probflac_\config^2}{2\refP_\config} - \sfrac12 \sum_{\config\in\statespace}\log\prob p_\config + \order{N^{-1/2}} \right\} 
\int \frac{\dd k}{2\pi} \exp\left\{i N k \sum_{\config\in\statespace}\probflac_\config\right\}
\\[1ex]
=
\sqrt{\frac{N}{2\pi}} \int \dd k \exp\left\{-\frac N2 k^2\right\}
=
1
\quads{\Rightarrow}\checkmark
\end{eqnarray*}
}\normalsize
For Dirac's delta function in the first line, we have used integral representation~\cite{NIST:DLMF}. 
After decomposing into fluctuations in the kernel of $\op N$,
\begin{equation*}
    \sum_{\config\in\statespace}\probflac_\config = 0 ~,
\end{equation*}
\eqref{eq:GaussianIntegrals:Sourced} was used to resolve the sourced Gaussian integral.
In the last line, we have exploited that $\vev{\refP}=1$ for a valid reference distribution. 
Latter two equalities are understood up to $1/N$-suppressed quantities.
Imposing next $D-1$ phenomenological constraints, 
we can analogously compute for each \totemplex constructed  by an irreducible element $\vecOp$ 
a probability mass
\begin{align}
 \label{eq:TOTEMplex:ProbabilityMass}
    \mathfrak P = &\,\,
    \log \int_{\prob p\in\equivclass} \multidistro(N\prob p; \refP) 
    \\[1ex]
    =&\,\,
    -N \infbraket{\refP}{\iprojSymbol\refP} - \frac{D-1}{2}\log N 
    - \sfrac12 \log\det \vev{\vecOp \iprojSymbol\refP \vecOp[{M^T}]}  + \order{N^{-1/2}}
        ~.
\nonumber
\end{align}
The $\order{1}$ corrections depend on the eigenvalue problem \eqref{eq:Pheno:Measurement}.

Theorem~\ref{th:Universality:GaussianApproximation} presents the $N$-leading estimate for the distribution over distributions on the \totemplex.
The positive-definiteness of kernel \eqref{eq:GaussianKernel} has already been established. 
Our aim is then to explicitly verify that the  multinomial distribution restricted to the \totemplex and normalized by the total probability mass $\mathfrak P$ reduces to the multivariate Gaussian distribution given in \eqref{eq:KnowledgeFluctuations:UniversalGaussian} to tractable order in the large-$N$ expansion.

Combining  Lemma~\ref{lm:UniversalKernel} substituted into expansion \eqref{eq:MutlinomialDistro:LargeN} together with \eqref{eq:TOTEMplex:ProbabilityMass} we indeed have 
\equ{
\label{eq:DistroOVerDistros:TOTEMplex_ProbabilityMassDensity}
\frac{\multidistro(N\prob p;\refP)}{\mathfrak P} = 
\sqrt{\frac{\det\left(2\pi N \vev{\vecOp \prob q \vecOp[{M^T}]}\right)}{\det\left(2\pi N \vev{\Id \prob q \Id}\right)}}
\exp\left\{-\frac N2 \mathcal K(\probflac) \right\}
\quads{\text{with}}
\prob q\equiv \iprojSymbol\refP
~.
}
Readily the 
factors of $(2\pi N)^{D/2-\abss{\statespace}/2} = (2\pi N)^{-\abss{\ker\vecOp}/2}$ from the determinant ratio combined with
$N^{\abss{\ker\vecOp}}$ from the integral measure (viz.\ \eqref{eq:FromDiscreteToContinuous}) produce the appropriate data-independent normalization for $\abss{\ker\vecOp}$-dimensional Gaussian kernel that scales linearly in $N$.
%
Hence, it only remains to show that the data-dependent ratio of determinants properly reduces to the normalizing factor of a multivariate  Gaussian distribution with $\order{1}$ precision $\vev{\op N_a\prob q^{-1}\op N_b}$ from \eqref{eq:GaussianKernel}. 

In fact, relation 
\equ{
\label{eq:constructingElements:det_Relation}
\det \left(\vev{\vecOp[N]\prob q^{-1}\vecOp[N^T]}\right) = \frac{\det \left(\vev{\vecOp\prob q\vecOp[M^T]}\right)}{\det\vev{\Id\prob q\Id}}~,
}
which holds along the lines of \eqref{eq:Schur:det_relation} between the determinants in the appropriate subspaces,  ensures the proper  normalization of data-dependent part in \eqref{eq:DistroOVerDistros:TOTEMplex_ProbabilityMassDensity}.
Based on \eqref{eq:VectorOperator:OrthogonalComplement},
we can arrange for 
\equ{
\vev{\op M_\alpha \op M_\beta} = \delta_{\alpha\beta}
\quads{,}
\vev{\op N_a \op N_b} = \delta_{ab}
\qand
\vev{\op M_\alpha \op N_a} = 0
~,
}
so that 
\equ{
\mathbf U 
= 
\begin{pmatrix}
\vecOp \\
\vecOp[N]
\end{pmatrix}
~.
}
defines an $\abss{\statespace}$-dimensional vector operator.
In turn, this enables us to introduce two block matrices 
\equ{
\label{eq:TOTEMplex:blockmatrices}
\mathbf A = \vev{\mathbf U \prob q \mathbf U^T} = 
\begin{pmatrix}
\vev{\vecOp \prob q \vecOp[M^T]} & \vev{\vecOp \prob q \vecOp[N^T]} \\
\vev{\vecOp[N] \prob q \vecOp[M^T]} & \vev{\vecOp[N] \prob q \vecOp[N^T]}
\end{pmatrix}
\qand
\mathbf B = \vev{\mathbf U \prob q^{-1} \mathbf U^T}
~.
}
Evidently, $\mathbf B$ is the inverse of $\mathbf A$ so that we can write 
\equ{
\begin{pmatrix}
\Id_{\rank\mathbf C} & \mathbf 0_{\rank\mathbf C\times\abss{\ker\mathbf C}} \\
\mathbf B_{21} &  \mathbf B_{22}
\end{pmatrix}
\begin{pmatrix}
\mathbf A_{11} & \mathbf A_{12} \\
\mathbf A_{21} & \mathbf A_{22} 
\end{pmatrix}
=
\begin{pmatrix}
\mathbf A_{11} & \mathbf A_{12} \\
\mathbf 0_{\abss{\ker\mathbf C}\times\rank\mathbf C} & \Id_{\abss{\ker\mathbf C}}
\end{pmatrix}
~,
}
where $\mathbf A_{\star\star}$ and $\mathbf B_{\star\star}$ compactly refer to  the corresponding block-matrices from \eqref{eq:TOTEMplex:blockmatrices}.
Taking the determinant on both sides of latter equality, 
\equ{
\det \mathbf B_{22} \det\mathbf A = \det \mathbf A_{11}~,
}
we arrive noting that $\det\mathbf A = \det\left(\Id\prob q\Id\right)$ at \eqref{eq:constructingElements:det_Relation}.

\end{appendix}


\bibliographystyle{paper} 
\bibliography{manuskript}       

\end{document}